\documentclass[10pt]{amsart}
\usepackage{latexsym}
\usepackage[all]{xy}
\usepackage{color}
\usepackage{srcltx}
\usepackage[normalem]{ulem}
\usepackage{amsmath}
\usepackage{amssymb}
\usepackage{amsfonts}

\usepackage[
	hypertexnames=false,
	hyperindex,
	pdftex,
        colorlinks,
	linkcolor=black,
	citecolor=blue,
	urlcolor=red,
]{hyperref}

\mathchardef\mhyphen="2D

\newtheorem{theorem}{Theorem}[section]

\theoremstyle{definition}
\newtheorem{definition}[theorem]{Definition}
\newtheorem{example}[theorem]{Example}
\newtheorem{remark}[theorem]{Remark}

\def\and{\ and }

\DeclareMathOperator{\Hom}{Hom}				
\DeclareMathOperator{\Ker}{Ker}

\newcommand{\infCat}{\hbox{$\sz\mhyphen\mathit{Cat}$}}
\newcommand{\bfinfCat}{\hbox{$\sz\mhyphen\mathbf{Cat}$}}
\newcommand{\bfinfgpd}{\hbox{$\sz\mhyphen\mathbf{Gpd}$}}
\newcommand{\Loc}{\mathop{\mathit{Loc}}}
\newcommand{\uLocs}{\underline{\mathit{Loc}}^*}
\newcommand{\erase}[1]{}
\newcommand{\Kvir}{^\textit{K-vir}}
\newcommand  {\sz}   {\infty}
\newtheorem*{defi}{Definition}

\newcommand  {\ecdga}     {\epsilon\mhyphen\mathbf{cdga}}

\newcommand  {\perf}   {\textit{perf}}
\newcommand  {\qcoh}   {\textit{qcoh}}
\newcommand  {\HH}   {\mathit{H\! H}}

\newtheorem{thm}{Theorem}[section]
\newtheorem{prop}[thm]{Proposition}

\newtheorem{df}[thm]{Definition}
\newtheorem{cor}[thm]{Corollary}
\newtheorem{conj}[thm]{Conjecture}
\newtheorem{war}[thm]{Warning}
\newtheorem{quest}[thm]{Question}

\newcommand {\fun} {Fun^{\sz}}
\newcommand {\Map} {\mathbb{R}\mathbf{Map}}
\newcommand {\Parf} {\mathbb{R}\mathbf{Perf}}

\newcommand {\OO} {\mathcal{O}}

\newcommand {\Spec} {\mathbf{Spec}}
\newcommand  {\dg}     {\mathbf{dg}}

\newcommand {\scat} {\infCat}   %%%{\infty-Cat}
\newcommand  {\uscat} {\bfinfCat} %%%{\infty-\mathbf{Cat}}

\newcommand  {\dAff}     {\mathbf{dAff}}
\newcommand  {\ncdga}     {\mathbf{cdga}^{\leq 0}}

\newcommand  {\length}   {length}

\newcommand  {\dSt}   {\mathbf{dSt}}
\newcommand  {\dSch}   {\mathbf{dSch}}
\newcommand  {\Sch}   {\mathbf{Sch}}
\newcommand  {\Top}   {\mathbb{S}}
\newcommand  {\scomm}     {\mathbf{sComm}}
\newcommand  {\spcomm}     {\mathbf{SpComm}}

\newcommand  {\Pol}   {\mathcal{P}ol}
\DeclareMathOperator*{\Zaff} {\operatorname{\mathit{Zaff}}}

\title[Derived Algebraic Geometry]{Derived Algebraic Geometry}

\author[Bertrand To\"en]{Bertrand To\"en}  \thanks{%CNRS
{Bertrand To\"en, Universit\'e de Montpellier 2, Case courrier 051,
  B\^at 9, Place Eug\`ene Bataillon, Montpellier Cedex 5, France. bertrand.toen@um2.fr }}

\begin{document}

\maketitle

\begin{abstract}
  This text is a survey of derived algebraic geometry. It covers a
  variety of general notions and results from the subject with a view
  on the recent developments at the interface with deformation
  quantization.
\end{abstract}

\iffalse

\usepackage{amsmath}
\usepackage{amssymb}
\usepackage{amsfonts}
\usepackage{theorem}
\usepackage{ifthen}

\fi

\tableofcontents

\section*{Introduction}

\emph{Derived algebraic geometry} is an extension of algebraic
geometry whose main purpose is to propose a setting to treat
geometrically \emph{special situations} (typically \emph{bad}
intersections, quotients by \emph{bad} actions,\dots), as opposed to
generic situations (transversal intersections, quotients by free and
proper actions,\dots). In order to present the main flavor of the
subject we will start this introduction by focussing on an
  emblematic situation in the context of algebraic geometry, or in
geometry in general: basic intersection theory. The setting is here a
smooth ambient algebraic variety $X$ (e.g., $X=\mathbb{C}^{n}$) and
two algebraic smooth subvarieties $Y$ and $Z$ in $X$, defined by some
system of algebraic equations. The problem is to understand, in a
refined and meaningful manner, the intersection $Y\cap Z$ of the two
subvarieties inside $X$. The nice, generic situation happens when they
meet transversally in $X$, i.e., when the tangent spaces of $Y$ and
$Z$ generate the whole tangent space of $X$. In this case their
intersection is itself a subvariety which possesses all the expected
properties. For example, its codimension is the sum of the
codimensions of the two subvarieties. Pathologies appear precisely
when the intersection ceases to be transversal, so that $Y$ and $Z$
may meet with higher order multiplicities or their intersection may
have some the components that are not of the expected dimension.  In
such cases, the naive intersection in $X$ fails to represent the
correct intersection.

The geometric treatment of these special situations is classically
based on cohomological methods, for which the \emph{correct
  intersection} is obtained as a cohomology class on $X$ (e.g., an
element in de~Rham cohomology, in complex cobordism, in algebraic
$K$-theory or in the intersection ring, possibly with a support
condition). In this approach the two varieties $Y$ and $Z$ must first
be slightly deformed in $X$, in order to obtain a generic situation
for which the intersection of the deformed subvarieties becomes nice.
This cohomological approach has shown itself to be extremely powerful,
particularly when concerned with questions of a numerical nature (as
is typical in enumerative geometry). However, its main drawback is
that the intersection $Y\cap Z$ is exhibited as a cohomology class,
thus depriving it of any geometric content.

Derived algebraic geometry offers a setting in which the intersection $Y\cap Z$ is realized as 
a \emph{derived scheme},, an object that encompasses the cohomological
and numerical aspects of the intersection, but at the same time
remains of a geometric nature. This \emph{derived intersection} is
obtained by a certain \emph{homotopical perturbation} of the naive
intersection $Y\cap Z$, which now comes equipped with an additional
structure, called the \emph{derived structure}, that reflects the
different possible pathologies such as the existence of multiplicities
and  defect of expected dimension. Intuitively, the derived intersection
consists of triples $(y,z,\alpha)$, where $y$ is a point in $Y$, $z$ a
point in $Z$, and $\alpha$ is a \emph{infinitesimally small continuous
  path on $X$ going from $y$ to $z$}.  The third component, the path
$\alpha$, is here the new feature, and is responsible for the derived
structure existing on the naive intersection $Y\cap Z$, which itself
sits inside the derived intersection as the locus where $\alpha$ is
constant. It is however difficult to provide a precise mathematical
meaning of the expression \emph{infinitesimally small continuous path
  on $X$}, without having to delve deep into the technical details of
the definition of a derived scheme (see definitions \ref{d1} and
\ref{d1'}). Indeed, the path $\alpha$ above is of \emph{higher
  categorical} nature, and consists of a homotopy (or equivalently a
$2$-morphism) in the \emph{$\sz$-category of derived schemes} that
will be introduced in our paragraph $\S2.2$ after we have reviewed the
basics of $\sz$-category theory.  We therefore kindly ask to the
reader to use his or her imagination, and to believe that these
concepts can be mathematically incarnated in a meaningful manner which
will be discussed later on in the main body of this paper.  It turns
out that this point of view on the specific problem of intersecting
two subvarieties is very general, and can also be applied to deal with
other instances of special situations encountered in algebraic
geometry, or in geometry in general.

Another important example we would like to mention is the problem of
considering the quotient $X/G$ of a nonfree action of an algebraic
group $G$ on a variety $X$.  In the same way that we have perturbed
the naive intersection by introducing a small path $\alpha$ as
explained above, we can introduce a refined quotient by requiring that
for $x\in X$ and $g\in G$, the points $x$ and $gx$ are not necessarily
\emph{equal} in $X/G$, but are \emph{homotopic}, or equivalently
linked by a path. Taking the quotient of $X$ by $G$ in this new sense
consists of formally adding a path between each pair of points
$(x,gx)$. The result is a well-known and already well-identified
object: the \emph{quotient stack} of $X$ by $G$ (see
\cite[2.4.2]{lm}), or equivalently the \emph{quotient groupoid} (see
\cite[2.4.3]{lm}).

This important idea to replace an intersection $Y\cap Z$, which
amounts of points $(y,z)$ with $y=z$, by the triples $(y,z,\alpha)$ as
above, or to treat a quotient $X/G$ as explained before, is not new
and belongs to a much more general stream of ideas often referred as
\emph{homotopical mathematics} (whose foundations could be said
  to be the homotopy theory of type developed by Voevodsky and al.,
see \cite{voev}). Very briefly, the expression
\emph{homotopical mathematics} reflects a shift of paradigm in which
the relation of equality relation is weakened to that of
homotopy\footnote{It is is very similar to the shift of paradigm that
  has appeared with the introduction of category theory, for which
  \emph{being equal} has been replaced by \emph{being naturally
    isomorphic.}}. Derived algebraic geometry is precisely what
happens to algebraic geometry when it is considered from the point of
view of homotopical mathematics. The purpose of this survey is on the
one hand to explain how these ideas have be realized mathematically,
and on the other to try to convince the reader that derived algebraic
geometry brings a new and interesting point of view on several aspects
and questions of algebraic geometry.\smallbreak

The mathematical foundations of derived algebraic geometry are
relatively recent. They date mostly from the first decade of this
century and appear in a series of works: \cite{hagdag}, \cite{hagI},
\cite{hagII}, \cite{luri3}, \cite{seat}, \cite{dag}.  It is built on
seminal work and important ideas in algebraic geometry, algebraic
topology and mathematical physics, some of which goes back to the
early 1950's (see $\S1$ for a selection of these pieces of
history). Derived algebraic geometry developed very
fast during the last decade, due to the works of various authors. Today
the subject possesses a very solid foundation and has a rather
large spectrum of interactions with other mathematical domains,
ranging from moduli spaces in algebraic geometry to aspects of
arithmetic geometry and number theory, not to mention geometric
representation theory and mathematical physics. Let us single out
recent progress in some of these areas, made possible by derived
algebraic geometry.

\begin{enumerate}

\item \textbf{Geometric Langlands.} The geometric version of the
  Langlands correspondence, as introduced by Beilinson and Drinfeld
  (see \cite{bd}), predicts the existence of an equivalence between
  two derived categories attached to a given smooth and proper complex
  curve $C$ and a reductive group $G$. On the one hand, we have the
  moduli space (it is really a stack, see \cite{lm}) $Bun_{G}(C)$ of
  principal $G$-bundles on $C$, and the corresponding derived category
  of $D$-modules $D(\mathcal{D}_{Bun_{G}(X)})$. On the other hand, if
  $G^{\vee}$ denotes the Langlands dual of $G$, we have the moduli
  space (again a stack) $\Loc_{G^{\vee}}(C)$ of principal
  $G^{\vee}$-bundles on $C$ endowed with flat connections, as well as
  its quasi-coherent derived category
  $D_{\qcoh}(\Loc_{G^{\vee}}(C))$. The geometric Langlands
  correspondence essentially predicts the existence of an equivalence
  of categories between $D_{\qcoh}(\Loc_{G^{\vee}}(C))$ and
  $D(\mathcal{D}_{Bun_{G}(X)})$ (see \cite{gait2}). Derived algebraic
  geometry interacts with the subject at various places. First
  of all, the moduli space $\Loc_{G^{\vee}}(C)$ naturally comes
  equipped with a nontrivial derived structure, concentrated around
  the points corresponding to flat $G^{\vee}$-bundles with many
  automorphisms. This derived structure must be taken into account for
  the expected equivalence to exist, as it modifies nontrivially the
  derived category $D_{\qcoh}(\Loc_{G^{\vee}}(C))$ (see for instance
  \cite[\S4.4(5)]{seat}). Moreover, the statement that the two above
  derived categories are equivalent is only a rough approximation of
  the correct version of the geometric Langlands correspondence, in
  which the notion of singular supports of bounded coherent complexes
  on $\Loc_{G^{\vee}}(C)$ must be introduced, a notion that is of a
  derived nature (see \cite{gaitarin}).

\item \textbf{Topological modular forms.} The notion of topological
  modular forms lies at the interface between stable homotopy theory,
  algebraic geometry and number theory (see \cite{hopk}). It has been
  observed that the formal deformation space of an elliptic curve $E$
  gives rise to a generalized cohomology theory $ell_E$, or
  equivalently a ring spectrum, called elliptic cohomology associated
  to $E$.  The spectrum of topological modular forms itself appears by
  an integration over the whole moduli space of elliptic curves of the
  spectra $ell_E$. The integration process has been considered as a
  very technical question for a long time, and has first been solved
  by deformation theory (see \cite{goer} for more about the
  subject). More recently, a completely new approach, based on derived
  algebraic geometry (or more precisely its topological analogue,
  \emph{spectral geometry}, see $\S3.4$) has been proposed in
  \cite{luri5}, in which the various spectra $ell_E$ are interpreted
  of the natural \emph{structure sheaf} of a certain spectral scheme
  (or rather stack) of elliptic curves. This approach not only 
  provided a natural and functorial point of view on elliptic
  cohomology, but also had important impact (e.g., the existence
  of equivariant version of elliptic cohomology, and later on the
  construction of the \emph{topological automorphic forms} in
  \cite{behrlaws}).

\item \textbf{Deformation quantization.}  In \cite{ptvv}, the authors
  have started developing a derived version of symplectic geometry
  motivated by the search of natural quantizations of moduli spaces
  such as Donaldson-Thomas moduli of sheaves on higher dimensional
  Calabi-Yau varieties.  This is the first step of derived Poisson
  geometry and opens up a new field of investigations related to a far
  reaching generalization of deformation quantization (see
  \cite{toen6}). This research direction will be partially presented
  in this manuscript (see \S5).  In a similar direction derived
  symplectic geometry has been used to construct and investigate
  quantum field theories (see \cite{gradgwil,cost}). In these works,
  derived algebraic geometry is essential. Many of the moduli spaces
  involved here are extremely singular (e.g., principal $G$-bundles on
  a Calabi-Yau 3-fold), and it is only when we consider them as
  derived schemes (or derived stacks) that we can notice the very rich
  geometric structures (such as a symplectic or a Poisson structure)
  that they carry.

\item \textbf{$p$-adic Hodge theory.} Finally, Bhatt (see
  \cite{bhat1,bhat2}), building on Beilinson's groundbreaking new
  proof of Fontaine's $\mathrm{C}_{\textrm{dR}}$ conjecture
  (\cite{beil}), has given strikingly short new proofs of the
  generalized Fontaine-Jannsen $\mathrm{C}_{\textrm{st}}$ and
  $\mathrm{C}_{\textrm{crys}}$, relating the algebraic de~Rham
  cohomology of algebraic varieties over $p$-adic local fields and
  their \'etale $p$-adic cohomology.  This work used in an essential
  manner the properties of the \emph{derived de~Rham cohomology},
  which computes the de~Rham cohomology in the setting of derived
  algebraic geometry (see our $\S4.4$ and $\S5.1$), and its relation
  with crystalline cohomology.

\end{enumerate}

This survey of derived algebraic geometry will, besides giving the
basic definitions and concepts of the theory, also touch on recent
developments with a particular focus on the interactions with
symplectic/Poisson geometry and deformation quantization. Our approach
is to present as far as possible \emph{mathematical facts}, without
insisting too much on formal aspects, matters of definition, or
technical issues. For instance, no proofs will be given or even
sketched. This text is therefore aimed at readers interested in having
a first look at derived algebraic geometry, but also at readers
familiar with the basics of the subject who wish to have an overview
that includes the most recent developments. In any case, the reader is
assumed to have familiarity of algebraic geometry, homological
algebra, as well as basic model category theory (as briefly recalled
in $\S2.1.1$). \smallbreak

The text is organized in 5 sections. In {Section 1}, I have
gathered some historical facts concerning the various ideas that led
to derived algebraic geometry. Its content does not pretend to be
exhaustive, and also reflects some personal taste and
interpretation. I have tried however to cover a large variety of
mathematical ideas that, I think, have influenced the modern
development of the subject. This first section is of course
independent of the sequel and can be skipped by the reader if he or
she wishes (the mathematical content truly starts in $\S2.1$), but I
have the feeling that it can explain at the same time the motivation
for derived algebraic geometry as well as some of the notions and the
definitions that will be presented in subsequent sections. In a
way, it can serve as an expanded introduction to the present paper.

{Section 2} is devoted to introducing the language of derived
algebraic geometry, namely higher category theory, and to present the
notion of \emph{derived schemes}.  The section starts with a first
paragraph on model category theory and $\sz$-category theory, by
presenting all the basic definitions and facts that will be used all
along this paper. I have tried to present the strict minimum needed
for the subject, and a priori no knowledge of higher category theory
is required for this. The second paragraph contains the first
mathematical definition of derived schemes as well as some basic
properties. More properties, such as base change, virtual classes,
tangent complexes, \ldots, are given in {Section 3}. This is
again not exhaustive and I have tried to focus on characteristic
properties of derived schemes (i.e., what makes them better behaved
than ordinary schemes). In the next two paragraphs I introduce the
functorial point of view: derived schemes are then considered as
certain $\sz$-functors from the $\sz$-category of simplicial rings. This
leads to more examples such as the derived Hilbert schemes of the
derived scheme of characters, and also leads to the notion of derived
Artin stacks, which is necessary in order to represent most of the
moduli problems appearing in derived algebraic geometry.  Finally, in
the last paragraph I have presented a short overview of derived
algebraic geometry in other contexts, such as derived analytic
geometry, spectral geometry, and the like.

The purpose of {Section 4} is to present the formal geometry of derived
schemes and derived stacks. It starts with a paragraph on cotangent complexes
and obstruction theory. The second paragraph concerns what I call \emph{formal
descent}, which is a purely derived phenomenon inspired by 
some previous work in stable homotopy theory, and which explains how 
formal completions appear by taking certain quotients by derived groupoids. 
The third paragraph presents the so-called tangent dg-lie algebra of
a derived scheme or more generally a derived stack, which is a global 
counterpart of formal geometry centered around a closed point.
The last paragraph focuses on the notion of derived loop schemes and derived
loop stacks, which are algebraic analogues of the free loop spaces studied in 
string topology. We also explain how these derived loop spaces are related
to differential forms and de~Rham theory. 

{Section 5} presents symplectic and Poisson structures in the derived
setting. It starts by a discussion of the notion of differential forms
and closed differential forms on derived schemes and on derived
stacks. In the next paragraph shifted symplectic and Lagrangian
structures are introduced, together with some basic examples coming
from classifying stacks and Lagrangian intersection theory. I have
also presented the relations with some classical notions such as
symplectic reduction and quasi-Hamiltonian actions. Paragraph 3
presents the existence results of symplectic and Lagrangian
structures, as well as some generalizations. The last paragraph of
this section contains the notion of polyvectors, Poisson structures
and their quantizations in the derived setting. This is work still in
progress and is presented here as it offers several open questions for
future research.\bigbreak

\noindent\textbf{Acknowledgements.} I am very grateful to B.~Keller
for bringing to me the idea to write a survey on derived algebraic
geometry, and for his enthusiastic support.  I would also like to
thank D.~Calaque, B.~Hennion, T.~Pantev, M.~Robalo and G.~Vezzosi for
numerous conversations that have helped me during the writing of this
paper.

%I would like to thank warmly V.~Drinfeld for having accepted to publish his letter
%to Schechtman

Finally,  I am also thankful to the referee for his or her suggestions and comments on
an earlier version of this manuscript.\bigbreak

\section{Selected pieces of history}

In this part we try to trace back a brief, and thus incomplete,
history of the mathematical ideas that have led to the modern
developments of derived algebraic geometry and more precisely to the
notion of \emph{derived Artin stacks}. Not only because we think that this might 
be of some general interest, but also because derived algebraic geometry as we
will describe later in this work is a synthesis of all these
mathematical ideas. As we will see the subject has been influenced by
ideas from various origins, such as intersection theory in algebraic
geometry, deformation theory, abstract homotopy theory, moduli and
stacks theory, stable homotopy theory, and so on. Derived algebraic
geometry incorporates all these origins, and therefore
possesses different facets and can be comprehended from different
angles. We think that knowledge of some of the key ideas  that we describe
below can help  to understand the subject from a philosophical as well
as from technical point of view.

The content of the next few pages obviously represents a personal
taste and pretends by no means to be exhaustive or very objective
(though I have tried to refer to available references as much as
possible). I apologize for any omission and misinterpretation that
this might cause. \smallbreak

\subsubsection*{The Serre intersection formula}
Serre's intersection formula \cite{serr} is commonly considered to be
the origin of derived algebraic geometry. It is probably more accurate
to consider it as the beginning of the prehistory of derived algebraic
geometry and to regard the later work of Grothendieck, Illusie and
Quillen (see below) as the starting point of the subject.

The famous formula in question expresses an intersection multiplicity
in the algebraic setting. For two irreducible algebraic subsets $Y$
and $Z$ inside a smooth algebraic variety $X$, the multiplicity
$i(X,Y.Z,W)$ is a number expressing the number of time that of $Y$ and
$Z$ meet along a fixed irreducible component $W$ of $Y\cap Z$.  For
us, the important property is that it is equal to one when $Y$ and $Z$
are smooth and meet transversely in $X$. For a nontransverse
intersection things become more complicated. Under a special condition
called \emph{Tor-independence}, the intersection number $i(X,Y.Z,W)$
can be recovered from the schematic intersection $Y\cap Z$: it is the
generic length of the structure sheaf along the component $W$
$$i(X,Y.Z,W) =  \length_{\mathcal{O}_{X,W}}(\mathcal{O}_{Y,W} \otimes_{\mathcal{O}_{X,W}}
\mathcal{O}_{Z,W}),$$ where $\mathcal{O}_{X,W}$ is the local ring of
functions on $X$ defined near $W$, similarly $\mathcal{O}_{Y,W}$
(resp. $\mathcal{O}_{Z,W}$) is the local ring of functions on $Y$
(resp. on $Z$) defined near $W$.

The Serre intersection formula explains that in general the above formula should be corrected
by higher order homological invariants  in a 
rather spectacular way:
$$i(X,Y.Z,W)=\sum_{i}(-1)^i \length_{\mathcal{O}_{X,W}}(Tor_{i}^{\mathcal{O}_{X,W}}(
\mathcal{O}_{Y,W}, \mathcal{O}_{Z,W})).$$ One possible manner to
understand this formula is that the schematic intersection of $Y$ and
$Z$ inside $X$ is not enough to understand the number $i(X,Y.Z,W)$,
and that the correcting terms
$Tor_{i}^{\mathcal{O}_{X,W}}(\mathcal{O}_{Y,W}, \mathcal{O}_{Z,W})$
should be introduced.  From the point of view of derived algebraic
geometry the presence of the correcting terms tells us that for
understanding intersection numbers, the notion of a scheme is not fine
enough. In generic situations, for instance under the assumption
that $Y$ and $W$ are smooth and meet transversely inside $X$, the
leading term equals the multiplicity $i(X,Y.Z,W)$ and the higher terms
vanish. The intersection formula is therefore particularly useful in
nongeneric situations for which the intersection of $Y$ and $Z$ has a
pathology along $W$, the \emph{worst} pathology being the case
$Y=Z$. As we shall see, the main objective of derived algebraic
geometry is to precisely understand nongeneric situations and
\emph{bad} intersections. The intersection formula of Serre is
obviously a first step in this direction: the schematic intersection
number $\length_{\mathcal{O}_{X,W}}(Tor_{0}^{\mathcal{O}_{X,W}}(
\mathcal{O}_{Y,W}, \mathcal{O}_{Z,W}))$ is corrected by the
introduction of higher terms of homological nature.

In modern terms, the formula can be interpreted by saying that
$i(X,Y.Z,W)$ can be identified with the generic length of $Y\cap Z$
(along $W$) considered as a \emph{derived scheme} as opposed as merely
a scheme, as this will be justified later on (see $\S2.2$).  However,
in the setting of the Serre intersection formula the object
$Tor_{*}^{\mathcal{O}_{X,W}}( \mathcal{O}_{Y,W}, \mathcal{O}_{Z,W})$
is simply considered as a module over $\mathcal{O}_{X,W}$, and the
whole formula is of linear nature and only involves homological
algebra. Derived algebraic geometry truly starts when
$Tor_{*}^{\mathcal{O}_{X,W}}( \mathcal{O}_{Y,W}, \mathcal{O}_{Z,W})$
is endowed with its natural multiplicative structure and is at the very
least considered as a graded algebra. In a way, the intersection
formula of Serre could be qualified as a statement belonging to
\emph{proto-derived algebraic geometry}: it contains some of the main
ideas of the subject but is not derived algebraic geometry  yet.

\subsubsection*{The cotangent complex}
In my opinion the true origin of derived algebraic geometry can be
found in the combined works of several authors, around questions
related to deformation theory of rings and schemes.  On the algebraic
side, Andr\'e and Quillen introduced a homology theory for commutative
rings, now called \emph{Andr\'e-Quillen} homology (\cite{andr,quil}),
which already had incarnations in some special cases in the work of
Harrison (\cite{harr}), and Lichtenbaum-Schlessinger
(\cite{lichschl}).  On the algebro-geometric side, Grothendieck
(\cite{grot})) and later Illusie (\cite{illu}) globalized the
definition of Andr\'e and Quillen and introduced the cotangent complex
of a morphism between schemes.  These works were motivated by the
study of the deformation theory of commutative rings and more
generally, of schemes. The leading principle is that affine smooth
schemes have a very simple deformation theory: they are rigid (do not
have nontrivial formal deformations), and their group of infinitesimal
automorphisms is determined by global algebraic vector fields.  The
deformation theory of a general scheme should then be understood by
performing an approximation by smooth affine schemes. Algebraically,
this approximation can be realized by simplicial resolving of
commutative algebras by smooth algebras, which is a multiplicative
analogue of resolving, in the sense of homological algebra, a module
by projective modules. For a commutative algebra A (say over some base
field $k$), we can choose a smooth algebra $A_0$ and a surjective
morphism $A_0 \rightarrow A$, for instance by choosing $A_{0}$ to be a
polynomial algebra. We can furthermore find another smooth algebra
$A_{1}$ and two algebra maps $A_{1} \rightrightarrows A_{0}$ in a way
that $A$ becomes the coequalizer of the above diagram of commutative
$k$-algebras.  This process can be continued further and provides a
simplicial object $A_{*}$, made out of smooth and commutative
$k$-algebras $A_{n}$, together with an augmentation $A_{*}
\longrightarrow A.$ This augmentation map is a resolution in the sense
that if we consider the total complex associated to the simplicial
object $A_{*}$, we find that the induced morphism $Tot(A_{*})
\rightarrow A$ induces isomorphisms in cohomology.  The deformation
theory of $A$ is then understood by considering the deformation theory
of the simplicial diagram of smooth algebras $A_{*}$, for which we
know that each individual algebra $A_{n}$ possesses a very simple
deformation theory.  For this, the key construction is the total
complex associated with the simplicial modules of K\"ahler
differentials
$$\mathbb{L}_{A}:=Tot(n\mapsto \Omega^{1}_{A_{n}}).$$
Up to a quasi-isomorphism this complex can be realized as a complex of
$A$-modules and is shown to be independent of the choice of the
simplicial resolution $A_{*}$ of $A$.  The object $\mathbb{L}_{A}$ is
the cotangent complex of the $k$-algebra $A$, and is shown to control
the deformation theory of $A$: there is a bijective correspondence
between infinitesimal deformations of $A$ as a commutative $k$-algebra
and $Ext^{1}_{A}(\mathbb{L}_{A},A)$. Moreover, the obstruction to
extend an infinitesimal deformation of $A$ to an order three
deformation (i.e., to pass from a family over $k[x]/x^2$ to a family
over $k[x]/x^3$) lies in $Ext^{2}(\mathbb{L}_{A},A)$.  Andr\'e and
Quillen also gave a formula for the higher extension groups
$Ext^{i}(\mathbb{L}_{A},M)$ by using the notion of derivations in the
setting of simplicial commutative algebras.

The algebraic construction of the cotangent complex has been
globalised for general schemes by Grothendieck (\cite{grot}) and
Illusie (\cite{illu}). The idea here is that the above construction
involving simplicial resolutions can be made at the sheaf level and is
then applied to the structure sheaf $\mathcal{O}_{X}$ of a scheme
$X$. To put things differently: a general scheme is approximated in
two steps, first by covering it by affine schemes and then by
resolving the commutative algebras corresponding to these affine
schemes. The important issue of how  these local
constructions are glued together is dealt with  by the use of \emph{standard
  simplicial resolutions} involving infinite dimensional polynomial
algebras. For a scheme $X$ (say over the base field $k$), the result
of the standard resolution is a sheaf of simplicial commutative
$k$-algebras $\mathcal{A}_{*}$, together with an augmentation
$\mathcal{A}_{*} \longrightarrow \mathcal{O}_{X}$ having the property
that over any open affine $U=Spec\, A \subset X$, the corresponding
simplicial algebra $\mathcal{A}_{*}(U)$ is a resolution of $A$ by
polynomial $k$-algebras (possibly with an infinite number of
generators).  Taking the total complex of K\"ahler differentials
yields a complex of $\mathcal{O}_{X}$-modules $\mathbb{L}_{X}$, called
the cotangent complex of the scheme $X$. As in the case of commutative
algebras, it is shown that $\mathbb{L}_{X}$ controls deformations of
the scheme $X$. For instance, first order deformations of $X$ are in
bijective correspondence with
$Ext^{1}(\mathbb{L}_{X},\mathcal{O}_{X})$, which is a far reaching
generalization of the Kodaira-Spencer identification of the first
order deformations of a smooth projective complex manifolds with
$H^{1}(X,T_{X})$ (see \cite{kodaspen}).  In a similar fashion the
second extension group $Ext^{2}(\mathbb{L}_{X},\mathcal{O}_{X})$
receives obstructions to extend first order deformations of $X$ to
higher order formal deformations.

I tend to consider the introduction of Andr\'e-Quillen cohomology as
well as cotangent complexes of schemes as the origin of derived
algebraic geometry. Indeed, the natural structure behind  this
construction is that of a pair $(X,\mathcal{A}_{*})$, where $X$ is the
underlying topological (Zariski) space of a scheme and
$\mathcal{A}_{*}$ is a sheaf of simplicial commutative algebras
together with an augmentation $\mathcal{A}_{*} \longrightarrow
\mathcal{O}_{X}$. Moreover, this augmentation is a resolution in the
sense that the induced morphism of complexes of sheaves
$Tot(\mathcal{A}_{*}) \rightarrow \mathcal{O}_{X}$ induces and
isomorphism on cohomology sheaves.   This makes the pair
$(X,\mathcal{A}_{*})$ a \emph{derived scheme} in the sense of
definition \ref{d1} (see also \ref{d1'}).  Here the derived scheme
$(X,\mathcal{A}_{*})$ is equivalent to a scheme, namely
$(X,\mathcal{O}_{X})$ itself, which reflects the fact that
$\mathcal{A}_{*}$ is a resolution of $\mathcal{O}_{X}$. But if we drop
the resolution condition and simply ask for an isomorphism
$H^{0}(Tot(\mathcal{A}_{*})) \simeq \mathcal{O}_{X}$, then we find the
general definition of a derived scheme. With this weaker condition,
the cohomology sheaves $H^{i}(Tot(\mathcal{A}_{*}))$ might not vanish
for $i\neq 0$:  they tell us how far the derived scheme
$(X,\mathcal{A}_{*})$ is from being equivalent to a scheme. Without 
trying to be very precise, we note here that in the context of Serre's
intersection formula above, the derived scheme obtained by
intersecting $Y$ and $Z$ in $X$ has the scheme intersection $Y\cap Z$
as its underlying scheme. The sheaf of simplicial commutative algebras
$\mathcal{A}_{*}$ will then be such that the module
$Tor_{i}^{\mathcal{O}_{X,W}}( \mathcal{O}_{Y,W}, \mathcal{O}_{Z,W}))$
is the stalk of the sheaf $H^{-i}(Tot(\mathcal{A}_{*}))$ over the
generic point of $W$.

In a way Andr\'e-Quillen went further than Grothendieck-Illusie in the
direction of what derived algebraic geometry is today, in the sense
that they did consider in an essential way simplicial commutative
algebras $A_{*}$ that might not be resolutions of algebras, and thus
can have nontrivial cohomology $H^{*}(Tot(A_{*}))$. This is a major
difference with the work of Grothendieck and Illusie in which all the
spaces endowed with a sheaf of simplicial commutative rings considered
are resolutions of actual schemes. In the context of Andr\'e-Quillen
homology the general simplicial rings appear in the important formula
(see \cite{quil}):
$$Ext^{i}(\mathbb{L}_{A},M) \simeq [A,A\oplus M[i]].$$ 
Here $A$ is a
commutative ring, $\mathbb{L}_{A}$ its cotangent complex, $M$ any
$A$-module and $A\oplus M[i]$ the simplicial algebra that is
the trivial square zero extension of $A$ by the Eilenberg-MacLane
space $K(M,i)$ ($H^{-i}(Tot(A\oplus M[i]))=M$). The bracket on the
right hand side stands for the set of maps in the homotopy category of
simplicial algebras over $A$.  This universal property of the
cotangent complex of $A$ does not appear in the works \cite{grot} and
\cite{illu}, even though the question of the interpretation of the
cotangent complex is spelled out in \cite[p.~4]{grot}.

To finish this brief review of work on cotangent complexes and
simplicial resolutions of commutative rings, we mention that there is
at least one other text in which ideas along the same line appear: a
letter from Grothendieck to Serre (\cite{grotserr}) and the manuscript
pursuing stacks (\cite{grot2}).  In \cite[p.~110]{grotserr}
Grothendieck suggests a construction of \emph{higher Jacobians} of an
algebraic variety. A first construction is proposed when the variety
is smooth, and is based on results in local cohomology. For a singular
variety Grothendieck suggests to use a simplicial resolution by smooth
algebras and to apply the construction in the smooth case degreewise.
This is of course very much in the style of the definition of the
cotangent complex that was conceived years later. Finally, in \cite[p.~554]{grot2}
(end of item~132) Grothendieck mentions the question of representing
complexes of projective $k$-modules geometrically: when the complex is
in the ``right quadrant'' (according to Grothendieck's own terms), the
answer is a linear higher stack (see $\S3.3$). However, in the ``wrong
quadrant'' case, i.e., for complex concentrated in negative
cohomological degrees Grothendieck asks the question as to what type
of structure that could represent such complexes: the answer is of
course a derived scheme \S3.3).

\subsubsection*{Derived deformation theory (DDT)}  As  mentioned
above, the introduction of cotangent complexes have been mainly
motivated by the deformation theory of algebras and schemes. The
interactions between deformation theory and derived techniques have
had a new impulse with a famous letter of Drinfeld to Schechtman
\cite{drin}. This letter is now recorded as the origin of what is
known as the \emph{derived deformation theory} (DDT for short), and
contains many of the key ideas and notions of derived algebraic
geometry at the formal level.  It has influenced a lot of work on the
subject, as for instance \cite{hini,mane}, culminating with the work
of Lurie on formal moduli problems \cite{luri}.  The main principle of
derived deformation theory stipulates that any reasonable deformation
theory problem (in characteristic zero) is associated to a
differential graded Lie algebra (dg-lie algebra for short).

A typical example illuminating this principle is the deformation
theory of a projective complex manifold $X$. The
corresponding dg-lie algebra is here $C^*(X,T_{X})$, the cochain
complex of cohomology of $X$ with coefficients in its holomorphic
tangent sheaf, turned into a dg-lie algebra using the Lie
bracket of vector fields. The space $H^{1}(X,T_{X})$ can be identified
with the first order deformation space of $X$. The space
$H^{2}(X,T_{X})$ is an obstruction space, for a given first order
deformation $\eta \in H^{1}(X,T_{X})$ the element $[\eta,\eta] \in
H^{2}(X,T_{X})$ vanishes if and only if the first order deformation of
$X$ corresponding to $\eta$ extends to a higher order deformation. More
is true, formal deformations of $X$ can all be represented by the
so-called solutions to the Mauer-Cartan equation: families of elements
$x_{i}$, for $i\geq 1$, all of degree $1$ in $C^*(X,T_{X})$ satisfying
the equation
$$d(x) + \hbox{\large$\frac12$}[x,x]=0,$$
where $x=\sum_{i}x_{i}.t^i$ is the formal power series (thus the equation
above is in $C^*(X,T_{X}) \otimes k[[t]]$).

This principle of derived deformation theory was probably already in
the air at the time of the letter \cite{drin}, as dg-lie algebras were
used already around the same time  to
describe formal completions of moduli spaces (see
\cite{schlstas,goldmill}) whose authors also refer to a letter of
Deligne. However, the precise relation between the formal deformation
theory and the dg-lie algebra was not clearly explained at
that time. For instance, various non quasi-isomorphic dg-lie algebras
could describe the same formal deformation problem (there are famous
examples with Quot schemes\footnote{If $Z \subset X$ is a closed
  immersion of smooth varieties, with sheaf of ideals $\mathcal{I}_{Z}$ the two dg-lie algebras
  $\mathbb{R}Hom_{\mathcal{O}_{Z}}(\mathcal{N}^{*}_{Z,X}[1],\mathcal{O}_{Z})$
  and
  are not quasi-isomorphic but determine the same
  functor on artinian (non-dg) rings, namely the deformation problem
  of $Z$ as a closed subscheme in $X$. The first dg-lie algebra
  considers $Z$ as a point in the Hilbert scheme of $X$ whereas the
  second dg-lie considers it as a point in the quot scheme
  $Quot(\mathcal{O}_{X})$.}).  It seems to me that one of the most
important points of the letter of Drinfeld \cite{drin} is to spell out
clearly this relation: a formal deformation problem is often defined
as a functor on (augmented) artinian algebras (see \cite{schl}), and
in order to get a canonical dg-lie algebra associated to a formal
deformation problem we need to extend this to functors defined on
\emph{artinian dg-algebras}. This is a highly nontrivial conceptual
step, which has had a very important impact on the subject.

To a dg-lie algebra $\frak{g}$ we can associate a functor
$F^{0}_{\frak{g}}$ defined on artinian algebras, by sending such an
algebra $A$ to the set of elements $x$ of degree $1$ in
$\frak{g}\otimes A$ satisfying the Mauer-Cartan equation
$d(x)+\frac{1}{2}[x,x]=0$. The main observation made in \cite{drin} is
that this functor extends naturally to a new functor $F_{\frak{g}}$
now defined on all \emph{artinian dg-algebras} (here these are the
commutative dg-algebras $A$ with finite total dimension over $k$), by
using the very same formula: $F_{\frak{g}}(A)$ consists of elements
$x$ of degree $1$ in $\frak{g}\otimes A$ such that
$d(x)+\frac{1}{2}[x,x]=0$. Moreover, Drinfeld introduces the Chevalley
complex $C^*(\frak{g})$ of the dg-lie algebra $\frak{g}$, which is by
definition the pro-artinian dg-algebra
$\widehat{Sym(\frak{g}^{*}[-1])}$ endowed with a total differential
combining the cohomological differential of $\frak{g}$ and its lie
structure. This pro-artinian dg-algebra pro-represents the functor
$F_{\frak{g}}$, and thus is thought of as the ring of formal functions
on the hypothetical formal moduli space associated to functor
$F_{\frak{g}}$.  These ideas has been formalized and developed by many
authors after Drinfeld, as for instance in \cite{hini,mane}. The
ultimate theorem subsuming these works is proven in \cite{luri} and
states that the construction $\frak{g}\mapsto F_{\frak{g}}$ can be
promoted to an equivalence between the category of dg-lie algebras up
to quasi-isomorphism, and a certain category of functors from
augmented artinian commutative dg-algebras to the category of
simplicial sets, again up to weak equivalences.

We should add a comment here, concerning the relation between the
functor $F^0_{\frak{g}}$ restricted to artinian non-dg algebras, and
the dg-lie algebra $\frak{g}$. It happens often that the functor
$F^{0}_{\frak{g}}$ is representable by a (pointed) scheme $M$, in
other words that a global moduli space $X$ exists for the moduli
problem considered (e.g.\ $\frak{g}$ can be $C^*(X,T_{X})$ for a
variety $X$ having a global moduli space of deformations $M$). By the
construction of \cite{grot,illu}, we know that $M$ is a tangent
complex $\mathbb{T}$. It is striking to notice that in general
$\mathbb{T}[-1]$ and $\frak{g}$ are not quasi-isomorphic, contrary to
what is expected. Even worse, the underlying complex $\frak{g}$ cannot
in general be the tangent complex of any scheme locally of finite
presentation: it has  usually finite dimensional total cohomology 
 and we know by \cite{avra} that this cannot happen for the
tangent complex of a scheme in general. Put differently: not
  only can we not reconstruct $\frak g$ from the functor
  $F^{0}_{\frak{g}}$, but if we try to extract a dg-lie algebra out of
  it the result will be not as nice as $\frak g$ (e.g., it will have
  infinite  dimensional total cohomology, and is probably impossible to
  describe).

What we have thus learned from derived deformation theory is that the
formal moduli space associated to a dg-lie algebra is itself a
(pro-artinian) commutative dg-algebra, and not merely a commutative
algebra. We also learned that in order to fully understand a
deformation problem, it is not enough to have a functor defined on
just artinian algebras, but that what is needed is one defined on all
artinian dg-algebras.  We are here extremely close to saying that
the main objects of study of derived deformation theory are
  commutative dg-algebras and more generally, functors on commutative
dg-algebras. This step of passing from the standard view of point of
deformation theory based on functors on artinian algebras to functors
on artinian dg-algebras is one of the most important steps in the
history of the subject: obviously the DDT has had an enormous
influence on the development of derived algebraic geometry.

\subsubsection*{Virtual classes, quasi-manifolds and dg-schemes} 
Among the most influential work concerning the global counterpart of
the derived deformation theory is \cite{kont}. The starting point is
the moduli space (or rather stack, orbifold, \ldots)
$\overline{M}_{g,n}(X,\beta)$, of stable maps $f : C \longrightarrow
X$ with $f_*[C]=\beta \in H_{2}(X,\mathbb{Z})$ fixed, where $C$ is a
curve of genus $g$ with $n$ marked points, and its relation to the
Gromov-Witten invariants of the projective manifold $X$. 
The moduli space $\overline{M}_{g,n}(X,\beta)$ is in general very
singular and trying to define the GW invariants of $X$ by performing
some integration on it would lead to the wrong answer. However,
following the DDT philosophy, the space $\overline{M}_{g,n}(X,\beta)$
can be understood locally around a given point $f : C \longrightarrow
X$ by a very explicit hypercohomology dg-lie algebra
$$\frak{g}_{f}:=C^*(C,T_{C}(-D) \rightarrow f^*(T_X)),$$ where
$T_{C}(-D)$ is the sheaf of holomorphic vector fields on $C$ vanishing
at the marked points, and the map $T_{C}(-D) \rightarrow f^*(T_X)$ is
the differential of the map $f$, which defines a two term complex of
sheaves on $C$. The dg-lie structure on $\frak{g}_{f}$ is not so
obvious to see, and is a combination of the Lie bracket of vector
fields on $C$ and the Atiyah class of the sheaf $T_{X}$. As in
\cite{drin} we can turn this dg-lie algebra $\frak{g}_{f}$ into a
pro-artinian dg-algebra by taking its Chevalley complex
$$\widehat{A}_{f}:=C^*(\frak{g}).$$ The stability of the map $f$
implies that the dg-algebra $\widehat{A}_{f}$ is cohomologically
concentrated in negative degrees and is cohomologically bounded.  The
algebra $H^{0}(\widehat{A}_{f})$ simply is the ring of formal
functions at $f\in \overline{M}_{g,n}(X,\beta)$. The higher
cohomologies $H^{i}(\widehat{A}_{f})$, which in general do not vanish
for $i< 0$, provide coherent sheaves locally defined around $f$ on the
space $\overline{M}_{g,n}(X,\beta)$.

In \cite{kont}, Kontsevich states that the local sheaves
$H^{i}(\widehat{A}_{f})$ can be glued when the point $f$ now varies
and can be defined globally as coherent sheaves $\mathcal{H}^{i}$ on
$\overline{M}_{g,n}(X,\beta)$. The family of sheaves $\mathcal{H}^{i}$
are called \emph{higher structures}, and is an incarnation of the
globalization of what the DDT gives us at each point $f \in
\overline{M}_{g,n}(X,\beta)$. Kontsevich then defines the
\emph{virtual K-theory class} of $\overline{M}_{g,n}(X,\beta)$ by the
formula
$$[\overline{M}_{g,n}(X,\beta)]\Kvir := \sum_{i\leq 0} (-1)^{i}
[\mathcal{H}^{i}] \in G_{0}(\overline{M}_{g,n}(X,\beta)).$$ In a
similar way, the complex $\frak{g}_{f}[1]$ has only nonzero cohomology
in degree $0$ and $1$, and thus defines a complex of vector bundles of
length $2$ locally around $f$. These local complexes can again be
glued to a global perfect complex of amplitude $[0,1]$, which is
called the virtual tangent sheaf of $\overline{M}_{g,n}(X,\beta)$. It
is not strictly speaking a sheaf but rather the difference of two
vector bundles and defines a class in $K$-theory
$[\mathbb{T}_{\overline{M}_{g,n}(X,\beta)}^{vir}] \in
K_{0}(\overline{M}_{g,n}(X,\beta)).$ Finally, the \emph{virtual
  homological class} of $\overline{M}_{g,n}(X,\beta)$ is defined by
the formula
$$[\overline{M}_{g,n}(X,\beta)]^{vir}\!\!:=\! 
\tau([\overline{M}_{g,n}(X,\beta)]\Kvir)\cap 
Td(\mathbb{T}_{[\overline{M}_{g,n}(X,\beta)]}^{vir})^{-1} 
\!\in H_{*}(\overline{M}_{g,n}(X,\beta),\mathbb{Q}),$$
where $Td$ is the Todd class and $\tau$ the homological Chern character (also 
called the Riemann-Roch natural transformation from the $K$-theory of coherent sheaves to 
homology, see \cite[\S18]{fult}). 

From the perspective of derived algebraic geometry the important point
is that Kontsevich not only introduced the above formula but also
provides an explanation for it based on the concept of
\emph{quasi-manifolds}. For an algebraic variety $S$ the structure of
a quasi-manifold on $S$ is a covering $\{U_{i}\}$ of $S$, together
with presentations of each $U_{i}$  as an intersection $\phi_{i}
: U_{i} \simeq Y_{i} \cap Z_{i}$ for $Y_{i}$ and $Z_{i}$ smooth
algebraic subvarieties inside a smooth ambient variety $V_{i}$.  The
precise way the various local data $Y_{i}$, $Z_{i}$, $V_{i}$ and
$\phi_{i}$ patch together is not fully described in \cite{kont}, but
it is noticed that it should involve an nontrivial notion of
equivalence, or homotopy, between different presentations of a given
algebraic variety $S$ as an intersection of smooth algebraic
varieties. These patching data, whatever they are, are certainly not
patching data in a strict sense: the local ambient smooth varieties
$V_{i}$ in which $U_{i}$ embed cannot be glued together to obtain a
smooth space $V$ in which $S$ would embed. For instance, the
dimensions of the various pieces $V_{i}$ can be nonconstant and depend
of $i$. The precise way to express these compatibilities is left
somewhat open in \cite{kont}. However, Kontsevich emphasizes that the
locally defined sheaves
$Tor_{i}^{\mathcal{O}_{Z_{i}}}(\mathcal{O}_{Y_{i}},\mathcal{O}_{Z_{i}})$,
which are coherent sheaves on $U_{i}$, glue to the globally defined
coherent sheaves $\mathcal{H}^{i}$ we mentioned before. Therefore, the
structure of a quasi-manifold on $\overline{M}_{g,n}(X,\beta)$ does
determine the higher structure sheaves $\mathcal{H}^{i}$ and thus the
K-theory virtual class $[\overline{M}_{g,n}(X,\beta)]\Kvir$.  The
virtual tangent sheaf can also be recovered from the quasi-manifold
structure, again by gluing local constructions. To each $U_{i}$, we
can consider the complex of vector bundles on $U_{i}$
$$\mathbb{T}_{i} := \left( T_{Y_{i}} \oplus T_{Z_{i}} \rightarrow
T_{V_{i}} \right).$$ Although the individual sheaves $T_{Y_{i}} \oplus
T_{Z_{i}}$ and $T_{V_{i}}$ do not glue globally on $S$ (again the
dimension of $V_{i}$ can jump when $i$ varies), the local complexes
$\mathbb{T}_{i}$ can be glued to a globally defined perfect complex
$\mathbb{T}_{S}^{vir}$, thus recovering  the virtual tangent sheaf
by means of the quasi-manifold structure. One technical aspect of this
gluing procedure is that the patching can only be done up to some
notion of equivalence (typically that of quasi-isomorphism between
complexes of sheaves) that requires a rather nontrivial formalism of
descent which is not discussed in \cite{kont}. However, the theory of
higher stacks, developed around the time by Simpson (see below),
suggests a natural way to control such gluing.

The notion of quasi-manifolds has been taken and declined by several
authors after Kontsevich. Behrend and Fantechi introduced the notion
of \emph{perfect obstruction theories} on a scheme $X$ (see
\cite{behrfant}) that consists of a perfect complex $\mathbb{T}$ on
$X$ that formally behaves as the virtual tangent sheaf of a structure
of quasi-manifold on $X$. In \cite{ciockapr} Kapranov and
Ciocan-Fontanine defined the notion of \emph{dg-schemes}, close to the
notion of supermanifolds endowed with a cohomological odd vector field
$Q$ used in mathematical physics, which by definition consists of a
scheme $X$ endowed a sheaf of commutative coherent
$\mathcal{O}_{X}$-dg-algebras. Later  a 2-categorical construction of
dg-schemes appeared in \cite{behr}. All of these notions are
approximations, more or less accurate, of the notion of a derived
scheme. They can all be used in order to construct virtual classes,
and hence suffice to define Gromov-Witten invariants in the
algebro-geometric context. However, they all suffer from a bad
functoriality behaviour and cannot be reasonably used as the
fundamental object of derived algebraic geometry (we refer to the end
of $\S3.1$ for a more  detailed discussion).

\subsubsection*{Interaction with homotopy theory, stacks and higher stacks} 
A stack is a categorical generalization of a sheaf. They are of
particular interest in moduli theory, as it is often the case that a
given moduli problem (e.g.\ the moduli of curves of a given genus)
involves objects with nontrivial automorphism groups. In such a
situation the moduli functor becomes a functor from schemes to
groupoids, and the sheaf condition on this functor is called the
descent or stack condition.

In the context of algebraic geometry stacks already appear in the late
1950s, as for instance in \cite{grot3}, as well as in
\cite{grot4}. They have been introduced to formalize the problem of
descent, but also in order to represent moduli problems with no fine
moduli spaces (it is already noted in \cite[Prop.~7.1]{grot4} that a
fine moduli space of curves does not exist). The formal definitions of
algebraic stacks appear in \cite{delimumf} in which it is shown that
the stack of stable curves of a fixed genus is an algebraic stack that
is smooth and proper over $Spec\, \mathbb{Z}$. It is interesting
however to note that many notions, such as fibered categories,
descent, quotient stack, stack of principal bundles, can be found in
\cite{grot4}. The definition of algebraic stack of \cite{delimumf} has
then been generalized by Artin in \cite{arti} in order to encompass
also moduli problems for which objects might admit algebraic groups of
automorphisms (as opposed to discrete finite groups).  In the
differential context stacks appeared even earlier in the guise
of differential groupoids for the study of foliations (see for
example \cite{heaf}, \cite{ehre}).

The insight that a notion of \emph{higher stack} exists and might be
useful goes back to \cite{grot2}. Higher stacks are higher categorical
analogues of stacks, and thus of sheaves, and presume to
classify objects for which not only nontrivial automorphisms exist,
but also automorphisms between automorphims, and automorphisms between
automorphisms of automorphisms (and so on) exist. These higher
automorphisms groups are now encoded in a higher groupoid valued
functor, together with a certain descent condition. In \cite{grot2}
Grothendieck stressed the fact that several constructions in rational
homotopy theory could be formalized by considering a nice class of
higher stacks called \emph{schematic homotopy types}, which are very
close to being a higher version of the Artin algebraic stacks of
\cite{arti}. One of the technical problem encountered in the theory of
higher stacks is the fact that it has to be based on a nice theory of
higher categories and higher groupoids that has not been fully
available until recently, and this aspect has probably delayed the
development of higher stack theory. However, in \cite{simp1} Simpson
proposed a definition of algebraic $n$-stacks based on the homotopy
theory of simplicial presheaves (due to Jardine and Joyal and largely
used in the setting of algebraic K-theory of rings and schemes), using
the principle that any good notion of $n$-groupoids should be
equivalent to the theory of homotopy $n$-types (a principle refered to
as \emph{the homotopy hypothesis} in \cite{grot2}, the higher
automorphisms groups mentioned above are then incarnated by the higher
homotopy groups of a simplicial set).  The definition of algebraic
$n$-stack of \cite{simp1} is inductive and based on a previous
observation of Walter that algebraic stacks in the sense of Artin can
be defined in any geometric context in which the notion of a smooth
morphism make sense. This simplicial approach has been extremely
fruitful, for instance in the context of higher non-abelian Hodge
theory (see \cite{simp2,simp3}), to pursue the \emph{schematization}
direction of Grothendieck's program \cite{grot2},  namely the
interpretation of rational homotopy theory and its extension over more
general bases (see \cite[Thm. 6.1]{simp1}, and also \cite{toen1}), or
to understand the descent problem inside derived and more generally
triangulated categories (see \cite{hirssimp}).

The introduction of simplicial presheaves as models for higher stacks
has had a huge impact on the subject. First of all it overcomes the
technical difficulties of the theory of $n$-groupoids and use the
power of Quillen's homotopical algebra in order to describe some of
the fundamental constructions (e.g.\ fiber products, quotients, stack
associated to a prestack, \dots). It also had the effect of bringing
the model category language in the setting of higher category theory
(see \cite{simp4}): it is interesting to note that most, if not all,
of the established theory of higher categories are based on the same
idea and rely in an essential way on model category theory
(\cite{simp5}, \cite{luri2}, \cite{rezk} to mention the most important
ones). Another aspect is that it has reinforced the interactions
between higher stack theory and abstract homotopy theory.  The
interrelations between ideas from algebraic geometry and algebraic
topology is one of the feature of derived algebraic geometry, and the
simplicial point of view of Simpson on higher stacks has contributed a
lot to the introduction of the notions of derived schemes and derived
stacks.

\subsubsection*{The influence of stable homotopy theory} Abstract homotopy theory, and 
the homotopy theory of simplicial presheaves in particular, have played an important role
in the development of higher stacks. Derived algebraic geometry has also been 
influenced by stable homotopy theory and to be more precise by the so-called
\emph{brave new algebra} (an expression introduced by Waldhausen, see \cite{vogt}). 
Brave new algebra is the study of ring spectra, also called \emph{brave new rings}, 
or equivalently of multiplicative generalized cohomology theories. It is based
on the observation that the homotopy theory of brave new rings behaves very similarly to 
the theory of rings, and that many notions and results of algebra and linear algebra
over rings extend to this setting. First of all, rings embed fully faifthfully 
into brave new rings and correspond to the discrete ring spectra. Moreover, 
the stable homotopy category possesses a nondegenerate $t$-structure whose heart is the
category of discrete spectra (abelian groups), 
so in a sense discrete rings generate the whole category of brave new rings. An efficient 
way of thinking consists of seeing the category of brave new rings as a kind of
small or infinitesimal perturbation of the category of rings, which is reflected by the
fact that the absolute brave new ring, the sphere spectrum $\mathbb{S}$ 
can be considered as an infinitesimal perturbation of the ring $\mathbb{Z}$.

One fundamental work in this direction is the work of Waldhausen on
algebraic K-theory of spaces \cite{wald}. This has been later pursued  
with the introduction of Hochschild and cyclic homology for ring
spectra, also called topological Hochschild and topological cyclic
homology, together with a Chern character map (see for instance
\cite{boks,bokshian,pirawald,vogtwald}). Another important impulse has
been the introduction of the Morava K-theories, and their
interpretations as the \emph{points} of the hypothetical object
$\Spec\, \mathbb{S}$, the spectrum (in the sense of algebraic
geometry!) of the sphere spectrum (in the sense of topology!). This
general philosophy has been spread by many authors, see for instance
\cite{mora,devihopk,hopksmit,rave}. It has also been pushed further 
with the introduction of the new idea that ring spectra should also
have an interesting Galois theory (see \cite{rogn,schwwald,mccamina}),
leading to the feeling that there should exists an \'etale topology
for ring spectra that extends the usual \'etale topology for schemes.
In a similar direction, the theory of topological modular forms (see
\cite{goer} for a survey) enhances the stack of elliptic curves with a
certain sheaf of generalized cohomology theories, and thus with a
sheaf of ring spectra, creating an even closer relation between stable
homotopy theory and algebraic geometry. It is notable that the modern
approach to the theory of topological modular forms is now based on
\emph{spectral algebraic geometry}, a topological version of derived
algebraic geometry (see \cite{luri5}).

I believe that all of these works and ideas from stable homotopy
theory have had a rather important impact on the emergence of derived
algebraic geometry, by spreading the idea that not only rings have
spectra (in the sense of algebraic geometry) but more general and
complicated objects such as ring spectra, dg-algebras and so on. Of
course, the fact that ring-like objects can be used to do geometry is
not  new here, as for instance a very general notion of relative
schemes appears already in \cite{haki}. However, the brand new idea
here was that the same general philosophy also applies to ring-like
objects of homotopical nature in a fruitful manner.

\subsubsection*{Mathematical physics} Last, but not least, derived
algebraic geometry has certainly benefited from a stream of ideas from
mathematical physics. \iffalse I am less informed of this aspect the
subject, but\fi It seems clear that some of the mathematical
structures introduced for the purpose of supersymmetry and string
theory have conveyed ideas and concepts closely related to the concept
of derived schemes.

A first instance can be found in the several 
generalizations of manifolds introduced for the purpose of supersymmetry: supermanifolds, 
$Q$-manifolds, $QP$-manifolds etc (see \cite{bere,kost,schw1,schw2}). Supermanifolds are manifolds endowed with 
the extra data of \emph{odd functions}, represented by a sheaf of $\mathbb{Z}/2$-graded 
rings. The $Q$-manifolds are essentially supermanifolds together with an 
vector field $Q$ of degree 1 (i.e. a derivation sending odd functions to even functions, 
and vice-versa), which squares to zero $Q^2=0$. The supermanifold together with the
differential $Q$ gives thus rise to a manifold endowed with a sheaf of ($\mathbb{Z}/2$-graded) 
commutative dg-algebras, which is quite close to what a derived scheme is.
One of the main differences between the theory of supermanifolds and the theory of derived schemes 
is the fact that supermanifolds were not considered up to quasi-isomorphism. In a way, 
$Q$-manifolds can be thought as \emph{strict models} for derived schemes. The influence of 
this  stream of ideas on derived algebraic geometry is not only found at the 
level of definitions, but also at the level of more advanced structures. For instance, the $QP$-manifolds of 
\cite{aksz} are certainly an avatar of the shifted symplectic structures recently introduced 
in the context of derived algebraic geometry in  \cite{ptvv}. In the same way, the BV-formalism 
of \cite{batavilk} has recent versions in the setting of derived algebraic geometry (see \cite{costgwil,vezz}). 

A second instance is related to mirror symmetry, and more particularly
to the homological mirror symmetry of Kontsevich (see
\cite{kontsoib}). Mirror symmetry between two Calabi-Yau varieties $X$
and $X^{\vee}$ is here realized as an equivalence of derived
categories $D(X) \simeq Fuk(X^{\vee})$, between on the one side the
bounded coherent derived category of $X$ (the $B$-side) and on the
other side the Fukaya category of the mirror $X^{\vee}$ (the
$A$-side). This equivalence induces an equivalence between the formal
deformation spaces of $D(X)$ and of $Fuk(X^{\vee})$, which can be
identified with the de~Rham cohomology of $X$ and the quantum
cohomology of $X^{\vee}$. Here we find again the DDT in action, as the
identification between the deformation spaces of these two categories
and the mentioned cohomologies requires to consider these moduli
spaces as formal derived moduli spaces.  This has led to the idea that
the \emph{correct} deformation space of $D(X)$ is the full de~Rham
cohomology of $X$ (and similarly for the Fukaya category), which again
convey the idea that the deformations of $D(X)$ live in a certain
\emph{derived moduli space}.

Finally, a third instance is deformation quantization. First of all, Kontsevich's proof of
the existence of deformation quantization of a Poisson manifold (see \cite{kont2}) is based on the identification
of two deformation problems (Poisson algebras and associative algebras), which is obtained
by the construction of an equivalence between the two dg-lie algebras controlling 
these deformations problems. This is once again an example of the DDT in action. The recent 
interactions between derived algebraic geometry and quantization (see our $\S5$, see also 
\cite{toen6}) also 
suggests that some of the concepts and ideas of derived algebraic geometry might have come
from  a part of quantum mathematics (my lack of knowledge of the subject prevents me
to make precise statements here).
 
\subsubsection*{Derived schemes and derived stacks} The modern
foundations of derived algebraic geometry have been set down in
a series of papers, namely \cite{hagdag}, \cite{hagI}, \cite{hagII},
\cite{luri3}, \cite{seat}, \cite{dag}. In my opinion all of the works
and ideas previously mentioned have had an enormous influence on their
authors. A particularly interesting point is that not
all of the ideas and motivations came from algebraic geometry itself,
as many important ideas also come from abstract homotopy theory and
stable homotopy theory. This probably explains many of the topological
flavours encountered in derived algebraic geometry, which adds to the
richness of the subject.

The subject has been developing fast in the last decade, thanks to the
work of many authors: B.~Antieau, D.~Arinkin, O.~Ben-Bassat,
D.~Ben-Zvi, B.~Bhatt, D.~Borisov, C.~Brav, V.~Bussi, D.~Calaque,
K.~Costello, J.~Francis, D.~Gaitsgory, D.~Gepner, G.~Ginot,
O.~Gwilliam, B.~Hennion, I.~Iwanari, D.~Joyce, P.~Lowrey, J.~Lurie,
D.~Nadler, J.~Noel, T.~Pantev, A.~Preygel, J.~Pridham, N.~Rozenblyum,
T.~Sch\"urg, M.~Spitzweck, D.~Spivak, M.~Vaqui\'e, G.~Vezzosi,
J.~Wallbridge, and others. We shall obviously not cover all
this work in the present survey, but will try to mention some of it
with emphasis on the interface with deformation quantization.

\section{The notion of a derived scheme}

In this first part we start by presenting the central object of study
of derived algebraic study: \emph{derived schemes}. The definition of
derived scheme will appear first and is rather
straightforward. However, the notions of morphisms between derived
schemes is a bit subtle and require first some notions of higher
category theory, or equivalently of homotopical algebra.  We will
start by a (very) brief overview of \emph{model category theory},
which for us will be a key tool in order to understand the notion of
\emph{$\sz$-category} presented in the second paragraph and used all
along this manuscript. We will then proceed with the definition of the
($\sz$-) category of derived schemes and provide basic examples. More
evolved examples, as well as the further notion of \emph{derived
  moduli problems} and \emph{derived algebraic stacks} are presented
in the next section. %\smallbreak

We start by extracting two principles mentioned from the variety of
ideas recalled in $\S1$, which are the foundation principles of
derived algebraic geometry. For this we begin by the following
metamathematical observation. A given mathematical theory often aims
to study a class of specific objects: algebraic varieties in algebraic
geometry, topological spaces in topology, modules over a ring in
linear algebra \dots. These objects are in general very complicated
(unless the theory might be considered as uninteresting), but it is
most often that a subclass of \emph{nice objects} naturally shows
up. As their names show the nice objects behave nicely, or at least
behave nicer than a generic object. In such a situation,
mathematicians want to believe that we fully understand the nice
objects and that a generic object should be approximated, with the
best approximation possible, by nice objects.  This metamathematical
observation can be seen in action in many concrete examples, two of
them are the following (there are myriads of other examples).

\begin{itemize}

\item (\emph{Linear algebra}) Let $A$ be a ring and we consider $A$-modules. 
The nice objects, for instance with respect to short exact sequences, 
are the projective modules. For a general $A$-module $M$ the best possible approximation 
of $M$ is a resolution of $M$ be means of projective modules
$$\dots P_{n} \rightarrow P_{n-1} \rightarrow \dots P_{0} \rightarrow M \rightarrow 0.$$

\item (\emph{Topology}) We consider topological spaces, and more particularly 
their cohomological properties. Spheres are the nice objects (for instance
from the cohomological point of view). For a space $X$ the best possible approximation
of $X$ is a cellular approximation, that is,  a CW complex $X'$ weakly equivalent
to $X$. 

\end{itemize}

These two examples possess many possible variations, for instance by replacing 
modules over a ring by objects in an abelian category, or topological 
spaces by smooth manifolds and cellular approximation by handle
body decompositions. The common denominator to all of these situations is the
behavior of the approximation construction. A delicate question is 
the uniqueness: approximations are obviously not unique in a strict sense (e.g.,  
for the two examples above, but this is a general phenomenon) and we have
to introduce a new notion of equivalence in order to control uniqueness 
and more generally functorial properties. In the examples above this notion is an 
obvious one: quasi-isomorphism of complexes in the
first example and weak homotopy equivalences in the second. As we will
see in the next paragraph, the introduction of  a new notion of equivalence automatically creates 
a higher categorical, or higher homotopical, phenomenon. This is one reason of the
ubiquity of higher categorical structures in many domains of mathematics. \smallbreak

Derived algebraic geometry is the theory derived from algebraic geometry (no pun intended)
by applying the same general principle as above, and by declaring that the good objects are 
the smooth varieties, smooth schemes and more generally smooth maps. The approximation
by smooth varieties is here the simplicial resolution of algebras by polynomial algebras
as mentioned in $\S1$. To sum up:

\begin{itemize}

\item (\emph{Principle 1 of derived algebraic geometry}) The smooth algebraic varieties,
or more generally smooth schemes and smooth maps, are good. Any nonsmooth variety, scheme
or map between schemes, must be replaced by the best possible approximation by 
smooth objects.

\item  (\emph{Principle 2 of derived algebraic geometry}) Approximations of 
varieties, schemes and maps of schemes, are expressed in terms of simplicial resolutions. 
The simplicial resolutions must only be considered up to the notion of weak equivalence, 
and are controlled by higher categorical, or homotopical, structures. 

\end{itemize}

Based on these two principles we can already extract a general definition of a derived scheme,
simply by thinking that the structure sheaf should now be a sheaf of simplicial commutative
rings rather than a genuine sheaf of commutative rings. However, principle 2 already tells
us that morphisms between these derived schemes will be a rather involved notion 
and must be defined with some care.

\begin{df}\label{d1}{(First definition of derived schemes)}
A \emph{derived scheme} consists of a pair $(X,\OO_{X})$, where $X$ is a topological space
and $\mathcal{A}_{X}$ is a sheaf of commutative simplicial rings on $X$ such that
the following two conditions are satisfied.
\begin{enumerate}

\item The ringed space $(X,\pi_{0}(\OO_{X}))$ is a scheme.

\item For all $i>0$ the homotopy sheaf $\pi_{i}(\OO_{X})$ is 
a quasi-coherent sheaf of modules on the scheme $(X,\pi_{0}(\OO_{X}))$.

\end{enumerate}
\end{df}

Some comments about this definition.

\begin{itemize}

\item A scheme $(X,\OO_{X})$ can be considered as a derived scheme in an obvious manner, 
by taking $\OO_{X}$ to be the constant simplicial sheaf of rings $\OO_{X}$. 

\item In the other way, a derived scheme $(X,\OO_{X})$ underlies a scheme
$(X,\pi_{0}(\OO_{X}))$ that is called the \emph{truncation of $(X,\OO_{X})$}.

\item A scheme can also be considered as a derived scheme $(X,\OO'_{X})$ where
now $\OO'_{X}$ is any simplicial resolution of $\OO_{X}$, that is 
$$\pi_{0}(\OO'_{X})\simeq \OO_{X} \qquad \pi_{i}(\OO'_{X})=0 \; \forall \; i>0.$$
As we will see, the derived scheme $(X,\OO'_{X})$ is equivalent to 
$(X,\OO_{X})$ (exactly as a resolution $P_{*}$ of a module $M$ is quasi-isomorphic to 
$M$ concentrated in degree $0$).

\item For a derived scheme $(X,\OO_{X})$, the truncation $(X,\pi_{0}(\OO_{X}))$
contains all the geometry. The sheaves $\pi_{i}(\OO_{X})$ on $(X,\pi_{0}(\OO_{X}))$ are pieces
reflecting the derived structure and should be thought as extraordinary nilpotent functions. 
The sheaves $\pi_{i}(\OO_{X})$ are analogous to the graded pieces 
$\mathcal{I}^{n}/\mathcal{I}^{n+1}$ where $\mathcal{I}$ is the nilradical of 
a scheme $Y$, which are sheaves on the reduced subscheme $Y_{red}$. It is  good 
intuition and also accurate  to think of $\OO_{X}$ as coming with a natural
filtration (incarnated by the Postnikov tower, see $\S2.2$) whose graded pieces
are the $\pi_{i}(\OO_{X})$. 

\end{itemize}

The above definition makes derived schemes an easy notion, at least at first sight. 
However, as already mentioned, morphisms between derived schemes require some care to be
defined in a meaningful manner. 
In the sequel of this section we will explain how to deal with derived schemes, how
to construct and define their ($\sz$)-category, but also how to work with them in practice.

\subsection{Elements of the language of $\infty$-categories}

In this paragraph we introduce the language of $\sz$-categories. The theory of $\sz$-categories 
shares very strong interrelations with the theory of model category, and most of the 
possible working definitions of $\sz$-categories available today haven been set down in the context
of model category theory. Moreover, model categories also provide a rich source of examples
of $\sz$-categories, and from the user point of view a given model category $M$ can be (should be?)
considered as a concrete model for an $\sz$-category. It is because of this strong interrelations between 
$\sz$-categories and model categories that this section starts with a brief overview of 
model category theory, before presenting some elements of the language of $\sz$-category theory. 
This language will be used in order to present the expected notion of maps between derived schemes, 
incarnated in the $\sz$-category $\dSch$ of derived schemes presented in the next paragraph.

\subsubsection*{A glimpse of model category theory}

Model category theory deals the localization problem, which consists
of inverting a certain class of maps $W$ in a given category $C$. The
precise problem is to find, and to understand, the category $W^{-1}C$
obtained out of $C$ by freely adding inverses to the maps in $W$. By
definition, the category $W^{-1}C$ comes equipped with a functor $l :
C \longrightarrow W^{-1}C$, which sends maps in $W$ to isomorphisms in
$W^{-1}C$, and which is universal with respect to this property. Up to
issues of set theory, it can be shown that $W^{-1}C$ always exists
(see \cite[\S I.1.1]{gabrzism}). However, the category $W^{-1}C$ is in
general difficult to describe in a meaningful and useful manner, and
its existence alone is often not enough from a practical point of view
(see \cite[\S2.2]{toen7} for more about the bad behavior of the
localization construction).

\begin{defi}
  A model category structure on a pair $(C,W)$ as above consists of
  extra pieces of data, involving two other class of maps, called
  \emph{fibrations} and \emph{cofibrations}, satisfying some standard
  axioms (inspired by the topological setting of topological spaces
  and weak equivalences), and ensuring that the sets of maps inside
  localized category $W^{-1}C$ possess a nice and useful description
  in terms of \emph{homotopy classes of morphisms between certain
    objects in $C$}. The typical example of such a situation occurs in
  algebraic topology: we can take $C=Top$, be the category of
  topological spaces and continuous maps, and $W$ be the class of weak
  homotopy equivalences (continuous maps inducing isomorphisms on all
  homotopy groups).  The localized category $W^{-1}C$ is equivalent to
  the category $[CW]$, whose objects are CW complexes and whose set of
  maps are homotopy classes of continuous maps. The upshot of model
  category theory is that this is not an isolated or specific example,
  and that there are many situations of different origins
  (topological, algebraic, combinatorial, \dots) in which some
  interesting localized categories can be computed in a similar
  fashion.
\end{defi}

By definition, a model category consists of a complete and cocomplete category $C$, together
with three classes of maps $W$ (called weak equivalences, or simply equivalences), 
\textit{Fib} (called fibrations), and \textit{Cof}  (called cofibrations), and satisfying the following axioms
(see \cite[\S I.1]{quil} or \cite[\S1.1]{hove} 
for more details\footnote{We use in this work the definition found in 
\cite{hove}, which differs  slightly  from the original one in \cite{quil}  (e.g., the functoriality of factorizations). 
The mathematical community seems to have adopted
the terminology of \cite{hove} as the standard one.}).

\begin{enumerate}

\item If $\xymatrix{X\ar[r]^-{f} & Y\ar[r]^-{g} & Z}$ are morphisms in $C$, then 
$f$, $g$ and $gf$ are all in $W$ if and only if two of them are in $W$. 

\item The fibrations, cofibrations and equivalences are all stable by compositions and retracts.

\item Let 
$$\xymatrix{
A \ar[r]^-{f} \ar[d]_-{i} & X\ar[d]^-{p} \\
B \ar[r]_-{g} & Y}$$
be a commutative square in $C$ with $i\in$\textit{Cof} and $p \in$\textit{Fib}. 
If either $i$ or $p$ is also in $W$ then there is
a morphism $h : B\longrightarrow X$ such that 
$ph=g$ and $hi=f$. 

\item Any morphism $f : X \longrightarrow Y$ can be factored in two ways
as $f=pi$ and $f=qj$, with $p\in $\textit{Fib},  $i\in$\textit{Cof}$\cap W$, 
$q\in$\textit{Fib}$\cap W$ and $j\in$\textit{Cof}. Moreover, the
existence of these factorizations are required to be functorial in $f$.

\end{enumerate}

The morphisms in \textit{Cof}$\cap W$ are usually called \emph{trivial cofibrations} and
the morphisms in \textit{Fib}$\cap W$ \emph{trivial fibrations}. Objects $x$ 
such that $\emptyset \longrightarrow x$ is a cofibration are called
\emph{cofibrant}. Dually, objects $y$ such that $y\longrightarrow *$ is a 
fibration are called \emph{fibrant}. The factorization axiom (4) implies that 
for any object $x$ there is a diagram
$$\xymatrix{
Qx \ar[r]^-{i} & x \ar[r]^-{p} & Rx,}$$
where $i$ is a trivial fibration, $p$ is a trivial cofibration, 
$Qx$ is a cofibrant object and $Rx$ is a fibrant object. Moreover, 
the functorial character of the factorization states that the above
diagram can be, and always will be, chosen to be functorial 
in $x$.

\subsubsection*{The homotopy category of a model category}
A model category structure is a rather simple notion, but in practice it is never
easy to check that three given classes \textit{Fib}, \textit{Cof} and $W$ satisfy the four axioms
above. This can be explained by the fact that the existence of a model category 
structure on $C$ has a very important consequence on the localized category $W^{-1}C$,
which is usually denoted by $Ho(C)$ and called the \emph{homotopy category}
in the literature\footnote{This is often misleading, as $Ho(C)$ is obtained by localization and
is \emph{not} a category obtained by dividing out the set of maps by a homotopy relation.}.
For this, we introduce the notion of homotopy between morphisms in $M$ in the following
way. Two morphisms $f,g : X \longrightarrow Y$ are called \emph{homotopic}
if there is a commutative diagram in $M$
$$\xymatrix{
X \ar[d]_-{i} \ar[dr]^-{f} & \\
C(X) \ar[r]^-{h} & Y \\
X \ar[u]^-{j} \ar[ur]_-{g} }$$
satisfying the following two properties: 
\begin{enumerate}
\item There exists
a morphism $p : C(X) \longrightarrow X$, which 
belongs to $Fib\cap W$, such that 
$pi=pj=id$.
\item The induced morphism
$$i\bigsqcup j : X\bigsqcup X \longrightarrow C(X)$$
is a cofibration.
\end{enumerate}

When $X$ is cofibrant and $Y$ is fibrant in $M$ (i.e. 
$\emptyset \longrightarrow X$ is a cofibration and
$Y \longrightarrow *$ is a fibration), it can be shown that 
being homotopic as defined above is an equivalence relation
on the set of morphisms from $X$ to $Y$. This equivalence relation
is shown to be compatible with composition, which implies the existence
of a category $C^{cf}/\sim$, whose objects are 
cofibrant and fibrant objects and whose morphisms are homotopy classes of
morphisms in $C$.

It is easy to see that if two morphisms $f$ and $g$ are homotopic in $C$ then
they are equal in $W^{-1}C$. Indeed, in the diagram above defining 
the notion of being homotopic, the image of $p$ in $Ho(C)$ is
an isomorphism. Therefore, so are the images of $i$ and $j$. Moreover, 
the inverses of the images of $i$ and $j$ in $Ho(C)$ are equal 
(because equal to the image of $p$), which implies that $i$ and
$j$ have the same image in $Ho(C)$. This implies that 
the image of $f$ and of $g$ are also equal. From this, we deduce that the localization functor
$C \longrightarrow Ho(C)$
restricted to the subcategory of cofibrant and fibrant objects $C^{cf}$ induces 
a well defined functor
$C^{cf}/\sim \longrightarrow Ho(C).$
One major statement of model category theory is that 
this last functor is an equivalence of categories. 

\begin{theorem} [{see \cite[\S I Thm.~$1'$]{quil2}, \cite[Thm.~1.2.10]{hove}}]
The above functor
$$C^{cf}/\sim \longrightarrow W^{-1}C=Ho(C).$$
is an equivalence of categories.
\end{theorem}

The above theorem is fundamental as it allows us to control, and to
describe in an efficient manner, the localized category $W^{-1}C$ in
the presence of a model structure.

\subsubsection*{Three examples}
Three major examples of model categories are the following.

\begin{itemize}

\item We set $C=Top$ be the category of topological spaces, and $W$ the class of 
weak equivalences (continuous maps inducing bijections on all homotopy groups). The class
$Fib$ is taken to be the \emph{Serre fibrations}, 
the morphism having the lifting property 
with respect to the inclusions $|\Lambda^{n}| \subset |\Delta^{n}|$, of a horn (the union 
of all but one of the codimension 1 faces) into a standard  $n$-dimensional simplex. 
The cofibrations are the retracts of the relative cell complexes. This defines a 
model category (see \cite[\S2.4]{hove}), and the theorem above states the well known fact that $Ho(Top)$
can be described as the category whose objects are CW complexes and morphisms
are homotopy classes of continuous maps.

\item For a ring $R$ we set $C(R)$ the category of (possibly
  unbounded) complexes of (left) $R$-modules. The class $W$ is taken
  to be the quasi-isomorphisms (the morphisms inducing bijective maps
  on cohomology groups). There are two standard possible choices for
  the class of fibrations and cofibrations, giving rise to two
  different model structures with the same class of equivalences,
  called the \emph{projective} and the \emph{injective} model
  structures. For the projective model structure the class of
  fibrations consists of the the epimorphisms (i.e. degreewise
  surjective maps) of complexes of $R$-modules, and the cofibrations
  are defined by orthogonality (see \cite[\S2.3]{hove}).  Dually, for
  the injective model structure the class of cofibrations consists of
  the monomorphisms (i.e. degreewise injective maps) of complexes of
  $R$-modules (see \cite[\S2.3]{hove}).  These two model categories
  share the same homotopy category $Ho(C(R))=D(R)$, which is nothing
  else than the derived category of (unbounded) complexes of
  $R$-modules. In this case, the theorem above states that the
  category $D(R)$ can also be described as the category whose objects
  are either K-injective, or K-projective, complexes, and morphisms
  are homotopy classes of maps between these complexes (see
  \cite[\S2.3]{hove}).

\item We set $C=sSet:=Fun(\Delta^{op},Set)$, 
the category of simplicial sets. For $W$ we take the class of weak equivalences of simplicial sets
(i.e., the maps inducing weak equivalences on the corresponding geometric realizations). The 
cofibrations are defined to be the monomorphisms (i.e., the  levelwise injective maps), and the fibrations
are the so-called Kan fibrations, defined as the maps having the lifting property of the
inclusions of the simplicial horns $\Lambda^{k,n} \subset \Delta^{n}$ (see \cite[II \S3]{quil}, \cite[3.2]{hove}). 
The homotopy category 
$Ho(sSet)$ is equivalent to the category $Ho(Top)$, via the geometric realization functor, and
can be described as the fibrant simplicial sets (also known as the \emph{Kan complexes}) 
together with the homotopy classes of maps.

\end{itemize}

\subsubsection*{Quillen adjunction, homotopy (co)limits and mapping spaces}
To finish this paragraph on model category theory we mention quickly
the notions of Quillen adjunctions (the natural notion of functors
between model categories), as well as the important notions of
homotopy (co)limits and mapping spaces.

First of all, for two model categories $C$ and $D$, a Quillen
adjunction between $C$ and $D$ consists of a pair of adjoint functors
$g : C \rightleftarrows D : f$ (here $g$ is the left adjoint), such
that either $f$ preserves fibrations and trivial fibrations, or
equivalently $g$ preserves cofibrations and trivial cofibrations. The
main property of a Quillen adjunction as above is to induce an
adjunction on the level of homotopy categories, $\mathbb{L}g : C
\rightleftarrows D : \mathbb{R}f$.  Here $\mathbb{L}g$ and
$\mathbb{R}f$ are respectively the left and right derived functor
deduced from $f$ and $g$, and defined by pre-composition with a
cofibrant (resp. fibrant) replacement functor (see \cite[I \S4
  Thm. 3]{quil}, \cite[\S1.3.2]{hove}).  The typical example of a
Quillen adjunction is given by the geometric realization, and the
singular simplex constructions $|-| : sSet \rightleftarrows Top :
Sing$, between simplicial sets and topological spaces.  This one is
moreover a Quillen equivalence, in the sense that the induced
adjunction at the level of homotopy categories is an equivalence of
categories (see \cite[\S1.3.3]{hove}. Another typical example is
given by a morphism of rings $R \rightarrow R'$ and the corresponding
base change functor $R'\otimes_{R}- : C(R) \longrightarrow C(R')$ on
the level of complexes of modules. This functor is the left adjoint of
a Quillen adjunction (also called \emph{left Quillen}) when the
categories $C(R)$ and $C(R')$ are endowed with the projective model
structures described before (it is no longer a Quillen adjunction for
the injective model structures, except in some very exceptional
cases).

For a model category $C$ and a small category $I$, we can form the category $C^{I}$ of functors 
from $I$ to $C$. The category $C^{I}$ possesses a notion of equivalences induced from the equivalences in 
$C$, and defined as the natural transformations that are levelwise in $W$ (i.e., 
their evaluations at each object $i\in I$ is an equivalence in $C$). With mild extra assumptions on $C$, 
there exists two possible definitions of a model structure on $C^{I}$ whose equivalences are 
the  levelwise equivalences: the projective model structure for which the fibrations are defined levelwise, 
and the injective model structure for which the cofibrations are defined levelwise (see
\cite[Prop.~A.2.8.2]{luri2}). 
We have a
constant diagram functor $c : C \longrightarrow C^{I}$, sending an object of $C$ to the corresponding
constant functor $I \rightarrow C$. The functor $c$ is left Quillen for the injective model structure
on $C^{I}$, and right Quillen for the projective model structure. We deduce a functor at the level 
of homotopy categories $c : Ho(C) \longrightarrow Ho(C^{I})$, which possesses both a right and a left
adjoint, called respectively the \emph{homotopy limit} and \emph{homotopy colimit} functors, 
and denoted by
$Holim_{I}, \; Hocolim_{I} : Ho(C^{I}) \longrightarrow Ho(C)$ (see \cite[Prop. A.2.8.7]{luri5}
as well as comments \cite[A.2.8.8,A.2.8.11]{luri5}).

The homotopy limits and colimits are the right notions of limits and colimits in the setting of 
model category theory and formally behave as the standard notions of limits and colimits. 
They can be used in order to see that the homotopy category $Ho(C)$ of any model category $C$ 
has a natural further enrichment in simpicial sets. For an object $x \in C$ and a simplicial set
$K \in sSet$, we can define an object $K\otimes x \in Ho(C)$, by setting
$$K\otimes x := Hocolim_{\Delta(K)}x \in Ho(C),$$
where $\Delta(K)$ is the category of simplices in $K$ (any category whose geometric realization
gives back $K$ up to a natural equivalence would work), and $x$ is considered as a
constant functor $\Delta(K) \longrightarrow C$. With this definition, it can be shown 
that for two objects $x$ and $y$ in $C$, there is a simplicial set $Map_{C}(x,y) \in sSet$, with natural 
bijections
$$[K,Map_{C}(x,y)] \simeq [K\otimes x,y],$$
where the left hand side is the set of maps in $Ho(sSet)$, and the right hand side the set of maps
in $Ho(C)$. The simplicial sets $Map_{C}(x,y)$ are called the \emph{mapping spaces} of the model 
category $C$, and can alternatively be described using the so-called simpicial and 
cosimplicial resolutions (see \cite[\S5.4]{hove}). 
Their existence implies that the localized category $Ho(C)$ inherits of an extra structure of a simplicial
enrichment, induced by the model category $C$. It is important to understand that this
enrichment only depends on the pair $(C,W)$ of a category $C$ and a class of equivalences $W$. We will
see in the next paragraph that this simplicial enrichment is part of an \emph{$\sz$-categorical structure}, 
and that the correct manner to understand it is by introducing the $\sz$-categorical version 
of the localization construction $(C,W) \mapsto W^{-1}C$. This refined version of the localization 
produces a very strong bridge between model categories and $\sz$-categories, part of which we will 
recall in below.

\subsubsection*{$\sz$-Categories}
An $\sz$-category\footnote{Technically speaking we are only considering
  here $(1,\sz)$-categories, which is a particular case of the more
  general notion of a $\sz$-category that we will not consider in this
  paper.}  is a mathematical structure very close to that of a
category. The main difference is that morphisms in an $\sz$-category
are not elements of a set anymore but rather points in a topological
space (and we think of a set as discrete topological space).  The new
feature is therefore that morphisms in an $\sz$-category can be
\emph{deformed} by means of continuous path inside the space of
morphisms between two objects, and more generally morphisms might come
in continuous family parametrized by an arbitrary topological space,
as for instance a higher dimensional simplex. This is a way to formalize
the notion of homotopy between morphisms often encountered, for
instance in homological algebra where two maps of complexes can be
homotopic.

In the theory of $\sz$-categories the spaces of morphisms are only considered 
up to weak homotopy equivalence for which it is very common to use the notion of 
simplicial sets as a combinatorial model (see \cite[\S3]{hove} for more about the 
homotopy theory of simplicial sets that we use below). This justifies the following definition.

\begin{df}\label{d2}
An \emph{$\sz$-category} $T$ is a simplicial enriched category.
\end{df}

Unfolding the definition, an $\sz$-category consists of the following
data. 

\begin{enumerate}
\item A set  $Ob(T)$, called the set of objects of $T$.

\item For two objects $x$ and $y$ in $T$ a simplicial set of morphisms
$T(x,y)$.

\item For any object $x$ in $T$ a $0$-simplex $id_{x} \in T(x,x)_{0}$.

\item For any triple of objects $x$, $y$ and $z$ in $T$ a 
map of simplicial sets, called the composition
$$T(x,y) \times T(y,z) \longrightarrow T(x,z).$$
\end{enumerate}
These data are moreover required to satisfy an obvious associativity and unit condition.

\begin{remark}\label{r1}
The notion of $\sz$-category of \ref{d2} is not the most general notion of $\sz$-category, 
and rather refers to \emph{semi-strict} $\sz$-categories. Semi-strict refers here to 
the fact that the associativity is strict rather than merely 
satisfied up to a natural homotopy, which itself would satisfy higher homotopy coherences. 
Various other notions of $\sz$-categories for which the compositions is only 
associative up to a coherent set of homotopies are gathered in \cite{berg,lein}. We
will stick to the definition above, as it is in the end equivalent to any other
notion of $\sz$-category and also because it is very easily defined. 
A drawback of this choice will be in the definition of $\sz$-functors and $\sz$-categories of
$\sz$-functors which will be described below and for which some extra care is necessary. 
The definition \ref{d2} seems to us the most efficient in terms  of energy one must spend
on learning the notion of the $\sz$-category theory, particularly for readers 
who do not wish to devote too much effort on the foundational  aspects of derived
algebraic geometry. The theory of $\sz$-categories as defined in \ref{d2} and 
presented below is however the minimum required in order to deal with
meaningful definitions in derived algebraic geometry. 
\end{remark}

There is an obvious notion of a morphism $f : T \longrightarrow T'$ between $\sz$-categories. It 
consists into the following data. 
\begin{enumerate}
\item A map of sets $f : Ob(T) \longrightarrow Ob(T')$.

\item For every pair of objects $(x,y)$ in $T$, a morphism of simplicial sets
$$f_{x,y} : T(x,y) \longrightarrow T'(f(x),f(y)).$$
\end{enumerate}

These data are required to satisfy an obvious compatibility with units and compositions
in $T$ and $T'$. 
These morphisms will be called \emph{strict $\sz$-functors}, as opposed
to a more flexible, and not equivalent, notion of $\sz$-functors that we will introduce later. 
The $\sz$-categories and strict $\sz$-functors form a category denoted by 
$\scat$. 
A category $C$ defines an $\sz$-category by considering the set of morphisms in $C$ as
constant simplicial sets. In the same way, a functor between categories induces
a strict $\sz$-functor between the corresponding $\sz$-categories. This defines
a full embedding
$Cat \hookrightarrow \scat$, 
from the category of categories to the category of $\sz$-categories and strict $\sz$-functors. 
This functor admits a left adjoint 
$$[-] : \scat \longrightarrow Cat.$$
This left adjoint sends an $\sz$-category $T$ to the category $[T]$ 
having the same set of objects as $T$ and whose set of morphisms
are the set of connected components of the simplicial sets of morphisms in $T$. With a formula:
$[T](x,y)=\pi_{0}(T(x,y))$. The category $[T]$ will be referred to the \emph{homotopy category of $T$}.

\subsubsection*{The $\sz$-categories of spaces and of complexes}
We mention here two major examples of $\sz$-categories, the $\sz$-category of Kan complexes, 
and the $\sz$-category of complexes of cofibrant modules over some ring $R$. 
Let $sSets$ be the category of simplicial sets which is naturally enriched over itself by using
the natural simplicial sets of maps, and thus is an $\sz$-category in the sense above. 
We let $\Top$ be the full sub-$\sz$-category of $sSets$ consisting of Kan simplicial sets (i.e. fibrant 
simplicial sets, see \cite[\S3.2]{hove}). 
The homotopy category $[\Top]$ is naturally equivalent to the usual homotopy category of spaces.

For a ring $B$, we let $C(B)$ be the category of (unbounded) cochain complexes of $B$-modules.
It has a natural enrichment in simplicial sets defined as follows. For two complexes
$M$ and $N$, we define the simpicial set $Map(M,N)$ by defining the formula
$$Map(M,N)_{n} := Hom_{C(B)}(M\otimes_{\mathbb{Z}} C_{*}(\Delta^{n}),N),$$
where $C_{*}(\Delta^{n})$ denotes the normalized chain complex 
of homology of the standard simplex $\Delta^{n}$. This makes $C(B)$ into an $\sz$-category. We consider
$L(B) \subset C(B)$ the full sub-$\sz$-category of $C(B)$ consisting of all cofibrant complexes
of $B$-modules (see \cite[\S2.3]{hove}). The homotopy category of
$L(B)$ is naturally equivalent to $D(B)$, the unbounded derived category of complexes of $B$-modules. 

\subsubsection*{The homotopy theory of $\sz$-categories} 
Before going further into $\sz$-category theory we fix some terminology. A \emph{morphism}
in a given $\sz$-category is simply a $0$-simplex of one of the simplicial sets
of maps $T(x,y)$. Such a morphism is an \emph{equivalence} in $T$ if its projection
as a morphism in the category $[T]$ is an isomorphism. Finally, we will sometimes
use the notation $Map_{T}(x,y)$ for $T(x,y)$. \smallbreak

A key notion is the following definition of an equivalence of $\sz$-categories. 

\begin{df}\label{d3}
A strict $\sz$-functor $f : T \longrightarrow T'$ is an \emph{equivalence of $\sz$-categories}
if it satisfies the two conditions below.
\begin{enumerate}
\item $f$ is fully faithful: for all objects $x$ and $y$ in $T$ the map $T(x,y) \longrightarrow 
T'(f(x),f(y))$ is a weak homotopy equivalence of simplicial sets.
\item $f$ is essentially surjective: the induced functor $[f] : [T] \longrightarrow [T']$
is an essentially surjective functor of categories.
\end{enumerate}
\end{df}

The theory of $\sz$-categories up to equivalence will be our 
general setting for derived algebraic geometry, it replaces the setting of categories and
functors that is customary in algebraic geometry. This is unfortunately not an easy theory and
it requires a certain amount of work in order to extend some of the standard constructions
and notions of usual category theory. The good news is that this work has been done and
written down by many authors, we refer for instance to \cite{luri2,simp5} 
(see also \cite[\S1]{toenvezz6}). In the paragraph below
we extract from these works the minimum required for the sequel of our exposition. 
These properties state that $\sz$-categories up to equivalence behave very much likely
as categories up to equivalence of categories, and thus that the basic categorical notions
such as categories of functors, adjunctions, limits and colimits, Yoneda embedding, 
Kan extensions\ldots, all have extensions to the setting of $\sz$-categories. %\smallbreak

 We have seen that the category 
$\scat$ of $\sz$-categories and strict $\sz$-functors possesses a class $W$ of equivalences
of $\sz$-categories. We set $Ho(\scat):=W^{-1}\scat$
the category obtained from $\scat$ by formally inverting the morphisms in $W$,
and call it the homotopy category of $\sz$-categories. The set of morphisms 
in $Ho(\scat)$ will be denoted by
$[T,T'] := Ho(\scat)(T,T').$
The category $H(\scat)$ is a reasonable object because of the existence 
of a model structure on the category $\scat$ which can be used in order to
control the localization along equivalences of $\sz$-categories (see \cite{berg}). 
The sets of morphisms in $Ho(\scat)$ also have explicit descriptions in terms
of equivalent classes of bi-modules (see for instance \cite[\S4.1 Cor. 1]{toen7}, 
for the statement on the setting of dg-categories).

\subsubsection*{(nonstrict) $\sz$-Functors}
By definition an (nonstrict)
$\sz$-functor between two $\sz$-categories $T$ and $T'$ is an element 
in $[T,T']$. This definition only provides a set $[T,T']$ of 
$\sz$-functors, which can be promoted to a full $\sz$-category as follows. 
It can be proved that the category $Ho(\scat)$ is cartesian closed: 
for any pair of $\sz$-categories $T$ and $T'$ there is an object
$\fun(T,T') \in Ho(\scat)$, together with functorial (with respect to 
the variable $U$) bijections
$$[U,\fun(T,T')] \simeq [U\times T,T'].$$
The $\sz$-category $\fun(T,T')$ is by definition the $\sz$-category of 
$\sz$-functors from $T$ to $T'$. It is only well defined up to a natural 
isomorphism as an object in $Ho(\scat)$. As for the case of the sets of maps in 
$Ho(\scat)$, the whole object $\fun(T,T')$ can be explicitly described using 
a certain $\sz$-category of fibrant and cofibrant bi-modules. 

\subsubsection*{Adjunctions, limits and colimits}
The existence of $\sz$-categories of $\sz$-functors can be used in order to 
define adjunctions between $\sz$-categories, and related notions such as 
limits and colimits. We say that an $\sz$-functor $f \in \fun(T,T')$ \emph{has a right adjoint}
if there exists $g$ an object in $\fun(T',T)$ and a morphism $h : id \rightarrow gf$ 
in $\fun(T,T)$, such that for all $x \in T$ and $y\in T'$ the composite
morphism
$$\xymatrix{T'(f(a),b) \ar[r]^-{g} & T(gf(a),g(b)) \ar[r]^-{h} & T(a,g(b))}$$
is a weak equivalence of simplicial sets. It can be shown that if 
$f$ has a right adjoint then the right adjoint $g$ is unique (up to equivalence). The notion
of left adjoint is defined dually.
We say that an $\sz$-category $T$ possesses (small) colimits (resp.\ limits) 
if for all (small) $\sz$-category 
$I$ the constant diagram $\sz$-functor $c : T \longrightarrow \fun(I,T)$
has a left (resp. right) adjoint. The left adjoint (resp. right adjoint), 
when it exists is simply denoted
by $colim_{I}$ (resp. $lim_{I}$). 

\subsubsection*{Yoneda, prestacks and left Kan extensions}
We remind the $\sz$-category $\Top$ consisting of Kan simplicial sets.
Any $\sz$-category $T$ has an $\sz$-category of prestacks $Pr(T)$, also denoted by
$\widehat{T}$, and defined to be $\fun(T^{op},\Top)$, the $\sz$-category of contravariant $\sz$-functors
from $T$ to $\Top$. There is a Yoneda $\sz$-functor
$h : T \longrightarrow \widehat{T},$
which is adjoint to the $\sz$-functor $T : T\times T^{op} \longrightarrow \Top$
sending $(x,y)$ to $T(x,y)$ (when $T$ does not have fibrant hom simplicial sets this
definition has to be pre-composed with chosing a fibrant replacement for $T$). The $\sz$-functor
$h$ is full faithful, and moreover, for all $F\in \widehat{T}$ we have a canonical
equivalence of simplicial sets 
$\widehat{T}(h_{x},F) \simeq F(x).$
The Yoneda embedding $h : T \longrightarrow \widehat{T}$ can also be characterized 
by the following universal property. For every $\sz$-category $T'$ that admit colimits,
the restriction $\sz$-functor
$$- \circ h : \fun_{c}(\widehat{T},T') \longrightarrow \fun(T,T')$$
is an equivalence of $\sz$-categories, where $\fun_{c}(\widehat{T},T')$ is
the full sub-$\sz$-category of $\fun(\widehat{T},T')$  consisting of 
$\sz$-functors that commute with colimits (see \cite[Thm. 5.1.5.6]{luri2}). The inverse
$\sz$-functor $\fun(T,T') \longrightarrow \fun_{c}(\widehat{T},T')$
is called the \emph{left Kan extension}. 

\subsubsection*{Localization and model categories} An important source of $\sz$-categories
come from localization, the process of making some morphisms to be invertible in 
a universal manner. For a category $C$ and a subset $W$ of morphisms in $C$, 
there is an $\sz$-category   $L(C,W)$ together with an $\sz$-functor
$l : C \longrightarrow L(C,W),$
such that for any $\sz$-category $T$, the restriction through $l$ induces an 
equivalence of $\sz$-categories
$$\fun(L(C,W),T) \simeq \fun_{W}(C,T),$$
where $\fun_{W}(C,T)$ denotes the full sub-$\sz$-category of $\fun(C,T)$ consisting
of all $\sz$-functors sending $W$ to equivalences in $T$. 
It can be shown that a localization always exists (see \cite[Prop. 8.7]{hirssimp}, 
see also \cite[\S4.3]{toen7} for dg-analogue), 
and is equivalent to the so-called Dwyer-Kan simplicial localization of \cite{dwyekan1}. 
The homotopy category $Ho(L(C,W))$ is canonically equivalent to the localized
category $W^{-1}C$ in the sense of Gabriel-Zisman (see \cite[\S1]{gabrzism}). In general
$L(C,W)$ is not equivalent to $W^{-1}C$, or in other words its mapping spaces are
not $0$-truncated. The presence of nontrivial higher homotopy in 
$L(C,W)$ is one justification of the importance of $\sz$-categories in many domains
of mathematics. 

When $C$ is moreover a simplical model category, and $W$ its subcategory of weak equivalence,
the localization $L(C,W)$ is simply denoted by $L(C)$, 
and can be described, up to a natural equivalence, as the
simplicially enriched category $\underline{C}^{cf}$ of fibrant and cofibrant 
objects in $C$ (see \cite{dwyekan2}). Without the simplicial assumption, for a general
model category $C$ a similar result is true but involves mapping spaces as defined in \cite{dwyekan2}
and \cite[\S5.4]{hove}
using simplicial and cosimplicial resolutions. 
For a model category $C$, and small category $I$, let $C^{I}$ be the
model category of diagrams of shape $I$ in $C$. It is shown in 
\cite[\S18]{hirssimp} (see also \cite[Prop. 4.2.4.4]{luri2}) that 
there exists a natural equivalence of $\sz$-categories
$$L(M^{I}) \simeq \fun(I,L(M)).$$
This is an extremely useful statement that can be used in order to provide
natural models for most of the $\sz$-categories encountered in practice. One important consequence
is that the $\sz$-category $L(M)$ always has limits and colimits, and moreover that these
limits and colimits in $L(M)$ can be computed using the well known 
homotopy limits and homotopy colimits of homotopical algebra (see \cite{dwyehirskan} for 
a general discussion about homotopy limits and colimits). 

\subsubsection*{The $\sz$-category of $\sz$-categories} The localization construction we just described can
be applied to the category $\infCat$, of $\sz$-categories and strict $\sz$-functors, together
with $W$ being equivalences of $\sz$-categories of our definition \ref{d3}. We thus have an $\sz$-category of $\sz$-categories
$$\uscat := L(\infCat).$$
The mapping spaces in $\uscat$ are closely related to the $\sz$-category of $\sz$-functors
in the following manner. For two $\sz$-categories $T$ and $T'$, we consider
$\fun(T,T')$, and the sub-$\sz$-category $\fun(T,T')^{eq}$ consisting
of $\sz$-functors and equivalences between them. The $\sz$-category $\fun(T,T')^{eq}$ has a geometric realization 
$|\fun(T,T')^{eq}|$, 
obtained by taking nerves of each categories of simplicies, and then the diagonal
of the corresponding bi-simplicial set (see below). We have a weak equivalence of simplicial sets
$$Map_{\uscat}(T,T') \simeq |\fun(T,T')^{eq}|,$$
which expresses the fact that mapping spaces in $\uscat$ are the spaces of $\sz$-functors
up to equivalence. 
Another important aspect is that $\uscat$ possesses all limits and colimits. This
follows for instance from the existence of a model structure on simplicially 
enriched categories (see \cite{berg}). We refer to \cite[Cor. 18.7]{hirssimp}
for more about how to compute the limits in $\uscat$ in an explicit manner.

A second important fact concerning the $\sz$-category of $\sz$-categories is
the notion of $\sz$-groupoids. The $\sz$-groupoids are defined to be 
the $\sz$-categories $T$ for which the homotopy category $[T]$ is a groupoid, or in other
words for which every morphism is an equivalence. If we denote by
$\bfinfgpd$ the full sub-$\sz$-category of $\uscat$ consisting of
$\sz$-groupoids, then the nerve construction (also called the geometric
realization) produces an equivalence of $\sz$-categories
$$|-| : \bfinfgpd\simeq \Top.$$
The inverse of this equivalence is the fundamental $\sz$-groupoid
construction $\Pi_{\sz}$ (denoted by $\Pi_{1,Se}$ in
\cite[\S2]{hirssimp}). The use of the equivalence above will be mostly 
implicit in the sequel, and we will allow ourselves to consider any simplicial sets
$K\in \Top$ as an $\sz$-category through this equivalence. 

\subsubsection*{$\sz$-Topos and stacks} For an $\sz$-category $T$, a
Grothendieck topology on $T$ is by definition a Grothendieck topology
on $[T]$ (see e.g.\ \cite[Exp.~II]{sga4}).  When such a topology
$\tau$ is given, we can define a full sub-$\sz$-category $St(T,\tau)
\subset Pr(T)$, consisting of prestacks satisfying a certain descent
condition. The descent condition for a given prestack $F :
T^{op}\rightarrow \Top$, expresses that for any augmented simplicial
object $X_{*} \rightarrow X$ in $Pr(T)$, which is a
$\tau$-hypercoverings, the natural morphism
$$Map_{Pr(T)}(X,F) \longrightarrow \displaystyle{\lim_{[n] \in \Delta}}Map_{Pr(T)}(X_{n},F)$$
is a weak equivalence in $\Top$ (and $lim$ is understood in the
$\sz$-categorical sense, or equivalently as a homotopy limit of simplicial sets). 
Here, $\tau$-hypercoverings are generalizations of nerves
of covering families and we refer to \cite[Def. 3.2.3]{hagI} for the precise definition in the context
of $\sz$-categories. The condition
above is the $\sz$-categorical analogue of the sheaf condition, and a prestack
satisfying the descent condition will be called a \emph{stack} (with respect to the 
topology $\tau$). The $\sz$-category $St(T,\tau)$, also denoted by $T^{\sim,\tau}$
is the $\sz$-category of stacks on $T$ (with respect to $\tau$) and is
an instance of an $\sz$-topos (all $\sz$-topos we will have to consider in this paper are of this form).
The descent condition can also be stated by an $\sz$-functor $F : T^{op} \longrightarrow
\mathcal{C}$, where $\mathcal{C}$ is another $\sz$-category with all limits. This will allow
us to talk about stacks of simplicial rings, which will be useful in the definition of derived
schemes we will give below (Def. \ref{d1'}). 
We refer to \cite{luri2} for more details about $\sz$-topos, and to \cite{hagI}
for a purely model categorical
treatment of the subject. 

\subsubsection*{Stable $\sz$-categories} Stable $\sz$-categories are the 
$\sz$-categorical counter-part of triangulated categories. We recall here the most basic
definition and the main property as they originally appear in \cite[\S7]{toenvezz4}, and we refer to \cite{dag} for more details. 

We say that an $\sz$-category $T$ is \emph{stable} if it has finite limits and colimits, 
if the initial object is also final, and if the loop endo-functor $\Omega_{*} : x \mapsto *\times_{x}*$
defines an equivalence of $\sz$-categories $\Omega_{*} : T \simeq T$. It is known that 
the homotopy category $[T]$ of a stable $\sz$-category $T$ possesses a canonical 
triangulated structure for which the distinguished triangles are the image of
fibered sequences in $T$. 
If $M$ is a stable model category (in the sense of \cite[\S7]{hove}), then $LM$ is a stable $\sz$-category.
For instance, if $M=C(k)$ is the model category of complexes of modules over some ring $k$, 
$LM$, the $\sz$-category of complexes of $k$-modules, is stable.

\begin{war}\label{war}
Before we start using the language of $\sz$-categories, we warn the reader that
we will use this language in a rather loose way and that most of our constructions will
be presented in a naively manner. Typically a $\sz$-category will be given by describing
its set of objects and simplicial sets of maps between two given objects without
taking care of defining compositions and units. Most of the time the compositions and
units are simply obvious, but it might also happen that some extra work is needed
to get a genuine $\sz$-category.  This happens for instance when 
the described mapping space is only well-defined up to  a weak equivalences of 
simplicial sets, or when composition is only defined up to a natural homotopy, 
for which it might be not totally obvious how to define things correctly. This is one
of the  technical difficulties of $\sz$-category theory that we will 
ignore in this exposition, but the reader must keep in mind that in some cases
this can be overcome only by means of a substantial amount of work.
\end{war}

\subsection{Derived schemes}

We are now coming back to our first definition of derived schemes \ref{d1}, 
but from the $\sz$-categorical point of view briefly reminded in the last paragraph. 
We start by considering $\scomm$, the $\sz$-category of simplicial commutative
rings, also called \emph{derived rings}. It is defined by 
$$\scomm := L(sComm),$$
the $\sz$-categorical localization of the category of simplicial commutative rings $sComm$
with respect to the weak homotopy equivalences. The weak homotopy equivalences are 
here morphisms of simplicial commutative rings 
$A \rightarrow B$ inducing a weak homotopy equivalence on the underlying
simplicial sets. 
Any derived ring $A$ provides a commutative graded ring $\pi_{*}(A)=\displaystyle{\oplus_{i\geq 0}} \pi_{i}(A)$, 
where the homotopy groups are all taken with respect to $0$ as a base point. The
construction $A \mapsto \pi_{*}(A)$ defines an $\sz$-functor from the $\sz$-category
$\scomm$ to the category of commutative graded rings. 

For a topological space $X$, there is an $\sz$-category $\scomm(X)$ of stacks 
on $X$ with coefficients in the $\sz$-category of derived rings. If we let 
$Ouv(X)$ the category of open subsets in $X$, $\scomm$ can be identified with the full
sub-$\sz$-category of $\fun(Ouv(X)^{op},\scomm)$, consisting of
$\sz$-functors satisfying the descent condition (see $\S2.1.2$). For a continuous
map $u : X \longrightarrow Y$  we have an adjunction of $\sz$-categories
$$u^{-1} : \scomm(Y) \rightleftarrows \scomm(X) : u_{*}.$$

We start by defining an $\sz$-category $\mathbf{dRgSp}$ of 
derived ringed spaces. Its objects are pairs $(X,\OO_{X})$, where $X$ 
is a topological space and $\OO_{X} \in \scomm(X)$ is a stack of derived rings on $X$.
For two derived ringed spaces $(X,\OO_{X})$ and $(Y,\OO_{Y})$ we set 
$$Map((X,\OO_{X}),(Y,\OO_{Y})) := \coprod_{u : X\rightarrow Y} 
Map_{\scomm(Y)}(\OO_{Y},u_{*}(\OO_{X})).$$
This definition can be promoted to an $\sz$-category $\mathbf{dRgSp}$ 
whose objects are derived ringed spaces and whose simplicial sets of maps are
defined as above. Technically speaking this requires the use of some more advanced
notion such as fibered $\sz$-categories, but can also be realized using
concrete model category structures of sheaves of simplicial commutative rings
(this is a typical example for our warning \ref{war} in a practical situation\footnote{From now on we will 
not refer to the warning \ref{war} anymore}). 

Any derived ringed space $(X,\OO_{X})$ has a truncation 
$(X,\pi_{0}(\OO_{X}))$ that is a (underived) ringed space, 
where $\pi_{0}(\OO_{X})$ denotes here the \emph{sheaf} of connected components. 
We define $\mathbf{dRgSp}^{loc}$ as a (nonfull) sub-$\sz$-category
of $\mathbf{dRgSp}$ consisting of objects $(X,\OO_{X})$ 
whose truncation $(X,\pi_{0}(\OO_{X}))$ is a locally ringed space, and
maps inducing local morphisms on the ringed spaces obtained by truncations. The inclusion
$\sz$-functor
$\mathbf{dRgSp}^{loc} \hookrightarrow \mathbf{dRgSp}$
is not fully faithful but it is faithful in the sense that the morphisms induced
on mapping spaces are inclusions of union of connected components. 

With this new language the $\sz$-category of derived schemes is defined as follows.

\begin{df}\label{d1'}
The \emph{$\sz$-category of derived schemes} is defined to be the full sub-$\sz$-category 
of $\mathbf{dRgSp}^{loc}$ consisting of all objects $(X,\OO_{X})$ 
with the two conditions above satisfied. 
\begin{enumerate}
\item The
truncation $(X,\pi_{0}(\OO_{X}))$ is a scheme.
\item For all $i$ the sheaf of $\pi_{0}(\OO_{X})$-modules $\pi_{i}(X)$
is quasi-coherent.
\end{enumerate}
The $\sz$-category of derived schemes is denoted by $\dSt$.  
\end{df}

A ring can be considered as a constant simplicial ring, and this defines a full embedding
$i : \mathbf{Comm} \hookrightarrow \scomm$, from the category of commutative of rings to the 
$\sz$-category
of derived rings. This inclusion has a left adjoint given by the $\sz$-functor $\pi_{0}$. 
This adjunction extends to an adjunction at the level of derived ringed spaces, derived locally
ringed spaces and derived schemes. We thus have an adjunction
$$t_{0} : \dSch \rightleftarrows \mathbf{Sch} : i,$$
between the $\sz$-category of schemes and the category of schemes. The functor $i$ is
moreover fully faithful, and therefore schemes sit inside derived schemes as a full sub-$\sz$-category. 
The $\sz$-functor $t_{0}$ sends a derived scheme $(X,\OO_{X})$ to the 
scheme $(X,\pi_{0}(\OO_{X}))$. We will often omit to mention the functor $i$ and simply 
considered $\Sch$ as sitting inside $\dSch$ as a full subcategory. By adjunction, 
for any derived scheme $X$ there is a natural morphism of derived schemes
$$j : t_{0}(X) \longrightarrow X.$$

\begin{remark}\label{rnil}
It is an accurate analogy to compare the morphism $j : t_{0}(X) \longrightarrow X$ 
with the inclusion $Y_{red} \hookrightarrow Y$, of the reduced subscheme $Y_{red}$
of a scheme $Y$. In this way, the truncation $t_{0}(X)$ sits inside
the derived scheme $X$, and $X$ can be thought as some sort of infinitesimal
thickening of $t_{0}(X)$, but for which the additional infinitesimals 
functions live in higher homotopical degrees. 
\end{remark}

The truncation $t_{0}$  possesses generalizations $t_{\leq n}$ for various
integers $n\geq 0$ (with $t_{0}=t_{\leq 0}$). Let $X$ be a derived scheme. The stack 
of derived rings $\OO_{X}$ has a Postnikov tower
$$
\OO_{X}\to\dots \to t_{\leq n}(\OO_{X}) \to  t_{\leq n-1}(\OO_{X}) \to 
\dots \to t_{0}(\OO_{X})=\pi_{0}(\OO_{X}).
$$
It is characterized, as a tower of morphisms in the $\sz$-category of stacks of derived 
rings on $X$, by the following two properties.
\begin{itemize}
\item For all $i>n$, we have $\pi_{i}(t_{\leq n}(\OO_{X}))\simeq 0$.
\item For all $i\leq n$, the morphism $\OO_{X} \longrightarrow t_{\leq n}(\OO_{X})$
induces isomorphisms $\pi_{i}(\OO_{X}) \simeq \pi_{i}(t_{\leq n}(\OO_{X}))$. 
\end{itemize}

Each derived ringed space $(X,t_{\leq n}(\OO_{X}))$ defines 
a derived scheme, denoted by $t_{\leq n}(X)$, and the above tower defines a diagram 
of derived schemes
$$
t_{0}(X) \to t_{\leq 1}(X) \to \dots \to t_{\leq n}(X) \to
t_{\leq n+1}(X) \to \dots \to X.
$$
This diagram exhibits $X$ as the colimit of the derived schemes $t_{\leq n}(X)$
inside the $\sz$-category $\dSch$. The Postnikov tower of derived schemes is a powerful tool
in order to understand maps between derived schemes and more generall mapping spaces. Indeed, 
for two derived schemes $X$ and $Y$ we have
$$Map_{\dSch}(X,Y) \simeq \displaystyle{\lim_{n\geq 0}} Map_{\dSch}(t_{\leq n}X,t_{\leq n}Y)
\simeq \displaystyle{\lim_{n\geq 0}} Map_{\dSch}(t_{\leq n}X,Y),$$
which presents that mapping spaces as a (homotopy) limit of simpler mapping spaces. 
First of all, for a given  $n$, the mapping space 
$Map_{\dSch}(t_{\leq n}X,t_{\leq n}Y)$ is automatically  $n$-truncated (its nontrivial
homotopy is concentrated in degree less or equal to  $n$). Moreover, the projection
$Map_{\dSch}(t_{\leq n+1}X,t_{\leq n+1}Y) \longrightarrow Map_{\dSch}(t_{\leq n}X,t_{\leq n}Y)$
can be understood using obstruction theory, as this will be explained in
 our section $\S4.1$): the description of the
fibers of this projection consists essentially into a linear problem of understanding
some specific extensions groups of sheaves of modules. 

The above picture of Postnikov towers 
is very analogous to the situation with formal schemes: any formal scheme
$X$ is a colimit of schemes $X_{n}$ together with closed immersions
$X_{n} \hookrightarrow X_{n+1}$ corresponding to a square zero ideal sheaf on $X_{n+1}$. 
This analogy with formal scheme is a rather accurate one. \smallbreak

To finish this paragraph we mention some basic examples of derived schemes and mapping
spaces between derived schemes. More advanced examples will be given in the next paragraph 
and later. 

\subsubsection*{Affine derived schemes} We let $\dAff$ be the full sub-$\sz$-category of $\dSch$
consisting of derived schemes $X$ whose truncation $\pi_{0}(X)$ is an affine scheme. Objects
in $\dAff$ are called \emph{affine derived schemes}. 
We have an $\sz$-functor of global functions
$$\mathbb{H}(-,\OO_{X}) : \dAff^{op} \longrightarrow \scomm,$$
sending an affine derived scheme $X$ to $\mathbb{H}(X,\OO_{X}):=p_{*}(\OO_{X})$, where
$p : X \longrightarrow *$ is the canonical projection, 
and $p_*$ is the induced $\sz$-functor on $\sz$-categories of stacks of derived rings. 
The $\sz$-functor $\mathbb{H}$ can be shown to be an equivalence of $\sz$-categories. The inverse
$\sz$-functor of $\mathbb{H}$ is denoted by $\Spec$, and can be described as follows. Let $A$ be 
a simplicial commutative ring. We consider the (underived) affine scheme $S=Spec\, A_{0}$, 
the spectrum of the ring of $0$-dimensional simplicies in $A$.
The simplicial ring $A$ is in a natural way a simplicial commutative $A_{0}$-algebra (through
the natural inclusion $A_{0} \longrightarrow A$) and thus
defines a sheaf of simplicial quasi-coherent $\OO_{S}$-modules $\mathcal{A}$ on $S$. 
This sheaf defines a stack of derived rings on $S$ and thus an object on 
$\scomm(S)$. We denote by $X\subset S$ the closed subset defined by $Spec\, \pi_{0}(A)$
(note that the ring $\pi_{0}(A)$ is a quotient of $A_{0}$).
By construction, the stack of derived rings $\mathcal{A}$ is supported on the closed
subspace $X$, in the sense that its restriction on $S-X$ is equivalent to
$0$. This implies that it is equivalent to a stack of derived rings of the form $i_{*}(\OO_{X})$
for a well defined object $\OO_{X} \in \scomm(X)$ ($i_{*}$ produces an equivalence between 
stacks of derived rings on $S$ supported on $X$ and stacks of derived rings on $X$).
The derived ringed space $(X,\OO_{X})$ is denoted by $\Spec\, A$ and is an affine derived scheme.
The two constructions $\mathbb{H}$ and $\Spec$ are inverse to each other. 

\subsubsection*{Fibered products} The $\sz$-category $\dSch$ of derived schemes has all finite
limits (see \cite[\S1.3.3]{hagII}). 
The final object is of course $*=Spec\, \mathbb{Z}$. On the level of affine derived
schemes fibered products are described as follows. A diagram of 
affine derived schemes $\xymatrix{X \ar[r] & S & \ar[l] Y}$ defines, by taking global functions, a 
corresponding diagram
of derived rings $\xymatrix{A & C \ar[r] \ar[l] & B}$. We consider 
the derived ring $D:=A\otimes^{\mathbb{L}}_{C}B \in \scomm$.
 From the point of view of $\sz$-categories
the derived ring $D$ is the push-out of the diagram $\xymatrix{A & C \ar[r] \ar[l] & B.}$ 
It can be constructed explicitly by replacing $B$ by a simplicial $C$-algebra $B'$ that is 
a cellular $C$-algebra (see \cite[\S2.1]{toenvaqu} for the general notion 
of cellular objects), and then considering the naive levelwise tensor product
$A\otimes_{C}B'$. For instance, when $A$, $B$ and $C$ are all commutative rings
then $D$ is a simplicial commutative ring with the property that 
$\pi_{n}(D) \simeq Tor_{n}^{C}(A,B).$
In general, for a diagram of derived schemes $\xymatrix{X \ar[r] & S & \ar[l] Y}$, 
the fibered product $X\times_{S}Y$ can be described by gluing the local affine pictures
as above. Again, when $X$, $Y$ and $S$ are merely (underived) schemes, 
$Z:=X\times_{S}Y$ is a derived scheme whose truncation is the usual 
fibered product of schemes. The homotopy sheaf of the derived structure sheaf 
$\OO_{Z}$ are the higher Tor's 
$$\pi_{n}(\OO_{Z})\simeq \mathcal{T}or_{n}^{\OO_{S}}(\OO_{X},\OO_{Y}).$$
We see here the link with Serre's intersection formula discussed at the beginning of
$\S1$. 

We note that the inclusion functor $i : \Sch\hookrightarrow \dSch$ does not 
preserve fibered products in general, except under the extra condition of Tor-independence (e.g., 
if one of the maps is flat). In contrast to this, the truncation $\sz$-functor 
$t_{0}$ sends fibered products of derived schemes to fibered products of schemes. 
This is a source of a lot of examples of interesting  derived schemes, 
simply by constructing a derived fibered product of schemes. A standard example is
the derived fiber of a nonflat morphism between schemes. 

\subsubsection*{Self-intersections} 
Let $Y \subset X$ be a closed immersion of schemes, and consider 
the derived scheme $Z:=Y\times_{X}Y \in \dSch$. The truncation $t_{0}(Z)$
is isomorphic to the same fibered product computed in $\Sch$, and thus is
isomorphic to $Y$. The natural morphism $t_{0}(Z)\simeq Y \longrightarrow Z$
is here induced by the diagonal $Y \longrightarrow Y\times_{X}Y$. The projection
to one of the factors produces a morphism of derived schemes $Z \rightarrow Y$ that is
a retraction of $Y \rightarrow Z$. This is an example of a \emph{split} derived scheme $Z$: 
the natural map $t_{0}(Z) \longrightarrow Z$ admits a retraction (this is not 
the case in general). For simplicity we assume that $Y$ is a local complete intersection
in $X$, and we let $\mathcal{I} \subset \OO_{X}$ be its ideal sheaf. The conormal bundle of
$Y$ inside $X$ is then $\mathcal{N}^{\vee}\simeq \mathcal{I}/\mathcal{I}^{2}$, which is
a vector bundle on $Y$.

When $X=\Spec\, A$ is affine, and $Y=\Spec\, A/I$, the derived scheme
$Z$ can be understood in a very explicit manner. Let $(f_1,\dots,f_{r})$ en regular
sequence generating $I$. We consider the derived ring $K(A,f_*)$, which is obtained
by freely adding a $1$-simplices $h_{i}$ to $A$ such that $d_{0}(h_{i})=0$ and 
$d_{1}(h_{1})=f_{i}$ (see \cite{toen3}, proof of proposition 4.9, 
for details). The derived ring $K(A,f)$ has a natural
augmentation $K(A,f) \longrightarrow A/I$ which is an equivalence because
the sequence is regular. It is moreover a cellular $A$-algebra, by construction, and thus
the derived ring $A/I \otimes^{\mathbb{L}}_{A} A/I$ can be identified with
$B=K(A,f)\otimes_{A}A/I$. This derived ring is an $A/I$-algebra such that $\pi_{1}(B)\simeq I/I^{2}$. 
As $I/I^2$ is a projective $A/I$-module we can represent the isomorphism 
$\pi_{1}(B)\simeq I/I^{2}$ by a morphism of simplicial $A/I$-modules
$I/I^{2}[1] \longrightarrow B$, where $[1]$ denotes the suspension in the $\sz$-category
of simplicial modules. This produces a morphism of derived rings
$Sym_{A/I}(I/I^{2}[1]) \longrightarrow B$, where 
$Sym_{A/I}$ denotes here the $\sz$-functor sending an $A/I$-module $M$
to the derived $A/I$-algebras it generates. 
This morphism is an equivalence in characteristic zero, and thus we have in this case
$$Z \simeq \Spec\, (Sym_{A/I}(I/I^{2}[1])).$$
In nonzero characteristic a similar but weaker statement is true, we have
$$Z \simeq \Spec\,B$$
where now the right hand side is not a free derived ring anymore, but 
satisfies
$$\pi_*(B)\simeq \oplus_{i\geq 0} \wedge^{i}(I/I^{2})[i].$$

The local 
computation we just made shows that the sheaf of graded $\OO_{Y}$-algebras $\pi_{*}(\OO_{Z})$ 
is isomorphic to $Sym_{\OO_{Y}}(\mathcal{N}^{\vee}[1])$. However, the sheaf of derived rings
$\OO_{Z}$ is not equivalent to $Sym_{\OO_{Y}}(\mathcal{N}^{\vee}[1])$ in general (i.e., 
in the non-affine case). It is
locally so in characteristic zero, 
but there are global cohomological obstructions for this to be globally
true. The first of these obstructions is a cohomology class in $\alpha_{Y} \in 
Ext^{2}_{Y}(\mathcal{N}^{\vee},
\wedge^2 \mathcal{N}^{\vee})$, which can be interpreted as follows. We have a natural
augmentation of stacks of derived rings $\OO_{Z} \longrightarrow \OO_{Y}$, which splits 
as $\OO_{Z} \simeq \OO_{Y} \times K$, where we consider this splitting in $D_{\qcoh}(Y)$, the
derived category of quasi-coherent complexes on $Y$. The complex $K$ is cohomogically
concentrated in degrees $]-\infty,1]$, and we can thus consider the exact triangle
$$\xymatrix{
H^{-2}(K)[2] \ar[r] & \tau_{\geq -2}(K) \ar[r] & H^{-1}(K)[1] \ar[r]^-{\partial} & 
H^{-2}(K)[3].}$$
The class $\alpha_{Y}$ is represented by the boundary map $\partial$. 

The obstruction class $\alpha_{Y}$ has been identified with the obstruction for
the conormal bundle $\mathcal{N}^{\vee}$ to extend to the second infinitesimal 
neighbourhood of $Y$ in $X$. The higher obstruction classes live in 
$Ext^{i}_{Y}(\mathcal{N}^{\vee},\wedge^i \mathcal{N}^{\vee})$ and can be shown to vanish if
the first obstruction $\alpha$ does so. We refer to \cite{arincald,griv} for more details on the subject,
and to \cite{calacald} for some refinement. 

\begin{remark}
One of the most important derived self-intersection is the derived loop scheme 
$X\times_{X\times X}X$, which we will investigate in more detail in our $\S4.4$. It behaves
in a particular fashion as the inclusion $X \longrightarrow X\times X$
possesses a global retraction. 
\end{remark}

\subsubsection*{Euler classes of vector bundles} We let $X$ be a scheme
and $V$ a vector bundle on $X$ (considered as a locally 
free sheaf of $\OO_{X}$-modules), together with a section $s \in \Gamma(X,V)$. 
We denote by $\mathbb{V}=\Spec\, Sym_{\OO_{X}}(V^{\vee})$ the total space of $V$, 
considered as a scheme over $X$. The section $s$ and the zero section define
morphisms $\xymatrix{X \ar[r]^-{s} & \mathbb{V} & X \ar[l]_-{0}}$, out of which we can 
form the derived fiber product $X \times_{\mathbb{V}} X$. This derived
scheme is denoted by $Eu(V,s)$, and is called the \emph{Euler class of $V$ with
respect to $s$}. The truncation $t_{0}(Eu(V,s))$ consists of the closed subscheme $Z(s) \subset X$
of zeros of $s$, and the homotopy sheaves of the derived structure sheaf $\OO_{Eu(V,s)}$
controls the defect of Tor-independence of the section $s$ with respect to the zero section.

Locally the structure of $Eu(V,s)$ can be understood using Koszul 
algebras as follows. We let $X=\Spec\, A$ and $V$ be given by a projective $A$-module $M$
of finite type. The section $s$ defines a morphism of $A$-modules $s : M^{\vee} \longrightarrow A$. 
We let $K(A,M,s)$ be the derived ring obtained out of $A$ by freely 
adding $M^{\vee}$ as $1$-simplicies, such that each $m\in M^{\vee}$ has boundary defined by 
$d_{0}(m)=s(m)$ and $d_{1}(m)=0$. This derived ring $K(A,M,s)$ is a simplicial 
version of Koszul resolutions in the dg-setting, and 
$\Spec\, K(A,M,s)$ is equivalent to $Eu(V,s)$. In characteristic zero, 
derived rings can also be modelled by commutative dg-algebras (see $\S3.4$), 
and $K(A,M,s)$ then becomes equivalent to the standard Koszul 
algebra $Sym_{A}(M^{\vee}[1])$ with a differential given by 
$s$. 

\section{Derived schemes, derived moduli\erase{problems}, and derived
  stacks}

In this section, we present the functorial point of view of derived algebraic geometry. 
It consists of viewing derived schemes as certain
($\sz$)-functors defined on simplicial algebras, similarly to the way 
schemes can be considered as functors on the category of algebras (see for instance 
\cite{eiseharr}). This will lead us to the notion 
of \emph{derived moduli problems} and to the representability by derived schemes and more generally,
by \emph{derived Artin stacks}, a derived analogue of 
algebraic stacks (see \cite{lm}), as well as to a far-reaching generalization of derived schemes
obtained by allowing certain quotients by groupoid actions. We will again provide basic examples, as well as
more advanced examples deduced from the Artin-Lurie representability theorem. Finally,
we will mention the existence of many variations of derived algebraic geometry, such 
as derived analytic and differential geometry, derived log geometry, spectral geometry, etc.

\subsection{Some characteristic properties of derived schemes}

We have gathered in this paragraph some properties shared by derived schemes
which are characteristic in the sense that they do not hold in general for schemes
without some extra and nontrivial conditions. They provide a first motivation
for the introduction of derived schemes and clearly show that the theory of derived
schemes has much more regularity than the theory of schemes. We will later see many more
examples. 

\subsubsection*{Base change} A scheme $X$ possesses a quasi-coherent derived 
category $D_{\qcoh}(X)$, which for us will be the derived category of (unbounded)
complexes of $\OO_{X}$-modules with quasi-coherent cohomology sheaves
(see for instance \cite{bondvand}). 
In the same way, a derived scheme $X$ possesses a quasi-coherent derived
$\sz$-category $L_{\qcoh}(X)$, defined as follows. 

 We consider  the 
$\sz$-category $\Zaff (X)$ of affine open derived subschemes $U \subset X$. This $\sz$-category
can be shown to be equivalent to a poset, in fact, through the
functor $U \mapsto t_{0}(U)$, to the poset of open subschemes in $t_{0}(X)$ (see \cite[Prop. 2.1]{reduce}).
For each object $U \in \Zaff(X)$ we have its derived ring of functions $A_{U}:=\mathbb{H}(U,\OO_{U})$.
The simplicial ring $A_{U}$ can be normalized to a commutative dg-algebra $N(A_{U})$, for
which we can consider the category \hbox{$N(A_{U})$-$Mod$} of (unbounded) $N(A_{U})$-dg-modules
(see \cite{shipschw} for more about the monoidal properties of the normalization functor).
Localizing this category along quasi-isomorphisms defines an $\sz$-category
$L_{\qcoh}(U):=L(N(A_{U})\textit{-Mod},\textit{quasi-isom})$.
For each inclusion of open subsets $V \subset U \subset X$, we have a morphism of commutative
dg-algebras $N(A_{U}) \longrightarrow N(A_{V})$ and thus an induced base change $\sz$-functor
$-\otimes_{N(A_{U})}^{\mathbb{L}}N(A_{V}) : L_{\qcoh}(U) \longrightarrow L_{\qcoh}(V).$
This defines an $\sz$-functor 
$L_{\qcoh}(-) : \Zaff(X)^{op} \longrightarrow \uscat,$
which moreover is a stack (i.e., satisfies the descent condition explained in $\S2.1.2$) 
for the Zariski topology. 
We set 
$$L_{\qcoh}(X):=\displaystyle{\lim_{U \in \Zaff(X)^{op}}}L_{\qcoh}(U) \in \uscat,$$
where the limit is taken in the $\sz$-category of $\sz$-categories, and call it the
\emph{quasi-coherent derived $\sz$-category of $X$}. When $X$ is a scheme, $L_{\qcoh}(X)$ is an 
$\sz$-categorical model for the derived category $D_{\qcoh}(X)$ of $\OO_{X}$-modules with
quasi-coherent cohomologies: we have a natural equivalence of categories
$$[L_{\qcoh}(X)] \simeq D_{\qcoh}(X).$$
When $X=\Spec\, A$ is affine for a derived ring $A$, then 
$L_{\qcoh}(X)$ is naturally identified with $L(A)$ the $\sz$-category of 
dg-modules over the normalized dg-algebra $N(A)$.  We will  for
$E \in L(A)$ often use the notation $\pi_{i}(E):=H^{-i}(E)$. Similarly, for a general 
derived scheme $X$, and $E \in L_{\qcoh}(X)$, we have cohomology sheaves $H^{i}(E)$ that
are quasi-coherent on $t_{0}(X)$, and which we  will also denote by $\pi_{i}(E)$ 

For a morphism between derived schemes $f : X \longrightarrow Y$ there is an natural
pull-back $\sz$-functor
$f^{*} : L_{\qcoh}(Y) \longrightarrow L_{\qcoh}(X)$, as well as its right adjoint the push-forward
$f_{*} : L_{\qcoh}(X) \longrightarrow L_{\qcoh}(Y)$. These are first defined locally on the level
of affine derived schemes: the $\sz$-functor $f^{*}$ is induced by the base
change of derived rings whereas the $\sz$-functor $f_{*}$ is a forgetful $\sz$-functor. 
The general case is done by gluing the local constructions
(see \cite[\S4.2]{seat}, \cite[\S1.1]{toen4} for details).

By the formal property of adjunction, for any commutative square of derived schemes
$$\xymatrix{
X' \ar[r]^-{g} \ar[d]_-{q} & X \ar[d]^-{p} \\
Y' \ar[r]_-{f} & Y}$$
there is a natural morphism between $\sz$-functors
$$h : f^{*}p_{*} \Rightarrow q_{*}g{^*} : L_{\qcoh}(X) \longrightarrow L_{\qcoh}(Y').$$
The base change theorem (see \cite[Prop. 1.4]{toen4}) insures that $h$ is an equivalence of $\sz$-functors as
soon as the square is cartesian and all derived schemes are quasi-compact and 
quasi-separated. When all the derived schemes are  schemes and moreover $f$ is flat, 
then $X'$ is again a scheme and the base change formula recovers the usual well known formula
for schemes. When $f$ and $p$ are not Tor-independent the derived scheme $X'$ is not
a scheme and the difference between $X'$ and its truncation $t_{0}(X')$ measures
the \emph{excess of intersection} (see e.g.\ \cite[\S6]{fultlang}). All the classical
excess intersection formulae can be recovered from the base change 
formula for derived schemes. 

\subsubsection*{Tangent complexes, smooth and \'etale maps} 
Let $A$ be a derived ring and $M$ a simplicial $A$-module. 
We can form the trivial square zero extension $A \oplus M$ of $A$ by $M$. It is a simplical
ring whose underlying simplicial abelian group is $A \times M$, and for which the multiplication
is the usual one $(a,m).(a',m')=(a.a',am'+a'm)$ (this formula holds degreewise in the
simplicial direction). If we denote by $X=\Spec\, A$, then 
$\Spec\, (A\oplus M)$ will be denoted by $X[M]$, and is by definition the trivial square 
zero extension of $X$ by $M$. We note here that $M$ can also be considered through its
normalization as a $N(A)$-dg-module and thus as an object in $L_{\qcoh}(X)$
with zero positive cohomology sheaves.

This construction can be globalized as follows. For $X$ a derived scheme and $E$ an object 
in $L_{\qcoh}(X)$ whose cohomology is concentrated in nonpositive degrees, we can form
a derived scheme $X[M]$ as the relative spectrum $\Spec\, (\OO_{X} \oplus E)$. Locally, 
when $X=\Spec\, A$ is affine, $E$ corresponds to a simplicial $A$-module, and 
$X[M]$ simply is $\Spec\, (A\oplus M)$. 
The derived scheme $X[M]$ sits under the derived scheme $X$ itself and is considered
in the comma $\sz$-category $X/\dSch$ of derived schemes under $X$. The mapping space
$Map_{X/\dSch}(X[M],X)$
is called the space of derivations on $X$ with coefficients in $M$, and when $E=\OO_{X}$
this can be considered as the space of vector fields on $X$. It is possible to show the
existence of an object $\mathbb{L}_{X}$ together with a universal derivation
$X[\mathbb{L}_{X}] \longrightarrow X.$
The object $\mathbb{L}_{X}$ together with the universal derivation are characterized by the
following universal property
$$Map_{X/\dSch}(X[M],X) \simeq Map_{L_{\qcoh}(X)}(\mathbb{L}_{X},M).$$
The object $\mathbb{L}_{X}$ is called the \emph{absolute cotangent complex of $X$}. Its 
restriction on an affine open $\Spec\, A \subset X$ is a quasi-coherent
complex $L_{\qcoh}(\Spec\, A)$ which corresponds to the simplicial $A$-module 
$\mathbb{L}_{A}$ introduced in \cite{quil}. 

The absolute notion has a relative version for any morphism of derived schemes
$f : X \longrightarrow Y$. There is a natural morphism $f^{*}(\mathbb{L}_{Y}) \longrightarrow \mathbb{L}_{X}$ in $L_{\qcoh}(X)$, 
and the \emph{relative cotangent
complex of $f$} is defined to be its cofiber
$$\mathbb{L}_{X/Y}=\mathbb{L}_{f} := cofiber \left( f^{*}(\mathbb{L}_{Y}) \longrightarrow 
\mathbb{L}_{X} \right).$$
It is an object in $L_{\qcoh}(X)$, cohomologically concentrated in nonpositive degrees, 
and equiped with a universal derivation  
$X[\mathbb{L}_{X/Y}] \longrightarrow X$ which is now a morphism 
in the double comma $\sz$-category $X/\dSch/Y$.

One of the characteristic properties of derived schemes is that cotangent complexes are
compatible with fibered products, as opposed to what is happening in the case of schemes.
For any cartesian square of derived schemes
$$\xymatrix{
X' \ar[r]^-{g} \ar[d]_-{q} & X \ar[d]^-{p} \\
Y' \ar[r]_-{f} & Y}$$
the natural morphism 
$g^{*}(\mathbb{L}_{X/Y}) \longrightarrow \mathbb{L}_{X'/Y'}$
is an equivalence in $L_{\qcoh}(X')$. This property is true in the setting of schemes only 
under some Tor-independence conditions insuring the pull-back square of schemes remains
a pull-back in derived schemes (e.g.\ when one of the morphism $f$ or $p$ is
flat, see \cite{illu}). 

We will see later on how cotangent complexes can also be used in order to
understand how morphisms decompose along Postnikov towers and more generally 
how they control obstruction theories (see our $\S4.1$). Let us simply mention here that 
for any derived scheme $X$, the inclusion of $i : t_{0}(X) \longrightarrow X$
induces a morphism on cotangent complexes
$i^{*}(\mathbb{L}_{X}) \longrightarrow \mathbb{L}_{t_{0}(X)},$
which is always an obstruction theory on $t_{0}(X)$ in the sense of \cite{behrfant} (see \cite{schu}
for more details about the relation between derived scheme and the obstruction theories induced
on the truncations).

Finally, cotangent complexes can be used in order to define smooth and \'etale morphisms
between derived schemes. A morphism $f : X \longrightarrow Y$ in $\dSch$ will be
called \emph{\'etale} (resp. \emph{smooth}) if it is locally of finite presentation
(see \cite[\S2.2.2]{hagII} for the definition of finite presentation in the homotopical context) 
and if $\mathbb{L}_{f}$ vanishes (resp. $\mathbb{L}_{f}$ is a vector
bundle on $X$). An \'etale
(resp. smooth) morphism $f : X \longrightarrow Y$ of derived schemes induces an \'etale (resp. smooth) 
morphism on the truncations $t_{0}(f) : t_{0}(X) \longrightarrow t_{0}(Y)$, 
which is moreover flat: for all $i$ the natural morphism 
$t_{0}(f)^{*}(\pi_{i}(\OO_{Y})) \longrightarrow \pi_{i}(\OO_{X})$
is an isomorphism of quasi-coherent sheaves on $t_{0}(X)$ (see \cite[\S2.2.2]{hagII}). We easily deduce from this
the so-called Whitehead theorem for derived schemes: a morphism of derived
scheme $f : X \longrightarrow Y$ is an equivalence if and only if 
it induces an isomorphism on the truncation and if it is moreover smooth. 

\subsubsection*{Virtual classes} For a derived scheme $X$
the sheaves $\pi_{i}(\OO_{X})$ define quasi-coherent sheaves on 
the truncation $t_{0}(X)$. Under the condition that $t_{0}(X)$
is locally noetherian, and that $\pi_{i}(\OO_{X})$ are coherent 
and zero for $i\gg 0$, we find a well defined class in the K-theory of coherent sheaves on $t_{0}(X)$
$$[X]\Kvir := \displaystyle{\sum_{i}} (-1)^{i}[\pi_{i}(\OO_{X})] \in G_{0}(t_{0}(X)),$$
called the \emph{K-theoretical virtual fundamental class of $X$}.

 The  class $[X]\Kvir$ possesses another
interpretation which clarifies its nature. We keep assuming that 
$t_{0}(X)$ is locally noetherian and that $\pi_{i}(\OO_{X})$ are coherent and vanish 
for $i$ big enough. An object $E\in L_{\qcoh}(X)$ will be 
called \emph{coherent} if
it is cohomologically bounded and if for all $i$, $H^{i}(E)$ 
is a coherent sheaf on $t_{0}(X)$. The $\sz$-category of coherent 
objects in $D_{\qcoh}(X)$ form a thick triangulated subcategory and thus
can be used in order to define $G_{0}(X)$ as their Grothendieck group. The 
group $G_{0}(X)$ is functorial in $X$ for morphisms whose push-foward preserves coherent 
sheaves. This is in particular the case for the natural map $j : t_{0}(X) \longrightarrow X$,
and we thus have a natural morphism
$j_{*} : G_{0}(t_{0}(X))\longrightarrow G_{0}(X).$
By devissage this map is bijective, and we have
$$[X]\Kvir=(j_{*})^{-1}([\OO_{X}]).$$
In other words, $[X]\Kvir$ simply is the (nonvirtual) fundamental class
of $X$, considered as a class on $t_{0}(X)$ via the bijection above. This interpretation
explains a lot of things: as $j_{*}([X]\Kvir)=[\OO_{X}]$, integrating over $t_{0}(X)$ with
respect to the class $[X]\Kvir$ is equivalent to integrating over $X$. Therefore, 
whenever a numerical invariant is  obtained by integration over
a virtual class (typically a Gromov-Witten invariant), then it is 
actually an integral over some naturally defined derived scheme.

The virtual class in K-theory can also provide a virtual class in homology. Assume
for instance that $X$ is a derived scheme which is of finite presentation over a field $k$, 
and that $\mathbb{L}_{X/k}$ is perfect of amplitude $[-1,0]$ (i.e., is locally the cone 
of a morphism between two vector bundles). Then $t_{0}(X)$ is automatically
noetherian and $\pi_{i}(\OO_{X})$ are coherent and vanish 
for $i$ big enough (see \cite[SubLem. 2.3]{toen4}). Moreover, the inclusion $j : t_{0}(X) \longrightarrow X$
produces a perfect complex $j^{*}(\mathbb{L}_{X/k})$ whose dual
will be denoted by $\mathbb{T}^{vir}$ and called the virtual tangent sheaf. It has a
Todd class in Chow cohomology $Td(\mathbb{T}^{vir}) \in A^{*}(t_{0}(X))$ (see \cite{fult}). 
We can define the virtual class in Chow homology by the formula
$$[X]^{vir}:=\tau([X]\Kvir).Td(\mathbb{T}^{vir})^{-1} \in A_{*}(t_{0}(X)),$$
where $\tau : G_{0}(t_{0}(X)) \longrightarrow A_{*}(t_{0}(X))$ is the Grothendieck-Riemann-Roch
transformation of \cite[\S18]{fult}. We refer to \cite{ciockapr2,lowrschu} for more on the subject. 

Finally, K-theoretical virtual classes can be 
described for some of the basic examples of derived schemes mentioned in $\S2.2$. For
$Y \hookrightarrow X$ a local complete intersection closed immersion of locally noetherian 
schemes, the virtual class of $Y\times_{X}Y$ is given by
$$[Y\times_{X}Y]\Kvir=\lambda_{-1}(\mathcal{N}^{\vee})=\displaystyle{\sum_{i}}
(-1)^i[\wedge^{i}\mathcal{N}^{\vee}] \in 
G_{0}(Y),$$
where $\mathcal{N}$ is the normal bundle of $Y$ in $X$. In the same way, 
for $V$ a vector bundle on a locally noetherian scheme $X$, with a section 
$s$, the virtual class of $Eu(V,s)$ is the usual K-theoretic Euler class of $V$
$$[Eu(V,s)]\Kvir=\lambda_{-1}(V^{\vee}) \in G_{0}(Z(s)).$$
From this we get virtual classes of these two examples in Chow
homology, as being (localized) top Chern classes of the normal bundle
$\mathcal{N}$ and of $V$.

\subsubsection*{Relations with obstruction theories and dg-schemes} Let $X$ be a 
derived scheme of finite presentation over some base ring $k$. The inclusion 
$j : t_{0}(X) \longrightarrow X$ provides a morphism in $L_{\qcoh}(t_{0}(X))$
$$j^{*} : j^{*}(\mathbb{L}_{X/k}) \longrightarrow \mathbb{L}_{t_{0}(X)/k}.$$
This morphism is a perfect obstruction theory in the sense of \cite{behrfant}, and we get 
this way a forgetful $\sz$-functor from the $\sz$-category of 
derived schemes locally of finite presentation over $k$ to a certain $\sz$-category 
of schemes (locally of finite presentation over $k$) together with perfect obstruction theories. 
This forgetful $\sz$-functor is neither full nor faithful, but is conservative (as this
follows from the already mentioned Whitehead theorem, or from obstruction theory, 
see $\S4.1$). The essential surjectivity of 
this forgetful $\sz$-functor has been studied in \cite{schu}. The notion of derived scheme is 
strictly more structured than the notion of schemes with a perfect obstruction theory. The later
is enough for enumerative purposes, typically for defining virtual classes (as explained
above), but is not enough to recover finer invariants such as the quasi-coherent derived
category of derived schemes. 

The situation with dg-schemes in the sense of \cite{ciockapr} is opposite, there is
a forgetful functor from dg-schemes to derived schemes, which is neither
full, nor faithful nor essentially surjective, but is conservative. The notion 
of dg-schemes is thus strictly more structured than the notion of derived schemes. 
To be more precise, a dg-scheme is by definition essentially 
a pair $(X,Z)$, consisting of a derived scheme $X$, a scheme $Z$, and together  
closed immersion $X \longrightarrow Z$. Maps between dg-schemes are given by 
the obvious notion of maps between pairs. The forgetful $\sz$-functor simply sends
the pair $(X,Z)$ to $X$. Thus, dg-schemes can only model 
\emph{embeddable derived schemes}, that is derived schemes that can be
embedded in a scheme, and maps between dg-schemes can only model 
embeddable morphisms. In general, derived schemes and morphisms between derived schemes
are not embeddable, in the exact same way that formal schemes and morphisms between 
formal schemes are not so. This explains why the notion of a \emph{dg-scheme} is too strict 
for dealing with certain derived moduli problems and only sees a tiny part of the 
general theory of derived schemes. 

Finally, there is also a forgetful $\sz$-functor from the category of derived schemes (maybe of
characteristic zero) to the $2$-category of differential graded schemes 
of \cite{behr}. This $\sz$-functor is again neither full, nor faithful, nor
essentially surjective, but is conservative. The major reason comes
from the fact that differential graded schemes of \cite{behr} are defined
by gluing derived rings only up to $2$-homotopy (i.e.,  in 
the $2$-truncation of the $\sz$-category of derived rings), and thus misses the higher homotopical
phenomenon. 

\subsection{Derived moduli problems and derived schemes}

In the last paragraph we have seen the $\sz$-category $\dSch$ of derived schemes, some
basic examples as well as some characteristic properties. In order to introduce
more advanced examples we present here the functorial point of view and embed
the $\sz$-category $\dSch$ into the $\sz$-topos $\dSt$ of derived stacks (for the
\'etale topology). Objects in $\dSt$ are also called \emph{derived moduli problems}, 
and one major question is their representability. We will see some examples
(derived character varieties, derived Hilbert schemes and derived mapping spaces) of derived moduli problems
representable by derived schemes.  In order to consider more examples we will
introduce the notion of derived Artin stacks in our next paragraph $\S3.3$, which will enlarge considerably the
number of examples of representable derived moduli problems. \smallbreak

We let $\dAff$ be the $\sz$-category of affine derived schemes, which is also
equivalent to the opposite $\sz$-category of derived rings $\scomm$. The $\sz$-category
can be endowed with the \'etale topology: a family of morphisms between affine derived schemes
$\{U_{i} \rightarrow X\}$ is defined to be an \'etale covering if 
\begin{itemize}

\item each morphism $U_{i} \longrightarrow X$ is \'etale (i.e. of finite presentation 
and $\mathbb{L}_{U_{i}/X}\simeq 0$),

\item the induced morphism on truncations $\coprod_{i}t_{0}(U_{i}) \longrightarrow t_{0}(X)$
is a surjective morphism of schemes.

\end{itemize}

The \'etale covering families define a Grothendieck topology on $\dAff$ and 
we can thus form the $\sz$-category of stacks (see $\S2.1.2$). We denote it by
$$\dSt:=\dAff^{\sim,et}.$$
Recall from $\S2.1.2$ that the $\sz$-category 
$\dSt$ consists of the full sub-$\sz$-category of $\fun(\scomm,\Top)$
of $\sz$-functors satisfying the etale descent condition.
By definition, $\dSt$ is the $\sz$-category of derived stacks, and is the $\sz$-categorical 
version of the category of sheaves on the big \'etale site of affine schemes
(see e.g.\ \cite[\S1]{lm}for the big \'etale site of underived schemes). Its objects
are simply called \emph{derived stacks} or \emph{derived moduli problems}, and 
we will be interested in their representability by geometric objects such as 
derived schemes, or more generally derived Artin stacks. \smallbreak

We consider the Yoneda embedding 
$$h : \dSch \longrightarrow \dSt,$$
which sends a derived scheme $X$ to its $\sz$-functor of points 
$Map_{\dSch}(-,X)$ (restricted to affine derived schemes). The $\sz$-functor
$h$ is fully faithful, which follows from a derived version of fpqc descent for schemes
(see \cite{hagII,dag}). A derived moduli problem $F \in \dSt$ is then representable by a derived scheme
$X$ if it is equivalent to $h_{X}$. Here are below three examples of derived moduli problems
represented by derived schemes. 

\subsubsection*{Derived character varieties} We first describe a derived version of 
character varieties and character schemes, which are derived extensions of
the usual affine algebraic varieties (or schemes) of linear representations 
of a given group (see e.g.\ \cite{culeshal}). 

We fix an affine algebraic group scheme $G$ over some base field $k$. We let $\Gamma$
be a discrete group and we define a derived moduli problem $\Map(\Gamma,G)$,
of morphisms of groups from $\Gamma$ to $G$ as follows. The group scheme $G$ 
is considered as a derived group scheme using the inclusion $\Sch_{k} \hookrightarrow \dSch_{k}$, 
of schemes over $k$ to derived schemes over $k$. The group object $G$ defines 
an $\sz$-functor
$$G : \dAff_{k}^{op} \longrightarrow \Top\mhyphen Gp,$$
from affine derived schemes to the $\sz$-category of group objects in $\Top$, or
equivalently the $\sz$-category of simplicial groups. 
We define $\Map(\Gamma,G) : \dAff_{k}^{op} \longrightarrow \Top$ by sending $S\in \dAff_{k}^{op}$
to $Map_{\Top\mhyphen Gp}(\Gamma,G(S))$. 

The derived moduli problem $\Map(\Gamma,G)$ is representable by an affine derived scheme. 
This can be seen as follows. When $\Gamma$ is free, then $\Map(\Gamma,G)$ is
a (maybe infinite) power of $G$, and thus is an affine scheme. In general, we can write $\Gamma$ has the colimit 
in $\Top\mhyphen Gp$ of free groups, by taking for instance a simplicial free resolution. Then
$\Map(\Gamma,G)$ becomes a limit of affine derived scheme and thus is itself an affine
derived scheme. When the group $\Gamma$ has a simple presentation by generators and relations
the derived affine scheme $\Map(\Gamma,G)$ can be described explicitly by means of simple
fibered products. A typical example appears when $\Gamma$ is fundamental group of 
a compact Riemann surface of genus $g$: the derived affine scheme 
$\Map(\Gamma,G)$ comes in a cartesian square
$$\xymatrix{
\Map(\Gamma,G) \ar[r] \ar[d] & G^{2g} \ar[d] \\
Spec\, k \ar[r] & G,}$$
where the right vertical map sends $(x_{1},\dots,x_{g},y_{1},\dots,y_{g})$
to the product of commutators $\prod_{i}[x_{i},y_{i}]$. 

The tangent complex of the derived affine scheme $\Map(\Gamma,G)$ can be
described as the group cohomology of $G$ with coefficients in the universal
representation $\rho : \Gamma \longrightarrow G(\Map(\Gamma,G))$. The 
morphism $\rho$ defines an action of $G$ on the trivial principal $G$-bundle 
on $\Map(\Gamma,G)$, and thus on the vector bundle $V$ associated to the
adjoint action of $G$ on its Lie algebra $\frak{g}$. The cochain complex
of cohomology of $\Gamma$ with coefficients in $V$ provides a 
quasi-coherent complex $C^{*}(\Gamma,V)$  on $\Map(\Gamma,G)$.  The tangent complex is then 
given by the part sitting in degrees $[1,\infty[$ as follows
$$\mathbb{T}_{\Map(\Gamma,G)} \simeq C^{\geq 1}(\Gamma,V)[1].$$
The algebraic group $G$ acts on $\Map(\Gamma,G)$, and when $G$ is linearly reductive
we can consider the derived ring of invariant functions $\OO(\Map(\Gamma,G))^{G}$. The
spectrum of this derived ring $\Spec\, \OO(\Map(\Gamma,G))^{G}$ is a derived 
GIT quotient of the action of $G$ on $\Map(\Gamma,G)$ and deserves the name
of \emph{derived character variety of $\Gamma$ with coefficients in $G$}. 

It is interesting to note here that the above construction can be modified in a
meaningful manner. We assume that $\Gamma$ is the fundamental group
of a connected CW complex $X$. We can modify the derived moduli problem
of representations of $\Gamma$ by now considering rigidified local systems on the
space $X$. In the underived setting these two moduli problems are equivalent, but it is
one interesting feature of derived algebraic geometry to distinguish them. We define
$\mathbb{R}\uLocs(X,G)$ as follows. We chose a simplicial group
$\Gamma_{*}$ with a weak equivalence $X\simeq B\Gamma_{*}$, that is 
$\Gamma_{*}$ is a simplicial model for the group of based loops in $X$. We then 
define $\Map(\Gamma_{*},G)$ by sending a derived scheme $S$ to 
$Map_{\Top\mhyphen Gp}(\Gamma_{*},G(S))$. This new derived moduli problem is again 
representable by a derived affine scheme $\mathbb{R}\uLocs(X,G)$. The truncations
of $\mathbb{R}\uLocs(X,G)$ and $\Map(\Gamma,G)$ are both equivalent to the
usual affine scheme of maps from $\Gamma$ to $G$, but the derived structures 
differ. This can be seen at the level of tangent complexes. As for $\Map(\Gamma,G)$
then tangent complex of $\mathbb{R}\uLocs(X,G)$ is given by
$$\mathbb{T}_{\mathbb{R}\uLocs(X,G)} \simeq C^{*}_{x}(X,V)[1],$$
where now $V$ is considered as a local system of coefficients on $X$ and we consider
the cochain complex of cohomology of $X$ with coefficients in $V$, and 
$C^{*}_{x}(X,V)$ denotes the reduced cohomology with respect to the base point $x$ of $X$
(the fiber of $C^{*}(X,V) \longrightarrow C^{*}(\{x\},V)\simeq \frak{g}$).
Interesting examples are already obtained with $\Gamma=*$ and $X$ higher dimensional spheres.
For $X=S^{n}$, $n>1$, and $k$ of caracteristic zero, we have
$$\mathbb{R}\uLocs(X,G)\simeq 
\Spec\, Sym_{k}(\frak{g}^{*}[n-1]).$$

\subsubsection*{The derived scheme of maps} We let $k$ be a commutative ring and 
$X$ be a scheme which is projective and flat over $Spec\, k$, and 
$Y$ a quasi-projective scheme over $Spec\, k$. We consider the 
derived moduli problem of maps of derived $k$-schemes from $X$ to $Y$, 
which sends $S\in \dSch_{k}$ to $Map_{\dSch_{k}}(X\times S,Y)$. This
is a derived stack (over $k$) $\Map_{k}(X,Y) \in \dSt_{k}$, which can be shown to be 
representable by a derived scheme $\Map_{k}(X,Y)$ which is locally of finite
presentation of $Spec\, k$ (see corollary \ref{c1} for a more general version). 
The truncation 
$t_{0}(\Map_k(X,Y))$ is the usual scheme of maps from $X$ to $Y$ as originally constructed by Grothendieck. 
Except in some very specific cases the derived scheme $\Map_k(X,Y)$ is not a scheme. This
can be seen at the level of tangent complexes already, as we have the following
formula for the tangent complex of the derived moduli space of maps
$$\mathbb{T}_{\Map_k(X,Y)} \simeq \pi_{*}(ev^{*}(\mathbb{T}_{Y})) \in L_{\qcoh}(\Map_k(X,Y)),$$
where $ev : \Map_k(X,Y)\times X \longrightarrow Y$ is the evaluation morphism, and 
$\pi : \Map_k(X,Y)\times X \longrightarrow \Map_k(X,Y)$ is the projection morphism. 
This formula shows that when $Y$ is for instance smooth, then 
$\mathbb{T}_{\Map_k(X,Y)}$ is perfect of amplitude contained in $[0,d]$ where
$d$ is the relative dimension of $X$ over $Spec\, k$. When this amplitude is 
actually strictly bigger than $[0,1]$, the main result of \cite{avra} implies that 
$\Map_k(X,Y)$ cannot be an (underived) scheme. 

One consequence of the representability of $\Map_k(X,X)$ is the
representability of the derived group of automorphisms of $X$, 
$\mathbb{R}\mathbf{Aut}_k(X)$, which is the open derived subscheme
of $\Map_k(X,X)$ consisting of automorphisms of $X$. The derived scheme
$\mathbb{R}\mathbf{Aut}_k(X)$ is an example of a derived group scheme
locally of finite presentation over $Spec\, k$. Its tangent complex at the
unit section is the complex of globally defined vector fields on $X$ over $Spec\, k$, 
$\mathbb{H}(X,\mathbb{T}_{X})$. We will see in $\S5.4$ that the complex $\mathbb{H}(X,\mathbb{T}_{X})$
always comes equipped with a structure of a dg-lie algebra, at least up to 
an equivalence, and this dg-lie algebra is here the tangent Lie algebra of the
derived group scheme $\mathbb{R}\mathbf{Aut}_k(X)$. For the same reasons as above, 
invoking \cite{avra}, the derived group scheme $\mathbb{R}\mathbf{Aut}_k(X)$ is in general
not a group scheme.  

\subsubsection*{Derived Hilbert schemes} We let again $X$ be a projective and flat 
scheme over $Spec\, k$ for some commutative ring $k$. For sake of simplicity we will only
be interested in a nice part of the derived Hilbert scheme of $X$, corresponding
to closed subschemes which are of local complete intersection (we refer
to \cite{ciockapr} for a more general construction). For any $S \in \dSch_{k}$ 
we consider the $\sz$-category $\dSch_{(X\times S)}$ of derived schemes
over $X\times S$. We let $\mathbb{R}\mathbf{Hilb}^{lci}(X)(S)$ 
be the (nonfull) sub-$\sz$-category of $\dSch_{(X\times S)}$ defined as follows.
\begin{itemize}

\item The objects of $\mathbb{R}\mathbf{Hilb}(X)(S)$ are the derived schemes 
$Z \longrightarrow X\times S$ which are flat over $S$, finitely presented
over $X\times S$, and moreover
induce a closed immersion on the truncation $t_{0}(Z) \hookrightarrow X \times t_{0}(S)$.

\item The morphisms are the equivalences in the $\sz$-category $\dSch_{(X\times S)}$. 

\end{itemize}

For a morphism of derived schemes $S' \rightarrow S$, the pull-back induces a morphism
of $\sz$-groupoids
$$\mathbb{R}\mathbf{Hilb}^{lci}(X)(S) \longrightarrow \mathbb{R}\mathbf{Hilb}^{lci}(X)(S').$$
This defines an $\sz$-functor from $\dSch_{k}^{op}$ to the $\sz$-category of 
$\sz$-groupoids which we can compose with the
nerve construction to get an $\sz$-functor from $\dSch_{k}^{op}$ to $\Top$, and thus
a derived moduli problem. This derived moduli problem is representable by 
a derived scheme $\mathbb{R}\mathbf{Hilb}^{lci}(X)$ which is locally of finite presentation
over $Spec\, k$. Its truncation
is the open subscheme of the usual Hilbert scheme of $X$ (over $k$) corresponding
to closed subschemes which are embedded in $X$ as local complete intersections. 

The tangent complex of $\mathbb{R}\mathbf{Hilb}^{lci}(X)$ can be described as follows. 
There is a universal closed derived subscheme
$j : \mathcal{Z} \longrightarrow X \times \mathbb{R}\mathbf{Hilb}^{lci}(X),$
with a relative tangent complex
$\mathbb{T}_{j}$ which consists of a vector bundle concentrated in degree $1$
(the vector bundle is the normal bundle of the inclusion $j$). If we denote by
$p : \mathcal{Z} \longrightarrow \mathbb{R}\mathbf{Hilb}^{lci}(X)$ the flat projection, we have
$$\mathbb{T}_{\mathbb{R}\mathbf{Hilb}^{lci}(X)}\simeq 
p_{*}(\mathbb{T}_{j}[1]).$$

In the same way, there exists a derived Quot scheme representing a derived version of
the Quot functor. We refer to \cite{ciockapr} for more on the subject. 

\subsection{Derived moduli problems and derived Artin stacks}

It is a fact of life that many interesting moduli problems are not representable by schemes, 
and algebraic stacks  have been introduced in order to extend the notion of representability
(see \cite{grot4,delimumf,arti,lm}). 
This remains so in the derived setting: many derived moduli problems are not representable
by derived schemes and it is necessary to introduce more general objects called
\emph{derived Artin stacks} in order to overcome this issue. \smallbreak

In the last paragraph we have embedded the $\sz$-category of derived schemes
$\dSch$ into the bigger $\sz$-category of derived stacks $\dSt$. We will
now introduce an intermediate $\sz$-category $\dSt^{Ar}$
$$\dSch \subset \dSt^{Ar} \subset \dSt,$$
which is somehow the closure of $\dSch$ by means of taking quotient by smooth groupoid objects
(see \cite[4.3.1]{lm} for the notion of groupoid objects in schemes in the nonderived setting). 

A groupoid object $\dSt$ (also called a \emph{Segal groupoid}), 
consists of an $\sz$-functor
$$X_{*} : \Delta^{op} \longrightarrow \dSt$$
satisfying the two conditions below.
\begin{enumerate}

\item For all  $n$ the Segal morphism
$$X_{n} \longrightarrow X_{1} \times_{X_{0}}X_{1} \times_{X_{0}} \dots
\times_{X_{0}}X_{1}$$
is an equivalence of derived stacks.

\item The composition morphism
$$X_{2} \longrightarrow X_{1}\times_{X_{0}}X_{1}$$
is an equivalence of derived stacks.
\end{enumerate}

In the definition above the object $X_{0} \in \dSt$ is the derived stack of objects
of the groupoid $X_{*}$, and $X_{1}$ the derived stack of morphisms. The morphism in the first 
condition above is induced by maps $[1] \rightarrow [n]$ in $\Delta$ sending 
$0$ to $i$ and $1$ to $i+1$. It provides the composition in the grupoid by means of the
following diagram
$$X_{1} \times_{X_{0}}X_{1} \simeq X_{2} \longrightarrow X_{1}$$
induced by the morphism $[1] \rightarrow [2]$ sending $0$ to $0$ and $1$ to $2$. The
morphism of the second condition insures that this composition is invertible up to 
an equivalence. We refer the reader to 
\cite[Def. 4.9.1]{hagI} \cite[\S1.3.4]{hagII} for more about Segal cateories and Segal groupoids
objects.

We say that a groupoid object $X_{*}$ is a \emph{smooth groupoid of derived schemes} if $X_{0}$ and $X_{1}$
are derived schemes and if the projections $X_{1} \longrightarrow X_{0}$ are smooth morphisms
of derived schemes. The colimit of the
simplicial object $X_*$ is denoted by $|X_{*}| \in \dSt$  and is called the \emph{quotient
derived stack} of the groupoid $X_{*}$.

\begin{df}\label{d4}
\begin{enumerate}
\item
A derived stack is a \emph{derived $1$-Artin stack} if it is of the form $|X_{*}|$ for
some smooth groupoid of derived schemes $X_{*}$. 
\item A morphism $f : X \longrightarrow Y$ between derived $1$-Artin stacks  
is \emph{smooth} if there exist smooth groupoids of derived schemes $X_{*}$ and $Y_{*}$ and 
a morphism of groupoid objects $f_{*} : X_{*} \longrightarrow Y_{*}$
with $f_{0} : X_{0} \longrightarrow Y_{0}$ smooth, such that $|f_{*}|$ is equivalent to $f$.
\end{enumerate}
\end{df}

The derived $1$-Artin stack form a full sub-$\sz$-category $\dSt^{Ar,1} \subset \dSt$ which
contains derived schemes (the quotient of the constant groupoid associated to a derived
scheme  $X$ gives back $X$). Moreover the definition above provides a notion of 
smooth morphisms between derived $1$-Artin stacks and the definition can thus be
extended by an obvious induction. 

\begin{df}\label{d4'}
\begin{enumerate}
\item
A derived stack is a \emph{derived  $n$-Artin stack} if it is of the form $|X_{*}|$ for
some smooth groupoid of derived $(n-1)$-Artin stacks $X_{*}$. 
\item A morphism between derived  $n$-Artin stacks $f : X \longrightarrow Y$ 
is \emph{smooth} if there exists smooth groupoid of  derived $(n-1)$-Artin stacks $X_{*}$ and $Y_{*}$, 
a morphism of groupoid objects $f_{*} : X_{*} \longrightarrow Y_{*}$
with $f_{0} : X_{0} \longrightarrow Y_{0}$ smooth, 
and such that $|f_{*}|$ is equivalent to $f$.
\end{enumerate}
A derived stack is a \emph{derived Artin stack} if it is a derived  $n$-Artin 
stack for some  $n$. The full sub-$\sz$-category of derived Artin stacks is denoted by $\dSt^{Ar}$.
\end{df}

Here are some standard examples of derived Artin stacks. More involved examples will be given 
later after having stated the representability theorem \ref{tartluri}. 

\subsubsection*{Quotients stacks} Let $G$ be a smooth group scheme over some base derived scheme
$S$. We assume that $G$ acts on a derived scheme $X \rightarrow S$. We can form
the quotient groupoid $B(X,G)$, which is the simplicial object  that is equal to 
$X\times_{S}G^{n}$ in degree  $n$, and with the faces and degeneracies 
 defined in the usual manner by means of the action of $G$ on $X$ and the multiplication in $G$. The 
groupoid $B(X,G)$ is a smooth groupoid of derived schemes over $S$, and its
quotient stack $|B(X,G)|$ is thus an example of a derived Artin stack which is denoted by 
$[X/G]$. It is possible to prove
that for a derived scheme $S' \longrightarrow S$, the simpicial set 
$Map_{\dSt/S}(S',[X/G])$ is (equivalent to) the nerve of the 
$\sz$-groupoid of diagrams of derived stacks over $S$ endowed with $G$-actions
$$\xymatrix{
P \ar[r] \ar[d] & X \ar[d] \\
S' \ar[r] & S,}$$
where the induced morphism $[P/G] \longrightarrow S'$ is moreover an equivalence 
(i.e., $P \longrightarrow S'$ is a principal $G$-bundle). 

In the example of the derived character scheme given in $\S3.2$, the group $G$ acts on
$\Map(\Gamma,G)$, and the quotient stack $[\Map(\Gamma,G)/G]$ is now the
derived Artin stack of representations of $\Gamma$  with coefficients in $G$ up to
equivalence. In the same way, 
$\mathbb{R}\underline{Loc}(X,G):=[\mathbb{R}\uLocs(X,G)/G]$
becomes the derived Artin stack of $G$-local systems on the topological space $X$ without 
trivialization at the base point. For higher dimensional spheres and $k$ a field of 
characteristic zero, we get an explicit presentation
$$\mathbb{R}\underline{Loc}(S^{n},G)\simeq [\Spec\, A/G],$$ where
$A=Sym_{k}(\frak{g}^{*}[n-1])$ and $G$ acts on $A$ by its co-adjoint
representation.

\subsubsection*{Eilenberg-MacLane and linear derived stacks} If $G$ is a smooth derived group scheme over 
some base derived scheme $S$, we have a classifying stack $BG:=[S/G]$
over $S$. When $G$ is abelian the derived Artin stack $BG$ is again a
smooth abelian group object in derived stacks. We can therefore
iterate the construction and set $K(G,n):=B(K(G,n-1))$,
$K(G,0)=BG$. The derived stack $K(G,n)$ is an example of a derived
$n$-Artin stack smooth over $S$. For each scheme $S' \longrightarrow
S$ we have
$$\pi_{i}(Map_{\dSt/S}(S',K(G,n))) \simeq H^{n-i}_{et}(S',G_{S'}),$$
where $G_{S'}$ is the sheaf of abelian groups represented by $G$ on the small \'etale site
of the derived scheme $S'$.

It is also possible to define
$K(G,n)$ with $n<0$ by the formula $K(G,n):=S\times_{K(G,n+1)}S$. With these notations, 
we do have $S\times_{K(G,n)}S\simeq K(G,n-1)$ for all $n\in \mathbb{Z}$, and these are all 
derived group schemes over $S$, smooth for $n\geq 0$. However, $K(G,n)$ are in general
not smooth for $n<0$. In the special case where $G$ is affine and smooth over  
over $S$ a scheme of characteristic zero, $K(G,n)$ can be described as a relative spectrum
$K(G,n) \simeq \Spec\, Sym_{\OO_{S}}(\frak{g}^{*}[-n])$
for $n<0$ and $\frak{g}$ is the Lie algebra of $G$ over $S$. 

A variation of the notion of Eilenberg-MacLane derived stack is the notion of
linear stack associated to perfect complexes. We let $S$ be a derived
scheme and $E \in L_{\qcoh}(S)$ be a quasi-coherent complex on $S$. We define a derived stack 
$\mathbb{V}(E)$
over $S$ as follows. For $u : S' \longrightarrow S$ a derived scheme we let 
$\mathbb{V}(E)(S'):=Map_{L_{\qcoh}(S')}(u^{*}(E),\OO_{S'})$. This defines an $\sz$-functor
$\mathbb{V}(E)$ on the $\sz$-category of derived schemes over $S$ and thus an 
object in $\dSt/S$, the $\sz$-category of derived stacks over $S$. 
The derived stack $\mathbb{V}(E)$ is a derived Artin stack over $S$ as soon 
as $E$ is a perfect $\OO_{S}$-module (i.e., is locally for the Zariski topology on 
$S$ a compact object in the quasi-coherent derived category, see \cite[\S2.4]{toenvaqu}
and \cite[SubLem. 3.9]{toenvaqu}
for more on perfect objects and derived Artin stacks). More is true, 
we can pull-back the relative 
tangent complex of $\mathbb{V}(E)$ along the zero section $e : S \longrightarrow
\mathbb{V}(E)$, and get (see $\S4.1$ for (co)tangent complexes of derived Artin stack)
$e^{*}(\mathbb{T}_{\mathbb{V}(E)/S}) \simeq E^{\vee},$
where $E^{\vee}$ is the dual of $E$. We thus see that $\mathbb{V}(E) \longrightarrow S$ is
smooth if and only if $E^{\vee}$ is of amplitude contained in $[-\infty,0]$, that is
if and only if $E$ is of nonnegative Tor amplitude (see \cite[\S2.4]{toenvaqu}, 
or below for the notion of amplitude). On the other  hand, 
$\mathbb{V}(E)$ is a derived scheme if and only if $E$ is of nonpositive Tor amplitude, 
in which case it can be written as a relative spectrum
$\mathbb{V}(E) \simeq \Spec\, (Sym_{\OO_{S}}(E))$. We note that when $E$ is $\OO_{S}[n]$, then 
$\mathbb{V}(E)$ simply is $K(\mathbb{G}_{a,S},-n)$, where $\mathbb{G}_{a,S}$ is the
additive group scheme over $S$. In general, the derived stack 
$\mathbb{V}(E)$ is obtained by taking twisted forms and certain 
finite limits of derived stacks of the form
$K(\mathbb{G}_{a,S},-n)$.

\subsubsection*{Perfect complexes} We present here a more advanced and less trivial 
example of a derived Artin stack. For this we fix two integers $a\leq b$ and 
we define a derived stack $\Parf^{[a,b]} \in \dSt$, classifying perfect complexes of amplitude 
contained in $[a,b]$. As an $\sz$-functor it sends a derived scheme $S$ 
to the $\sz$-groupoid (consider as a simplicial set by the 
nerve construction, see $\S2.1.2$) of perfect objects in $L_{\qcoh}(S)$ with amplitude contained
in $[a,b]$. We remind here that the amplitude of a perfect complex $E$ on $S$ is
contained in $[a,b]$ if its cohomology sheaves are universally concentrated in 
degree $[a,b]$:  for every derived scheme $S'$ and every morphism 
$u : S' \rightarrow S$, 
we have $H^{i}(u^{*}(E))=0$ for $i\notin [a,b]$ (this can be tested for all $S'=Spec\, L$ with 
$L$ a field). The following theorem has been announced in \cite{hirssimp}, at least in the nonderived setting, 
and has been proved in \cite{toenvaqu}.

\begin{theorem}\label{t1}
The derived stack $\Parf^{[a,b]}$ is a derived Artin stack locally of finite presentation
over $\Spec\, \mathbb{Z}$. 
\end{theorem}

There is also a derived stack $\Parf$, classifying all perfect complexes, without any
restriction on the amplitude.
The derived stack $\Parf$ is covered by open derived substacks $\Parf^{[a,b]}$, 
and is itself an increasing union of open derived Artin substacks. Such derived 
stacks are called \emph{locally geometric} in \cite{toenvaqu} but we will allow ourselves to 
keep using the expression derived Artin stack.

The derived stack $\Parf$ is one of the most foundamental example of derived Artin stacks. First of
all it is a far-reaching generalization of the varieties of complexes 
(the so-called \emph{Buchsbaum-Eisenbud varieties}, see e.g.,\ \cite{decostri}). 
Indeed, the variety of complexes, suitably derived, can be shown to produce a smooth atlas
for the derived stack $\Parf$. In other words, $\Parf$ is the quotient of the (derived) varieties 
of complexes by the subtle equivalence relation identifying two complexes which are quasi-isomorphic. 
The fact that this equivalence relation involves dividing out by 
quasi-isomorphisms instead of isomorphisms
is responsible for the fact that $\Parf$ is 
only a derived Artin stack in a higher sense. To be more precise, $\Parf^{[a,b]}$ 
is a derived  $n$-Artin stack where $n=b-a+1$. This reflects the fact that morphisms
between complexes of amplitude in $[a,b]$ has homotopies and higher homotopies up to 
degree $n-1$, or equivalently that the $\sz$-category of complexes of amplitude
in $[a,b]$ has $(n-1)$-truncated mapping spaces.

The derived Artin stack $\Parf$ also possesses some extension, for instance 
by considering perfect complexes with an action of some nice dg-algebra, 
or perfect complexes over a given smooth and proper scheme. We refer 
to \cite{toenvaqu} in which the reader will find more details. 

\subsubsection*{Derived stacks of stable maps} Let $X$ be a smooth and projective
scheme over the complex numbers. We fix $\beta \in H_{2}(X(\mathbb{C}),\mathbb{Z})$ a curve class.
We consider $\overline{M}^{pre}_{g,n}$ the Artin stack of pre-stable curves
of genus $g$ an  $n$ marked points. It can be considered as a derived Artin stack and
thus as an object in $\dSt$. We let 
$\mathcal{C}_{g,n} \longrightarrow \overline{M}^{pre}_{g,n}$ the universal pre-stable curve. 
We let 
$$\mathbb{R}\overline{M}^{pre}_{g,n}(X,\beta)=\Map_{\dSt/\overline{M}^{pre}_{g,n}}(\mathcal{C}
_{g,n},X),$$
be the relative derived mapping stack of $\mathcal{C}_{g,n}$ to $X$ (with fixed
class $\beta$). The derived stack $\mathbb{R}\overline{M}^{pre}_{g,n}(X,\beta)$
is a derived Artin stack, as this can be deduced from the representability of
the derived mapping scheme (see $\S3.2$). It contains an open derived Deligne-Mumford 
substack $\mathbb{R}\overline{M}_{g,n}(X,\beta)$ which consists of 
stable maps. The derived stack $\mathbb{R}\overline{M}_{g,n}(X,\beta)$ is proper 
and locally of finite presentation over $Spec\, \mathbb{C}$, and can be used in order to 
recover Gromov-Witten invariants of $X$. We refer to \cite{reduce} for some works in that 
direction, as well as \cite{toen5} for some possible application to the categorification
of Gromov-Witten theory. \smallbreak

We now state the representability theorem of Lurie, an extremely powerful tool 
for proving that a given derived stack is a derived Artin stack, thus extending
to the derived setting the famous
Artin's representability theorem. The first proof appeared 
in the thesis \cite{luri3} and can now be found in the series \cite{dag}. There also are some
variations in \cite{prid} and \cite[App.~C]{hagII}.

\begin{theorem}\label{tartluri}
Let $k$ be a noetherian commutative ring.
A derived stack $F \in \dSt_{k}$ is a derived Artin stack locally of finite presentation
over $Spec\, k$ if and only if the following conditions are satisfied.
\begin{enumerate}

\item There is an integer $n\geq 0$ such that 
for any underived affine scheme $S$ over $k$ the simplicial set $F(S)$
is  $n$-truncated.

\item For any filtered system of derived $k$-algebras $A=colim A_{\alpha}$ the natural morphism
$$\displaystyle{colim_{\alpha}} F(A_{\alpha}) \longrightarrow F(A)$$
is an equivalence (where $F(A)$ means $F(\Spec\, A)$). 

\item For any derived $k$-algebra $A$ with Postnikov tower 
$$
\xymatrix{A \ar[r] & \dots \ar[r] &  A_{\leq k} \ar[r] & A_{\leq k-1} \ar[r] & \dots \ar[r] & \pi_{0}(A)},
$$
the natural morphism
$$F(A) \longrightarrow \displaystyle{lim}_{k} F(A_{\leq k})$$
is an equivalence.

\item The derived stack $F$ has an obstruction theory (see \cite[\S1.4.2]{hagII} for details).

\item For any local noetherian $k$-algebra $A$ with maximal ideal $m\subset A$, the natural morphism
$$F(\hat{A}) \longrightarrow \displaystyle{lim}_{k} F(A/m^{k})$$
is an equivalence (where $\hat{A}=lim A/m^k$ is the completion of $A$).

\end{enumerate}
\end{theorem}

We extract one important corollary of the above theorem.

\begin{cor}\label{c1}
Let $X$ be a flat and proper scheme over some base scheme $S$ and $F$ be a derived Artin stack
which is locally of finite presentation over $S$. Then the derived mapping stack
$\Map_{\dSt/S}(X,F)$ is again a derived Artin stack locally of finite presentation over $S$. 
\end{cor}

\subsection{Derived geometry in other contexts}

The formalism of derived schemes and derived Artin stacks we have described in this section
admits several modifications and generalizations that are worth mentioning. 

\subsubsection*{Characteristic zero} When restricted to 
zero characteristic, derived algebraic geometry admits a slight conceptual simplification
due to the fact that the homotopy theory of simplicial commutative $\mathbb{Q}$-algebras
become equivalent to the homotopy theory of nonpositively graded commutative
dg-algebras over $\mathbb{Q}$. This fact can be promoted to an equivalence of $\sz$-categories
$$N : \scomm_{\mathbb{Q}} \simeq \ncdga_{\mathbb{Q}},$$
induced by the normalization functor $N$. The normalization functor $N$ from
simplicial abelian groups to cochain complexes sitting in nonpositive degrees 
has a lax symmetric monoidal structure given by the so-called Alexander-Whithney morphisms
(see \cite{shipschw}), and thus always induces a well defined $\sz$-functor
$$N : \scomm \longrightarrow \ncdga_{\mathbb{Z}}.$$
This $\sz$-functor is not an equivalence in general but induces an equivalence on the
full sub-$\sz$-categories of $\mathbb{Q}$-algebras. 

The main consequence is that the notion of a derived scheme, and more generally of 
a derived Artin stack can, when restricted over $Spec\, \mathbb{Q}$,  also be modelled
using commutative dg-algebras instead of simplicial commutative rings: we can formally
replace $\scomm$ by $\ncdga_{\mathbb{Q}}$ in all the definitions and all the constructions, 
and  obtain a theory of derived $\mathbb{Q}$-schemes  and derived
Artin stacks over $\mathbb{Q}$ equivalent to the one
we have already seen. This simplifies a bit the algebraic manipulation at the level of derived
rings. For instance, the free commutative dg-algebras are easier to understand than
free simplicial commutative algebras, as the latter involves divided powers (see e.g.\ \cite{fres}). 
One direct consequence is that explicit computations involving generators
and relations tend to be easier done in the dg-algebra setting.
A related phenomenon concerns explicit models, using model category of commutative dg-algebras
for example. The 
cofibrant commutative dg-algebras are, up to a retract, the quasi-free
commutative dg-algebras (i.e., free as graded, non-dg, algebras), and are easier to understand than 
their simplicial counterparts. 

The reader can see derived algebraic geometry using dg-algebras in action for instance in \cite{bravjoyc,bogr,ptvv}. 

\subsubsection*{$E_{\infty}$-Algebraic geometry} The theory of derived rings $\scomm$ can 
be slightly modified by using other homotopical notions of the notion of rings. One possibility,
which have been explored in \cite{dag,hagII}, is to use the $\sz$-category of $E_{\infty}$-algebras
(or equivalently $H\mathbb{Z}$-algebras)
instead of $\sz$-category of simplicial commutative rings. The $\sz$-category $E_{\infty}\mhyphen \mathbf{dga}^{\leq 0}$, 
of nonpositively graded $E_{\infty}$-dg-algebras (over $\mathbb{Z})$, behaves formally
very similarly to the $\sz$-category $\scomm$. It contains the category of commutative rings
as a full sub-$\sz$-category of $0$-truncated objects, and more generally 
a given $E_{\infty}$-dg-algebra has a Postnikov tower as for the case of commutative simplicial rings, 
whose stage are also controlled by a cotangent complex. Finally, the normalization functor induces
an $\sz$-functor
$$N : \scomm \longrightarrow E_{\infty}\mhyphen \mathbf{dga}^{\leq 0}.$$
The $\sz$-functor is not an equivalence, expect when restricted to $\mathbb{Q}$-algebras
again. Its main failure of being an equivalence is reflected in the fact that it
does not preserve cotangent complexes in general. To present things differently, 
simplicial commutative rings and $E_{\infty}$-dg-algebras are both generated by the same
elementary pieces, namely commutative rings, but the manner these pieces are glued together
differs (this is typically what is happening in the Postnikov towers). 

As a consequence, there is a very well established algebraic geometry over $E_{\infty}$-dg-algebras, 
which is also a natural extension of algebraic geometry to the homotopical setting, but it 
differs from the derived algebraic geometry we have presented. The main difference
between the two theories can be found in the notion of smoothness: the affine line
over $Spec\, \mathbb{Z}$ is smooth as a derived scheme but it is not smooth as
an $E_{\infty}$-scheme (simply because the polynomial ring $\mathbb{Z}[T]$ differs
from the free $E_{\infty}$-ring on one generator in degree $0$, the latter involving homology of symmetric
groups has nontrivial cohomology). Another major difference 
is that derived algebraic geometry is the universal
derived geometry generated by algebraic geometry (this sentence can be
made into a mathematical theorem, expressing a universal property of $\dSch$), whereas $E_{\infty}$-algebraic
geometry is not. From a general point of view, $E_{\infty}$-algebraic geometry is
more suited to treat questions and problems of topological origin and derived algebraic geometry
is better suited to deal with questions coming from algebraic geometry. 

\subsubsection*{Spectral geometry} Spectral geometry is another modification of derived algebraic 
geometry. It is very close to $E_{\infty}$-algebraic geometry briefly mentioned above, and
in fact the $E_{\infty}$ theory is a special case of the spectral theory. This time it 
consists of replacing the $\sz$-category $\scomm$ with 
$\spcomm$, the $\sz$-category of commutative ring spectra. This is a generalization 
of $E_{\infty}$-algebra geometry, which is recovered as spectral schemes
over $\Spec\, H\mathbb{Z}$, where $H\mathbb{Z}$ is the Eilenberg-McLane ring spectrum. 
Spectral geometry is mainly 
developed in \cite{dag} (see also \cite[\S2.4]{hagII}), and has found an 
impressive application to the study of topological modular forms
(see \cite{luri5}).

\subsubsection*{Homotopical algebraic geometry} Homotopical algebraic geometry
is the general form of derived algebraic geometry, $E_{\infty}$-algebraic geometry
and spectral geometry. It is a homotopical version of relative geometry of \cite{haki},
for which affine schemes are in one-to-one correspondence with 
commutative monoids in a base symmetric monoidal model category (or more generally
a symmetric monoidal $\sz$-category). Most of the basic notions, such 
schemes, the Zariski, etale or flat topology, Artin stacks \dots have versions in this
general setting. This point of view is developed in \cite{hagII} as well as in \cite{souz}, 
and makes it possible to do geometry in non-additive contexts. Recently, this has lead to 
a theory of derived logarithmic geometry as exposed in \cite{schuvezz}. 

\subsubsection*{Derived analytic geometries} Finally, let us also mention the existence
of analytic counter-parts of derived algebraic geometry, but which 
are out of the scope of this paper. We refer to \cite{dag} in which 
derived complex analytic geometry is discussed. 

\section{The formal geometry of derived stacks}

As we have seen any derived scheme $X$, or more generally a derived Artin stack, 
has a truncation $t_{0}(X)$ and a natural morphism $j : t_0(X) \longrightarrow X$. 
We have already mentioned that $X$ behaves like a formal thickening of $t_0(X)$, and
in a way the difference between derived algebraic geometry and algebraic geometry
is concentrated at the formal level. We explore this furthermore in the present section, 
by explaining the deep interactions between derived algebraic geometry and 
formal/infinitesimal geometry. 

\subsection{Cotangent complexes and obstruction theory}

In $\S3.1$ we have seen that any derived scheme $X$ possesses a cotangent complex
$\mathbb{L}_{X}$. We will now explain how this notion extends to the more general
setting of derived Artin stacks, and how it controls obstruction theory. \smallbreak

Let $X$ be a derived Artin stack. We define its quasi-coherent derived $\sz$-category
$L_{\qcoh}(X)$ by integrating all quasi-coherent derived $\sz$-categories of derived
schemes over $X$. In a formula
$$L_{\qcoh}(X):=\displaystyle{\lim_{S\in \dSch/X}}L_{\qcoh}(S),$$
where the limit is taken along the $\sz$-category of all derived schemes
over $X$. By using descent, we could also restrict to affine derived scheme over $X$ and
get an equivalent definition. 

We define the cotangent complex of $X$ in a similar fashion as for derived schemes. 
For $M\in L_{\qcoh}(X)$, with cohomology sheaves concentrated in nonpositive degrees, 
we set $X[M]$, the trivial square zero infinitesimal  extension of $X$ by $M$. It is
given by the relative spectrum
$X[M] := \Spec\, (\OO_{X}\oplus M).$
The object $X[M]$ sits naturally under $X$, by means of the augmentation
$\OO_{X} \oplus M \longrightarrow \OO_{X}$. The cotangent complex of $X$ is the
object $\mathbb{L}_{X} \in L_{\qcoh}(X)$ 
such that for all $M \in L_{\qcoh}(X)$ as above, we have functorial equivalences
$$Map_{X/\dSt}(X[M],X) \simeq Map_{L_{\qcoh}(X)}(\mathbb{L}_{X},M).$$
The existence of such the object $\mathbb{L}_{X}$ is a theorem, whose proof can 
be found in \cite[Cor. 2.2.3.3]{hagII}. 

Cotangent complexes of derived Artin stacks behave similarly to the case
of derived schemes: functoriality and stability by base-change. In
particular, for a morphism between derived Artin stacks $f : X \longrightarrow Y$, 
we define the relative cotangent complex $\mathbb{L}_{f} \in L_{\qcoh}(X)$
has being the cofiber of the morphism $f^{*}(\mathbb{L}_{Y}) \rightarrow \mathbb{L}_{X}$. 
The
smooth and \'etale morphisms between derived Artin stacks have similar
characterizations using cotangent complexes (see \cite[\S2.2.5]{hagII}). A finitely presented
morphism $f : X \rightarrow Y$ between derived Artin stacks is \'etale if and only if
the relative cotangent complex $\mathbb{L}_{f}$ vanishes. The same morphism
is smooth if and only if the relative cotangent complex $\mathbb{L}_{f}$ 
has positive Tor amplitude. 

We note here that cotangent complexes of derived Artin stacks might not be
themselves cohomologically concentrated in  nonpositive degrees. It is a general
fact that if $X$ is a derived  $n$-Artin stack, in the sense of the inductive definition \ref{d4'}, 
then $\mathbb{L}_{X}$ is cohomologically concentrated in degree $]-\infty,n]$. This
can be seen inductively by using groupoid presentations as follows. Suppose that 
$X$ is the quotient of a smooth groupoid in derived $(n-1)$-Artin stacks $X_{*}$. 
We consider the unit section $e : X_{0} \longrightarrow X_{1}$, as
well as the natural morphism $X_{1} \longrightarrow X_{0}\times X_{0}$. We get a morphism of quasi-coherent 
complexes
on $X_{0}$
$$\mathbb{L}_{X_{0}}\simeq e^{*}(\mathbb{L}_{X_{0}\times X_{0}/X_{0}}) \longrightarrow
e^{*}(\mathbb{L}_{X_{1}/X_{0}}),$$
which is the infinitesimal action of $X_{1}$ on $X_{0}$. The fiber of this map 
is $\pi^{*}(\mathbb{L}_{X})$, where $\pi : X_ {0} \longrightarrow |X_{*}|\simeq X$ is 
the natural projection. This provides an efficient manner to understand
cotangent complexes of derived Artin stack by induction using presentations by quotient 
by smooth groupoids. \smallbreak

Let $Y=Spec\, A$ be an affine derived scheme and $M \in L_{\qcoh}(Y)\simeq L(A)$ a quasi-coherent
complexes cohomologically concentrated in strictly negative degrees. 
A morphism $d : \mathbb{L}_{Y} \longrightarrow M$ in $L_{\qcoh}(Y)$ corresponds to 
a morphism $(id,d) : A \longrightarrow A \oplus M$ of derived rings augmented to $A$.
We let $A\oplus_{d}M[-1]$ be the derived ring defined by the following
cartesian square in $\scomm$
$$\xymatrix{
A\oplus_{d}M[-1] \ar[r] \ar[d] & A \ar[d]^-{(id,d)} \\
A \ar[r]_-{(id,0)} & A\oplus M,}$$
and $Y_{d}[M[-1]]=\Spec\, (A\oplus_{d}M[-1])$ be the corresponding affine derived scheme. 
It comes equipped with a natural morphism $Y \longrightarrow Y_{d}[M[-1]]$, which by definition 
is the square zero extension of $Y$ by $M[-1]$ twisted by $d$. 

Assume now that $X$ is a derived Artin stack, and consider the following lifting problem. 
We assume given a morphism $f : Y \longrightarrow X$, and we consider the space
of all possible lifts of $f$ to $Y_{d}[M[-1]]$
$$L(f,M,d):=Map_{Y/\dSt}(Y_{d}[M[-1]],X).$$

The next proposition subsumes the content of 
the derived algebraic geometry approach to obstruction theory (see \cite[\S1.4.2]{hagII}).

\begin{prop}\label{p1}
With the above notations there is a canonical 
element $o(f,M,d) \in Ext^{0}(f^{*}(\mathbb{L}_{X}),M)$
such that $o(f,M,d)=0$ if and only if $L(f,M,d)$ is non-empty. Moreover, if
$o(f,M,d)=0$ then the simplicial set $L(f,M,d)$ is a
torsor over the simplicial abelian group $Map_{L_{\qcoh}(Y)}(f^{*}(\mathbb{L}_{X}),M[-1])$.
\end{prop}

A key feature of derived algebraic geometry is that the element $o(f,M,d)$ is functorial
in $X$ and in $M$. It can also be generalized to the case where 
$Y$ is no more affine and is itself a derived Artin stack. \smallbreak

The proposition \ref{p1} is an extremely efficient tool  to understand the
decomposition of the mapping spaces between derived schemes and derived Artin stacks
obtained by Postnikov decomposition. For this, let $X$ and $Y$ be
derived Artin stacks, and let $t_{\leq n}(X)$ and 
$t_{\leq n}(Y)$ be their Postnikov truncations. It can be shown that for each  $n$
the natural morphism $t_{\leq n}(X) \longrightarrow t_{\leq n+1}(X)$
is of the form
$$t_{\leq n}(X) \longrightarrow t_{\leq n}(X)_{d}[\pi_{n+1}(X)[n+1]],$$
for some map $d : \mathbb{L}_{t_{\leq n}(X)} \longrightarrow \pi_{n+1}(X)[n+2]$
(so here $M=\pi_{n+1}(X)[n+2]$). From this we deduce the shape of the fibers of the
morphism of spaces
$$Map_{\dSt}(t_{\leq n+1}(X),Y) \longrightarrow Map_{\dSt}(t_{\leq n}(X),Y).$$
For each $f : t_{\leq n}(X) \longrightarrow Y$, there is an obstruction element
$$
o(f,n) \in Ext^{n+2}(f^{*}(\mathbb{L}_{Y}),\pi_{n+1}(X)), 
$$ vanishing precisely
when $f$ lifts to a morphism from the next stage of the Postnikov tower $t_{\leq n+1}(X)$.
Moreover, when such a lift exists, the space of all lifts is, noncanonically, equivalent 
to $Map(f^{*}(\mathbb{L}_{Y}),\pi_{n+1}(X)[n+1])$. In particular, if non-empty, the
equivalence classes of lifts are in one-to-one correspondence with 
$Ext^{n+1}(f^{*}(\mathbb{L}_{Y}),\pi_{n+1}(X))$. 

One immediate consequence is the following co-connectivity statement: 
{any  $n$-truncated derived Artin stack $X$ such that $X=t_{\leq n}(X)$ has the property
that for any  derived $m$-Artin stack $Y$ (see definition \ref{d4'}), $Map_{\dSt}(X,Y)$ is an 
$(n+m)$-truncated simplicial set. 

\subsection{The idea of formal descent}

Because the $\sz$-category $\dSt$ of derived stacks is an $\sz$-topos, all the epimorphisms
between derived stacks are effective (see \cite[Thm. 4.9.2 (3)]{hagI}). One instance of this fact is that for
a smooth and surjective morphism of derived schemes (more generally 
of derived Artin stacks) $f : X \longrightarrow Y$, the object $Y \in \dSt$ can be recovered
by the following formula
$$Y \simeq |N(f)|,$$
where $N(f)$ is the nerve of the morphism $f$, that is the simplicial object
$[n] \mapsto X\times_{Y}X\times_{Y}\dots \times_{Y}X$, and $|N(f)|$ is the 
colimit of the simplicial object $N(f)$. 

Derived algebraic geometry proposes another form of the descent property in a rather
unusual context, namely when $f$ is now a closed immersion of schemes. This descent for 
closed immersions goes back to some fundamental results of Carlsson 
concerning completions in stable homotopy theory (see \cite{carl})
and is more subtle than the smooth descent just mentioned. It is however
an extremely nice and characteristic property of derived algebraic geometry which 
do not have any underived counterpart. 

We let $f : X \longrightarrow Y$ be a closed immersion of locally noetherian (underived) schemes.
We let $\hat{Y}_{X}$ be the formal completion of $Y$ along $X$. From the functorial point 
of view it is defined as the $\sz$-functor on the $\sz$-category of derived rings as follows.
For a derived ring $R$ we let $R_{red} :=\pi_{0}(R)_{red}$, the reduced
ring obtained from $\pi_{0}(R)$. The $\sz$-functor $\hat{Y}_{X}$ is then defined by the formula
$$\hat{Y}_{X}(R):=Y(R) \times_{Y(R_{red})}X(R_{red}).$$
As such, $\hat{Y}_{X}$ is a subobject of $Y$ because the morphism
$X(R_{red}) \longrightarrow Y(R_{red})$ is an injective map of sets. As a stack, 
$\hat{Y}_{X}$ is representable by a formal scheme, namely the formal completion of 
$Y$ along $X$. 

To the map $f$, we can form its nerve $N(f)$. This is the simplicial object 
in derived schemes obtained by taking the multiple fiber products of $X$ over $Y$.
In degree  $n$, $N(f)$ is the $n+1$-fold fibered product
$X\times_{Y}X\times_{Y} \dots \times_{Y}X$. The simplicial diagram of derived schemes $N(f)$ 
comes equipped with an augmentation to $Y$ which naturally factors through the subobject 
$\hat{Y}_{X}$. The following theorem is a direct re-interpretation of \cite[Thm. 4.4]{carl}.

\begin{theorem}\label{t3}
The augmentation morphism $N(f) \longrightarrow \hat{Y}_{X}$ exhibits $\hat{Y}_{X}$ 
as the colimit of the diagram $N(f)$ inside the $\sz$-category of derived schemes: 
for any derived scheme $Z$ we have an equivalence
$$Map_{\dSt}(\hat{Y}_{X},Z) \simeq \displaystyle{\lim_{[n] \in \Delta}} 
Map_{\dSch}(N(f)_{n},Z).$$
\end{theorem}

The above statement possesses a certain number of subtleties. First of all the
noetherian hypothesis is necessary, already in the affine case. Another subtle point,
which differs from the smooth descent we have mentioned, is that the colimit of 
$N(f)$ must be taken inside the $\sz$-category of derived schemes. The statement 
is wrong if the same colimit is considered in the $\sz$-category of derived stacks
for a simple reason: it is not true that any morphism $S \longrightarrow \hat{Y}_{X}$
factors locally for the \'etale topology through $f : X \longrightarrow Y$. Finally,
the possible generalizations of this statement to the setting of derived schemes
and derived Artin stacks require some care related to the size of the derived structures
sheaves. We refer to \cite[\S2.3]{gait} for more about formal completions in the general context of derived
Artin stacks. \smallbreak

On simple but instructing example of theorem \ref{t3} in action is the case
where $f$ is a closed point inside a smooth variety $Y$ over a field $k$, 
$x : Spec\, k \longrightarrow Y.$
The nerve of $x$ can be computed using Koszul resolutions obtained from
the choice of a system of local parameters at $x$ on $Y$. This gives 
$$N(x)_n \simeq \Spec\, A^{\otimes n},$$
where $A=Sym_{k}(V[1])$, where $V=\Omega^{1}_{V,x}$ is the cotangent space of $Y$ at $x$. 
Functions on the colimit of the simplicial derived scheme $N(x)$
is the the limit of the cosimplicial object $n \mapsto Sym_{k}(V^{n}[1])$, which can be
identified with $\widehat{Sym_{k}(V)}$, the completed symmetric algebra of $V$, 
or equivalently with the formal local ring $\OO_{V,x}$. It is interesting to note here that 
the limit of the co-simplicial diagram $n \mapsto Sym_{k}(V^{n}[1])$ lies in the
wrong quadrant and thus involves a nonconverging spectral sequence a priori. This nonconvergence 
is responsible for the completion of the symmetric algebra $\widehat{Sym_{k}(V)}$ as a final 
result. 

The case of a closed point $x : Spec\, k \longrightarrow Y$, for $Y$ a scheme of finite type
over $k$, possesses also an interpretation in terms of classifying spaces
of derived group schemes. This point of view, more topological, makes a clear link between 
derived algebraic geometry and algebraic topology. The basic observation here 
is that the nerve $N(x)$ is the nerve of a derived group scheme $\Omega_{x}Y=k\times_{Y}k$, 
called the based derived loop group of $Y$ at $x$. A reformulation of theorem \ref{t3}
is the existence of an equivalence of formal schemes
$$B(\Omega_{x}Y) \simeq Spf\, \widehat{\OO_{Y,x}}=\hat{Y}_{x},$$
or equivalently that $\hat{Y}_{x}$ is a classifying object inside formal schemes 
for the derived group scheme $\Omega_{x}Y$. This last equivalence should be understood
as an geometrico-algebraic version of the well known fact in homotopy theory, recovering 
a connected component of a topological space $Y$ containing a point $x \in Y$
 as the classifying space of $\Omega_{x}Y$, the based loop group of $Y$. 

As a final comment, the theorem \ref{t3} has also a another interpretation in terms of 
derived de~Rham theory of $X$ relative to $Y$. It is a special case of
a more general theorem relating derived de~Rham cohomology and algebraic de~Rham 
cohomology in characteristic zero (see \cite{bhat1}), which itself 
possesses $p$-adic and finite characteristic versions (see \cite{bhat2}).

\subsection{Tangent dg-lie algebras}

We assume now that $k$ is a base commutative ring of characteristic zero, 
and we work in $\dSt_{k}$, the $\sz$-category of derived stacks over $k$. For
a derived Artin stack $X$ locally of finite presentation over $k$, we have seen 
the existence of a cotangent complex
$\mathbb{L}_{X/k} \in L_{\qcoh}(X)$. Because $X$ is locally of finite presentation over $k$
the complex $\mathbb{L}_{X/k}$ is perfect and can be safely dualized to another
perfect complex $\mathbb{T}_{X/k}:=\mathbb{L}_{X/k}^{\vee}$, called the 
tangent complex. The cotangent complex $\mathbb{L}_{X/k}$ controls 
obstruction theory for $X$, but we will see now that 
$\mathbb{T}_{X/k}$ comes equipped with an extra structure of a (shifted) Lie algebra over $X$, 
which controls, in some sense,  the family of all formal completions of $X$ taken at various points. 
The existence of the Lie structure on $\mathbb{T}_{X/k}[-1]$ has been 
of a folklore idea for a while, with various attempts of construction. 
As an object in the non-$\sz$ derived category $D_{\qcoh}(X)$, and for
$X$ a smooth variety, this Lie structure is constructed in \cite[Thm. 2.6]{kapr}
(see also \cite{calacald} for a generalization).
More general approaches using derived loop spaces (see our next paragraph $\S4.4$)
appear in \cite{benznadl}: $\mathbb{T}_{X/k}[-1]$ is identified
with the Lie algebra of the derived loop stack $\mathcal{L}(X) \longrightarrow X$,
but assume that the relations between Lie algebras and formal groups
extend to the general setting of derived Artin stack (which is today not yet fully established). 
The very first general
complete construction appeared recently in \cite{henn}, following 
the general strategy of \cite{luri}. The main result of \cite{henn} can be
subsumed in the following theorem.

\begin{theorem}\label{t4}
With the above conditions and notations, there is a well defined structure of
an $\OO_{X}$-linear Lie algebra structure on $\mathbb{T}_{X/k}[-1]$. Moreover,
any quasi-coherent complex $M\in L_{\qcoh}(X)$ comes equipped with a canonical
action of $\mathbb{T}_{X/k}[-1]$. 
\end{theorem}

It is already noted in \cite{kapr} that the Lie algebra
$\mathbb{T}_{X/k}[-1]$ is closely related to the geometry of the diagonal
map $X \longrightarrow X \times X$, but a precise statement would require
a further investigation of the formal completions in the setting of derived
Artin stacks (see \cite[\ S 2.3]{gait}). However, 
it is possible to relate $\mathbb{T}_{X/k}[-1]$ with 
the various formal moduli problems represented by $X$ at each of its 
points. For this, let $x : Spec\, k \longrightarrow X$ be 
a global point, and let $l_{x}:=x^{*}(\mathbb{T}_{X/k}[-1])$, which is
a dg-lie algebra over $k$. By the main theorem of \cite{luri}, $l_{x}$
determines a unique $\sz$-functor
$F_{x} : \mathbf{dgart}^{*}_{k} \longrightarrow \Top,$
where $\mathbf{dgart}^{*}_{k}$ is the $\sz$-category of locak augmented
commutative dg-algebras over $k$ with finite dimensional total
homotopy (also called \emph{artinian commutative dg-algebras} over $k$). This
$\sz$-functor possesses several possible descriptions, one of them being very well known 
involving spaces of Mauer-Cartan elements. For 
each $A \in \mathbf{dgart}^{*}_{k}$ let $m_{A}$ be the kernel of the augmentation
$A \longrightarrow k$, and let us consider the space 
$\underline{MC}(l_{x}\otimes_{k}m_{A})$, of Mauer-Cartan elements in 
the dg-lie algebra $l_{x}\otimes_{k}m_{A}$ (see \cite{hini} for details). 
One possible definition for $F_{x}$ is $F_{x}(A)=\underline{MC}(l_{x}\otimes_{k}m_{A})$.

It can be checked that $F_{x}$ defined as above is equivalent to the 
formal completion $\hat{X}_{x}$ of $X$ at $x$ defined as follows.
The $\sz$-functor $\hat{X}_{x}$ simply is the restriction of the derived stack $X$ as an 
$\sz$-functor over 
$\mathbf{dgart}^{*}_{k}$, using $x$ as a base point: for all $A \in  \mathbf{dgart}^{*}_{k}$ we have
$$\hat{X}_{x}(A):=X(A)\times_{X(k)}\{x\}.$$
The equivalence $F_{x}\simeq \hat{X}_{x}$ can be interpreted as the statement that 
the dg-lie algebra $l_{x}$ does control the formal completion of $X$ at $x$. 
As $l_x$ is the fiber of the sheaf of Lie algebras $\mathbb{T}_{X/k}[-1]$, 
it is reasonable to consider that $\mathbb{T}_{X/k}[-1]$ encodes the family of 
formal completions $\hat{X}_{x}$, which is a family of formal moduli problems parametrized
by $X$. The total space of this family, which is still undefined in general, should of course be the 
formal completion of $X\times X$ along the the diagonal map. 

In the same way, for an object $M\in L_{\qcoh}(X)$, the $\mathbb{T}_{X/k}[-1]$-dg-module
structure on $M$ of theorem \ref{t4} can be restricted at a given 
global point $x : Spec\, k \longrightarrow X$. It provides an $l_{x}$-dg-module
structure on the fiber $x^{*}(M)$, which morally encodes the restriction of 
$M$ over $\hat{X}_{x}$, the formal completion of $X$ at $x$. 

At a global level, the Lie structure on $\mathbb{T}_{X/k}[-1]$ includes a bracket morphism in 
$L_{\qcoh}(X)$
$$[-,-] : \mathbb{T}_{X/k}[-1] \otimes_{\OO_{X}}
\mathbb{T}_{X/k}[-1] \longrightarrow \mathbb{T}_{X/k}[-1],$$
and thus a cohomology class $\alpha_{X} \in H^{1}(X,\mathbb{L}_{X/k} \otimes_{\OO_{X}} 
\underline{End}(\mathbb{T}_{X/k}))$. In the same way, a quasi-coherent module $M$, 
together with its $\mathbb{T}_{X/k}[-1]$-action provides a class 
$\alpha_{X}(M) \in H^{1}(X,\mathbb{L}_{X/k} \otimes_{\OO_{X}} 
\underline{End}(M))$. It is strongly believed that the class $\alpha_{X}(M)$ is
the Atiyah class of $M$, though the precise comparison is under investigation and not established yet. 

\subsection{Derived loop spaces and algebraic de~Rham theory}

We continue to work over a base commutative ring $k$ of characteristic zero. \smallbreak

Let $X$ be a derived Artin stack locally of finite presentation over $k$. 
We let $S^{1}:=B\mathbb{Z}$, the simplicial circle
considered as a constant derived stack $S^{1} \in \dSt_{k}$. 

\begin{df}\label{d}
The \emph{derived loop stack
of $X$ (over $k$)} is defined by 
$$\mathcal{L}X:=\Map(S^{1},X) \in \dSt_{k},$$
the derived stack of morphisms from $S^{1}$ to $X$.
\end{df}

The derived loop stack $\mathcal{L}X$ is an algebraic counter-part of the
free loop space appearing in string topology. Intuitively it consists of 
\emph{infinitesimal loops in $X$} and encodes many of the de~Rham theory 
of $X$, as we will going to explain now. \smallbreak

The constant derived stack $S^{1}$ can be written a push-out 
$S^{1}\simeq *\displaystyle{\coprod_{*\coprod *}}*,$
which implies the following simple formula for the derived loop stack
$$\mathcal{L}X\simeq \displaystyle{X\times_{X\times X}X},$$
from which the following descriptions of derived loop stacks follows.

\begin{itemize}

\item If $X=\Spec\, A$ is an affine derived scheme (over $k$), then 
so is $\mathcal{L}X$ and we have 
$$\mathcal{L}X \simeq \Spec\, (A\otimes_{A\otimes_{k}A}A).$$

\item For any derived scheme $X$ over $k$ the natural base point $* \in S^{1}=B\mathbb{Z}$
provides an affine morphism of derived schemes $\pi : \mathcal{L}X \longrightarrow X$. The
affine projection $\pi$ identifies $\mathcal{L}X$ with the relative spectrum 
of the symmetric algebra $Sym_{\OO_{X}}(\mathbb{L}_{X/k}[1])$ (see \cite{toenvezz5})
$$\mathcal{L}X \simeq \Spec\, (Sym_{\OO_{X}}(\mathbb{L}_{X/k}[1])).$$

\item Let $X$ be a (nonderived) Artin stack (e.g.\ in the sense of \cite{arti}), considered
as an object $X \in \dSt_{k}$. Then the truncation $t_{0}(\mathcal{L}X)$
is the so-called \emph{inertia stack of $X$} (also called \emph{twisted sectors}) which classifies
objects endowed with an automorphism.  The derived stack 
$\mathcal{L}X$ endows this inertia stack with a canonical derived structure.

\item For $X=BG$, for $G$ a smooth group scheme over $k$, we have 
$\mathcal{L}BG\simeq [G/G]$, where $G$ acts by conjugation on itself. 
More generally, for a smooth group scheme $G$ acting on a scheme $X$, 
we have 
$$\mathcal{L}[X/G] \simeq [X^{ts}/G],$$
where $X^{ts}$ is the derived scheme of fixed points defined as derived fiber product
$X^{ts}:=(X \times G)\times_{X\times X}X$. 

\end{itemize}

The first two properties above show that the geometry of the derived loop scheme 
$\mathcal{L}X$ is closely 
related to differential forms on the derived scheme $X$. When $X$ is no more a derived scheme but a derived Artin 
stack, similar  relations still hold but these are  more subtle. However, using descent 
for forms (see proposition \ref{p3}) 
it is possible to see the existence of a natural morphism
$$\mathbb{H}(\mathcal{L}X,\OO_{\mathcal{L}X}) \longrightarrow 
\mathbb{H}(X,Sym_{\OO_{X}}(\mathbb{L}_{X/k}[1]))\simeq
\mathbb{H}(X,\oplus_{p} (\wedge^{p}\mathbb{L}_{X/k})[p]).$$

The interrelations between derived loop stacks and differential forms become even more interesting when 
$\mathcal{L}X$ is considered equipped with the natural action of the group $S^{1}=B\mathbb{Z}$
coming from the $S^{1}$-action on itself by translation. In order to explain this we first need to 
remind some equivalences of $\sz$-categories in the context of \emph{mixed} and
\emph{$S^{1}$-equivariants} complexes.

It has been known for a while that the homotopy theory of simplicial $k$-modules endowed with 
an $S^{1}$-action is equivalent to the homotopy theory of nonpositively graded mixed complexes. 
This equivalence, suitably generalized, provide an equivalence between the $\sz$-category of commutative 
simplicial rings with an $S^{1}$-action, and the $\sz$-category of nonpositively graded mixed
commutative dg-algebras (see \cite{toenvezz5}). To be more precise, we let 
$S^{1}\mhyphen \scomm_{k}:=\fun(BS^{1},\scomm_{k})$, the the $\sz$-category of derived rings over $k$ together 
with an 
action of $S^{1}$, where $BS^{1}$ is the $\sz$-category with a unique object and $S^{1}$ as its 
simplicial monoid of 
endomorphisms. On the other hand, we set $k[\epsilon]$ be the dg-algebra over $k$ generated by 
a unique element $\epsilon$ in degree $-1$ and satisfying $\epsilon^2=0$. The dg-algebra
$k[\epsilon]$ can be identified with the homology algebra of the circle $H_{*}(S^{1})$ (with
the algebra structure induced from the group structure on $S^{1}$). The category of dg-modules over
$k[\epsilon]$ is, by definition and observation, the category of mixed complexes. The category 
$k[\epsilon]\mhyphen Mod$ comes equipped with a symmetric monoidal structure induced from the natural 
cocommutative bi-dg-algebra structure on $k[\epsilon]$. The commutative monoids in $k[\epsilon]\mhyphen Mod$
are called \emph{mixed commutative dg-algebras}, and consists of a 
commutative dg-algebra $A$ over $k$ together with a $k$-derivation $\epsilon : A \longrightarrow A[-1]$
of cohomological degree $-1$. By localization along the quasi-isomorphisms, the 
mixed commutative dg-algebras form an $\sz$-category $\ecdga_{k}$, and we restrict to its
full sub-$\sz$-category $\ecdga_{k}^{\leq 0} \subset \ecdga_{k}$ consisting of nonpositively graded 
objects. 

The normalization functor, from simplicial modules to nonpositively graded complexes, is shown in 
\cite{toenvezz5} to naturally extend to an equivalence of $\sz$-categories
$$\phi : S^{1}\mhyphen \scomm_{k} \simeq \ecdga_{k}^{\leq 0}.$$
This equivalence is not a formal result, and is achieved through a long sequence of equivalences 
between 
auxiliaries $\sz$-categories. One possible interpretation of the existence of the equivalence $\phi$ 
is the statement that the Hopf dg-algebra $C_{*}(S^{1},k)$, of chains on $S^{1}$ is 
formal (i.e., quasi-isomorphic to its cohomology). The nonformal nature of the equivalence $\phi$
can also be seen in the following result, which is a direct consequence of its existence. 
Let $A \in \scomm_{k}$ be a derived ring over $k$, and denote by $A$ again the 
commutative dg-algebra obtained out of $A$ by normalization. Then $S^{1}\otimes_{k} A$
defines an object in $S^{1}\mhyphen \scomm_{k}$, where $S^{1}$ acts on itself by translation. 
The dg-algebra $A$ also possesses a de~Rham complex 
$DR(A/k)=Sym_{A}(\mathbb{L}_{A}[1])$, which is an object in $\ecdga_{k}^{\leq 0}$ for which 
the mixed structure is induced by the de~Rham differential. Then 
we have $\phi(S^{1}\otimes A)\simeq DR(A/k)$, and this follows directly from the existence
of $\phi$ and the universal properties of the two objects $S^{1}\otimes A$ and $DR(A/k)$. 
Globally, on a general derived scheme, this result reads as below.

\begin{prop}{(See \cite{toenvezz5})}\label{p2}
Let $X$ be a derived scheme over $k$ and $\mathcal{L}X=\Map_{\dSt_{k}}(S^{1},X)$  its derived loop 
scheme over $k$, 
endowed with its natural action of $S^{1}$. Then there is an equivalence of stacks of 
mixed commutative dg-algebras over $X$
$$\phi(\OO_{\mathcal{L}X}) \simeq DR(\OO_{X}/k),$$
where $DR(\OO_{X}/k):=Sym_{\OO_{X}}(\mathbb{L}_{X/k}[1])$, with the mixed structure 
induced by the de~Rham 
differential.
\end{prop} 

When $X$ is no longer a derived scheme but rather a derived Artin stack, the proposition above
fails, simply because $X \mapsto \mathbb{H}(\mathcal{L}X,\OO_{\mathcal{L}X})$ does not satisfy descent for
smooth coverings. However, the equivalence above can be localized on the smooth topology in order
to obtain an analogous result for derived Artin stacks, for which the left hand side is replaced
by a suitable stackification of the construction $X \mapsto \OO_{\mathcal{L}X}$. In some
cases, for instance for smooth Artin stacks with affine diagonal, this stackification can be 
interpreted using a modified derived loop stack, such as the formal completion of
the derived loop stack (see for instance \cite{benznadl}).
A general corollary of proposition \ref{p2} is the following, concerning 
$S^{1}$-invariants functions on derived loop stacks. 

\begin{cor}\label{c2}
Let $X$ be a derived Artin stack over $k$ and $\mathcal{L}X$  its derived loop stack over $k$.
\begin{enumerate}

\item There exists a morphism of complexes of $k$-modules
$$\phi : \mathbb{H}(\mathcal{L}X,\OO_{\mathcal{L}X})^{S^{1}} \longrightarrow 
H_{DR}^{\mathbb{Z}/2}(t_{0}(X)/k),$$
where $H_{DR}^{\mathbb{Z}/2}(t_{0}(X)/k)$ is the $2$-periodic de~Rham cohomology of the truncation 
$t_{0}(X)$ of $X$.

\item If $X$ is a quasi-compact and quasi-separated derived scheme, then 
the morphism $\phi$ induces an equivalence
$$\phi : \mathbb{H}(\mathcal{L}X,\OO_{\mathcal{L}X})^{S^{1}}[\beta^{-1}] \simeq 
H^{\mathbb{Z}/2}_{DR}(t_{0}(X)/k).$$

\end{enumerate}
\end{cor}

In the corollary above, the de~Rham cohomology $H_{DR}^{\mathbb{Z}/2}(t_{0}(X)/k)$ is 
simply defined by the cohomology of $t_{0}(X)$ with coefficients in the 2-periodized 
algebraic de~Rham complex (in the sense of \cite{hart}). It can also be computed using the (negative) periodic cyclic 
of the derived stack $X$. In the second point, we use that $k^{S^{1}}\simeq k[\beta]$, 
whith $\beta$ in degree $2$, acts on (homotopy) fixed points of $S^{1}$ on any 
$S^{1}$-equivariant complex. \smallbreak

The above corollary can be used for instance to provide a new interpretation of the Chern 
character with coefficients in de~Rham cohomology. Indeed, for 
a derived Artin stack $X$ over $k$, and $E$ a vector bundle on $X$, or more generally
a perfect complex on $X$, the pull-back $\pi^{*}(E)$ on the derived loop stack
$\mathcal{L}X$ possesses a natural automorphism $\alpha_{E}$, obtained as the monodromy
operator along the loops. The formal existence of $\alpha_{E}$ 
follows from the fact that the projection $\pi : \mathcal{L}X \longrightarrow X$, 
as morphism in the $\sz$-category of derived stacks, possesses 
a self-homotopy given by the evaluation map
$S^{1} \times \mathcal{L}X \longrightarrow X$. This self-homotopy induces an
automorphism of the pull-back $\sz$-functor $\pi^{*}$. 

The automorphism $\alpha_{E}$ is another incarnation of the Atiyah class of $E$, and its
trace $Tr(\alpha_{E})$, as a function on $\mathcal{L}E$, can be shown to be 
naturally fixed by the $S^{1}$-action. This $S^{1}$-invariance is 
an incarnation of the well known cyclic invariance of traces, but its conceptual 
explanation is a rather deep phenomenon closely related to \emph{fully extended 1-dimensional
topological field theories in the sense of Lurie} (see \cite{toenvezz6}
for details). The trace $Tr(\alpha_{E})$ therefore
provides an element in $\mathbb{H}(\mathcal{L}X,\OO_{\mathcal{L}X})^{S^{1}}$. It is shown in \cite[App.~B]{toenvezz6}, 
as least when $X$ is a smooth and quasi-projective $k$-scheme, that $\phi(Tr(\alpha_{E})) \in 
H^{\mathbb{Z}/2}_{DR}(X/k)$
is the Chern character of $E$ in algebraic de~Rham cohomology. \smallbreak

The above interpretation of de~Rham cohomology classes in terms of $S^{1}$-equivariant
functions on the derived loop stack possesses a categorification relating quasi-coherent and 
$S^{1}$-equivariant sheaves on derived loop stacks and $\mathcal{D}$-modules. This relation
is again a direct consequence of the proposition \ref{p2} and can be stated as follows. 

\begin{cor}\label{c3}
Let $X$ be a smooth Artin stack locally of finite presentation 
over $k$ and $\mathcal{L}X$ be its derived loop stack over $k$.
\begin{enumerate}

\item There exists a natural $\sz$-functor
$$\phi : L_{\qcoh}^{S^{1}}(\mathcal{L}X) \longrightarrow 
L_{\qcoh}^{\mathbb{Z}/2}(\mathcal{D}_{X/k})$$
where $L_{\qcoh}^{\mathbb{Z}/2}(\mathcal{D}_{X/k})$ is the $2$-periodic derived
$\sz$-category of complexes of $\mathcal{D}_{X/k}$-modules with quasi-coherent cohomologies.

\item If $X$ is moreover a quasi-compact and quasi-separated smooth scheme over $k$, then 
the morphism $\phi$ induces an equivalence
$$\phi : L_{coh}^{S^{1}}(\mathcal{L}X)[\beta^{-1}] \simeq 
L_{coh}^{\mathbb{Z}/2}(\mathcal{D}_{X/k}).$$

\end{enumerate}
\end{cor}

In the above corollary $L_{\qcoh}^{S^{1}}(\mathcal{L}X)$ denotes the $S^{1}$-equivariant quasi-coherent 
derived
$\sz$-category of $\mathcal{L}X$, which can be defined for instance as 
the $\sz$-category of $S^{1}$-fixed points of the natural $S^{1}$-action on $L_{\qcoh}(\mathcal{L}X)$. 
The symbol $L_{coh}^{S^{1}}(\mathcal{L}X)$ denotes the full sub-$\sz$-category of
$L_{\qcoh}^{S^{1}}(\mathcal{L}X)$ consisting of bounded coherent objects ($E$ with
$\pi_{*}(E)$ coherent on $X$), and $L_{coh}^{\mathbb{Z}/2}(\mathcal{D}_{X/k})$
consists of 2-periodic complexes of $\mathcal{D}_{X/k}$-modules with coherent
cohomology (as sheaves of $\mathcal{D}_{X/k}$-modules). Finally, $L_{coh}^{S^{1}}(\mathcal{L}X)$
is an $\sz$-category which comes naturally enriched over $L_{coh}^{S^{1}}(k)\simeq L_{\perf}(k[\beta])$, 
which allows us to localize along $\beta$. We refer to \cite{benznadl,prey} for more details about the objects
involved in the corollary \ref{c3}, as well as for possible generalizations and modifications (e.g.\ 
for the non-periodic version). 

To finish this paragraph, we mention that the interpretation of the Chern  character as the trace
of the universal automorphism on $\mathcal{L}X$ can be also categorified in an interesting manner. 
We now start with a stack of dg-cateories $\mathcal{T}$ over $X$ (see \cite{toen3,toenvezz6}), which is
a categorification of a quasi-coherent sheaf. The pull-back of $\mathcal{T}$ over
$\mathcal{L}X$ also possesses a universal automorphism $\alpha_{\mathcal{T}}$, 
which itself has a well-defined trace. This trace is no more a function but rather is a 
quasi-coherent sheaf on $\mathcal{L}X$, which also turns out to carry a natural 
$S^{1}$-equivariant structure. It is therefore an object in $L_{\qcoh}^{S^{1}}(\mathcal{L}X)$, 
and by the corollary \ref{c3} its image by $\phi$ becomes a $2$-periodic $\mathcal{D}_{X/k}$-module
over $X$. This $\mathcal{D}_{X/k}$-module must be interpreted as the family of periodic 
homology of the family $\mathcal{T}$, endowed with a noncommutative version of the Gauss-Manin 
connection (see \cite{toenvezz6}). This is the first step in the general construction of variations of 
noncommutative
Hodge structures in the sense of \cite{katzkontpant} and is a far-reaching generalization of the noncommutative 
Gauss-Manin
connection constructed on flat families of algebras in \cite{tsyg}. In the same way that the Chern character
of perfect complexes can be understood using 1-dimensional fully extended TQFT in the sense of Lurie, 
this noncommutative Gauss-Manin connection, which is a categorification of the Chern character, 
can be treated using 2-dimensional fully extended TQFT's. 

\section{Symplectic, Poisson and Lagrangian structures in the derived setting}

At the end of the section $\S4$ we have seen the relations between derived loop stacks and
algebraic de~Rham theory. We now present further materials about differential forms on derived
Artin stacks, and introduce the notion of shifted symplectic structure. We will 
finish the section by some words concerning the dual notion of shifted Poisson structures, and
its possible importance for deformation quantization in the derived setting. \smallbreak

All along this section $k$ will be a noetherian commutative base ring, assumed to be of characteristic zero.

\subsection{Forms and closed forms on derived stacks}

Let $X$ be a derived Artin stack locally of finite presentation over $k$. We have seen in $\S4.1$ that 
$X$ admits a (relative over $k$) cotangent complex  
$\mathbb{L}_{X/k}$, which is the derived version of the sheaf of $1$-forms. For $p\geq 0$,  the complex
of $p$-forms on $X$ (relative to $k$) can be naturally defined as follows
$$\mathcal{A}^{p}(X):=\mathbb{H}(X,\wedge^{p}_{\OO_{X}}\mathbb{L}_{X/k}).$$
By definition, for $n\in \mathbb{Z}$, a $p$-form of degree  $n$ on $X$ is an element in $H^{n}
(\mathcal{A}^{p}(X))$, or
equivalently an element in $H^{n}(X,\wedge^{p}_{\OO_{X}}\mathbb{L}_{X/k})$. When $X$ is smooth over $k
$, 
all the perfect complexes $\wedge^{p}_{\OO_{X}}\mathbb{L}_{X/k}$ have nonnegative Tor amplitude, and 
thus
there are no nonzero $p$-forms of negative degree on $X$. In a dual manner, 
if $X$ is an affine derived scheme, then $X$ does not admit any nonzero $p$-form of positive degree. 
This is not true anymore without these hypothesis, and in general
a given derived Artin stack might have nonzero $p$-forms of arbitrary degrees. 

An important property of forms on derived Artin stacks is the smooth descent property, which is
a powerful computational tool as we will see later on with some simple examples. It can be stated as 
the following
proposition. 

\begin{prop}\label{p3}
Let $X$ be a derived Artin stack over $k$ (locally of finite presentation by our assumption), and let 
$X_{*}$ be a smooth Segal groupoid in derived Artin stack whose quotient $|X_{*}|$ is equivalent to $X$ 
(see $\S3.3$). Then, 
for all $p\geq 0$, the natural morphism
$$\mathcal{A}^{p}(X) \longrightarrow \displaystyle{\lim_{[n] \in \Delta}}\mathcal{A}^{p}(X_{n})$$
is an equivalence (in the $\sz$-category of complexes).
\end{prop}

Here are two typical examples of complexes of forms on some fundamental derived stacks. 

\subsubsection*{Forms on classifying stacks, and on $\Parf$}  We let $G$ be a smooth group scheme over 
$Spec\, k$, and 
$X=BG$ be its classifying stack. The cotangent complex of $X$ is $\frak{g}^{\vee}[-1]$, where 
$\frak{g}$ is the 
Lie algebra of $G$, considered as a quasi-coherent sheaf on $BG$ by the adjoint representation. 
We thus have $\wedge^{p}_{\OO_{X}}\mathbb{L}_{X/k}\simeq Sym^{p}_{k}(\frak{g}^{\vee})[-p]$, and the 
complex
of $p$-forms on $X$ is then the cohomology complex $\mathbb{H}(G,Sym^{p}_{k}(\frak{g}^{\vee}))[-p]$, 
of the group scheme $G$ with values in the representation $Sym^{p}_{k}(\frak{g}^{\vee})$. When $G$ is 
a reductive group scheme over $k$, its cohomology vanishes and the complex of $p$-forms
reduces to $Sym^{p}_{k}(\frak{g}^{\vee})^{G}[-p]$, the $G$-invariant symmetric $p$-forms on $\frak{g}$
sitting in cohomological degree $p$. In other words, when $G$ is reductive, there are no nonzero
$p$-forms of degree $n\neq p$ on $X$, and $p$-forms of degree $p$ are given by 
$Sym_k^{p}(\frak{g}^{\vee})$.

We let $X=\Parf$ be the derived stack of perfect complexes (see theorem \ref{t1})  
The cotangent complex $\mathbb{L}_{X/k}$ has the following description. There is a 
universal
perfect complex $\mathcal{E} \in L_{\qcoh}(X)$, and we have $\mathbb{L}_{X/k}\simeq \underline{End}
(\mathcal{E})[-1]$,
where $\underline{End}(\mathcal{E})=\mathcal{E}\otimes_{\OO_{X}}\mathcal{E}^{\vee}$ is the stack of 
endomorphisms
of $\mathcal{E}$. We obtain the following description of the complex of $p$-forms on the derived stack 
$\Parf$
$$\mathcal{A}^{p}(\Parf)\simeq \mathbb{H}(\Parf,Sym_{\OO_{\Parf}}(\underline{End}(\mathcal{E})))[-p].$$

\subsubsection*{Forms on a derived quotient stack}  Let $G$ be a reductive smooth group scheme over $k$ acting
on an affine derived scheme $Y=\Spec\, A$, and let $X:=[Y/G]$ be the quotient derived stack. 
The cotangent complex $\mathbb{L}_{X/k}
$, pulled back to
$Y$, is given by the fiber $\mathbb{L}$ of the natural morphism
$\rho : \mathbb{L}_{Y/k} \longrightarrow \OO_{Y} \otimes_{k}\frak{g}^{\vee}$, dual to the infinitesimal 
action of $G$ on $Y$. 
The group $G$ acts on $Y$ and on the morphism above, and thus on $\mathbb{L}$. 
The complex of $p$-forms on $X$ is then given by
$\mathcal{A}^{p}(X) \simeq (\wedge^{p}\mathbb{L})^{G}.$
Concretely, the fiber of the morphism $\rho$ can be described as the complex 
$\mathbb{L}_{Y/k} \oplus (\OO_{Y} \otimes_{k}\frak{g}^{\vee}[-1])$ endowed with a suitable differential 
coming from
the $G$-action on $Y$. The complex of $p$-forms $\mathcal{A}^{p}(X)$ can then be described as
$$\mathcal{A}^{p}(X) \simeq \left(\bigoplus_{i+j=p}(\wedge^{i}_{A}\mathbb{L}_{A}) 
\otimes_{k}Sym^{j}_{k}(\frak{g}^{\vee})[-j]\right)^{G},$$
again with a suitable differential. \smallbreak

We now define closed $p$-forms on derived Artin stacks. For this we start to treat the affine
case and then define the complex of closed $p$-forms on a derived stack $X$ by
taking a limit over all affine derived schemes mapping to $X$. 

Let $A$ be a derived ring over $k$, and $N(A) \in \ncdga_{k}$ be its normalization, which 
is a nonpositively graded commutative dg-algebra over $k$. We let $A'$ be a cofibrant
model for $A$ and we consider $\Omega_{A'}^{1}$ the $A'$-dg-module 
of K\"alher differential over $A'$ over $k$ (see e.g.\ \cite{behr,ptvv}). Note that under the equivalence of 
$\sz$-categories $L(N(A)) \simeq L(A')$, the object $\Omega_{A'}^{1}$ 
can be identified with $\mathbb{L}_{A/k}$, the cotangent complex of $A$ over $k$. 
We set $\Omega^{i}_{A'}:=\wedge^{i}_{A'}\Omega^{1}_{A'}$ for all $i\geq 0$. As for the case
of non-dg algebras, there is a de~Rham differential $dR : \Omega^{i}_{A'} \longrightarrow
\Omega_{A'}^{i+1}$, which is a morphism of complexes of $k$-modules and satisfies
$dR^{2}=0$. The differential $dR$ is also characterized by the
property that it endows $Sym_{A'}(\Omega^{1}_{A'}[1])$ with a structure of a graded mixed
commutative dg-algebra, which coincides with the universal derivation
$A' \longrightarrow \Omega^{1}_{A'}$ in degree $1$. 

We define a complex of $k$-modules $\mathcal{A}^{p,cl}(A)$, of closed $p$-forms
over $A$ (relative to $k$) as follows. The underlying graded $k$-module is given by
$$\mathcal{A}^{p,cl}(A):=\displaystyle{\prod_{i\geq 0}}\Omega^{p+i}_{A'}[-i].$$
The differential $D$ on $\mathcal{A}^{p,cl}(A)$ is defined to be the total differential combining 
the cohomological differential $d$ on $\Omega^{i}_{A'}$, and the de~Rham differential $dR$. 
In formula we have, for an element of degree  $n$, 
$\{\omega_{i}\}_{i\geq 0} \in \displaystyle{\prod_{i\geq 0}}(\Omega^{p+i}_{A'})^{n-i}$
$$D(\{\omega_{i}\}) := \{dR(\omega_{i-1})+d(\omega_{i})\}_{i\geq 0} \in 
\displaystyle{\prod_{i\geq 0}}(\Omega^{p+i}_{A'})^{n-i+1}.$$
The complex of closed $p$-forms $\mathcal{A}^{p,cl}(A)$ is functorial in $A$ and provides
an $\sz$-functor $\mathcal{A}^{p,cl}$ from the $\sz$-category of derived rings over $k$ to the
$\sz$-category $\dg_{k}$ of complexes of $k$-modules. This $\sz$-functor satisfies \'etale descent
and can then be left Kan extended to all derived stacks (see $\S2.1.2$)
$\mathcal{A}^{p,cl} : \dSt_{k}^{op} \longrightarrow \dg_{k}.$
For a derived stack $X$ we have by definition
$\mathcal{A}^{p,cl}(X)\simeq \displaystyle{\lim_{\Spec\, A \rightarrow X}}\mathcal{A}^{p,cl}(A)$. \smallbreak

The relation between closed $p$-forms and $p$-forms is based on the descent property 
\ref{p3}. The projection to the first factor $\displaystyle{\prod_{i\geq 0}}\Omega^{p+i}_{A'}[-i]
\longrightarrow \Omega_{A'}^{p}$ provides a morphism of $\sz$-functors 
$\mathcal{A}^{p,cl} \longrightarrow \mathcal{A}^{p}$ defined on derived rings over $k$. For a derived
Artin stack $X$ over $k$, and because of proposition \ref{p3}, we obtain a natural morphism
$$\mathcal{A}^{p,cl}(X) \simeq \displaystyle{\lim_{\Spec\, A \rightarrow X}}\mathcal{A}^{p,cl}(A)
\longrightarrow \displaystyle{\lim_{\Spec\, A \rightarrow X}}\mathcal{A}^{p}(A) \simeq 
\mathcal{A}^{p}(X).$$

\begin{remark}\label{r2}
It is important to note that the morphism above $\mathcal{A}^{p,cl}(X) \longrightarrow
\mathcal{A}^{p}(X)$ can have rather complicated fibers, contrary to the intuition 
that closed forms form a subspace inside the space of all forms. In the derived setting, 
\emph{being closed} is no more a property but becomes an extra structure: a given 
form might be closed in many nonequivalent manners. This degree of freedom will be 
essential for the general theory, but will also create some technical complications, 
as the construction a closed form will require in general much more work than 
constructing its underlying nonclosed form. 
\end{remark}

\begin{remark}\label{r2'}
Another comment concerns the relation between complexes of closed forms and
negative cyclic homology of commutative dg-algebra. For a commutative dg-algebra
$A$, we have its complex of negative cyclic homology $HC^-(A)$, as well
as its part of degre $p$ for the Hodge decomposition $HC^{-}(A)^{(p)}$
(see e.g.\ \cite{loda}). The so-called HKR theorem implies that we have a natural equivalence of complexes
$$\mathcal{A}^{p,cl}(A) \simeq HC^{-}(A)^{(p)}[-p].$$
\end{remark}

\subsubsection*{Closed forms on smooth schemes} Let $X$ be a smooth scheme over $k$ of relative dimension $d$. 
The complex $\mathcal{A}^{p,cl}(X)$ is nothing else than the standard truncated 
de~Rham complex of $X$, and is given by
$$\mathcal{A}^{p,cl}(X) \simeq \mathbb{H}(X,\Omega_{X/k}^{p} \rightarrow \Omega_{X/k}^{p+1} 
\rightarrow \dots \rightarrow \Omega_{X/k}^{d}).$$
In particular, $H^{0}(\mathcal{A}^{p,cl}(X))$ is naturally isomorphic to the space 
of closed $p$-forms on $X$ in the usual sense. Note also that when 
$X$ is moreover proper, the morphism $\mathcal{A}^{p,cl}(X) \longrightarrow
\mathcal{A}^{p}$ is injective in cohomology because of the degeneration of the Hodge to de~Rham
spectral sequence. This is a very special behavior of smooth and proper schemes, for which 
\emph{being closed} is indeed a well defined property. 

\subsubsection*{Closed forms on classifying stacks} Let $G$ be reductive smooth group scheme
over $k$ with Lie algebra $\frak{g}$. 
The complex of closed $p$-forms on $BG$ can be seen, using for instance
the proposition \ref{p3}, to be naturally equivalent to 
$\oplus_{i\geq 0}Sym^{p+i}(\frak{g}^{\vee})^{G}[-p-2i]$
with the zero differential. In particular, any element in $Sym^{p}(\frak{g}^{\vee})^{G}$
defines a canonical closed $p$-form of degree $p$ on $BG$. In this example
the projection $\mathcal{A}^{p,cl}(BG) \longrightarrow \mathcal{A}^{p}(BG)$
induces an isomorphism on the $p$-th cohomology groups: any $p$-form of degree
$p$ on $BG$ is canonically closed. This is again a specific property of 
classifying stacks of reductive groups. 

\subsubsection*{The canonical closed $2$-form of degree $2$ on $\Parf$} Let $\Parf$ be the
derived stack of perfect complexes over $k$ (see theorem \ref{t1}). We have seen that its tangent complex
$\mathbb{T}_{\Parf}$ is given by $\underline{End}(\mathcal{E})[1]$, the shifted
endomorphism dg-algebra of the universal perfect complex $\mathcal{E}$ on $\Parf$. We can therefore
define a $2$-form of degree $2$ on $\Parf$ by considering
$$Tr : \underline{End}(\mathcal{E})[1] \otimes \underline{End}(\mathcal{E})[1] \longrightarrow 
\OO_{X}[2],$$
which is, up to a suitable shift, the morphism obtained by taking the trace of the
multiplication in $\underline{End}(\mathcal{E})$.

The above $2$-form of degree $2$ has a canonical lift to a closed $2$-form
of degree $2$ on $\Parf$ obtained as follows. We consider the Chern character 
$Ch(\mathcal{E})$ of the universal object, as an object 
of negative cyclic homology $HC^{-}_{0}(\mathbb{\Parf})$. The part of
weight $2$ for the Hodge decomposition on $HC^{-}_{0}(\mathbb{\Parf})$ provides
a closed $2$-form $Ch_{2}(\mathcal{E})$ of degree $2$ on $\Parf$. It can be checked that 
the underlying $2$-form of $Ch_{2}(\mathcal{E})$ is $\frac{1}{2}.Tr$ (see \cite{ptvv} for details).
\smallbreak

By letting  $p=0$ in the definition of closed $p$-forms we
obtain the derived de~Rham complex of derived Artin stack. More explicitely, 
for a derived ring $A$ over $k$, with normalization $N(A)$ and cofibrant model 
$A'$, we set 
$\mathcal{A}^{*}_{DR}(A):=\prod_{i\geq 0}\Omega^{i}_{A'}[-i],$
with the exact same differential $D=dR+d$, sum of the cohomological and
de~Rham differential. For a derived Artin stack $X$ over $k$ we set 
$$\mathcal{A}_{DR}^{*}(X):=\displaystyle{\lim_{\Spec\, A \rightarrow X}}\mathcal{A}^{*}_{DR}(A),$$
and call it the \emph{derived de~Rham complex of $X$ over $k$}. It is also the complex
of closed $0$-forms on $X$, and we will simply denote
by $H^{*}_{DR}(X)$ the cohomology of the complex $\mathcal{A}^{*}_{DR}(X)$. 
There are obvious inclusions of subcomplexes
$\mathcal{A}^{p,cl}(A)[-p] \subset \mathcal{A}^{p-1,cl}(A)[-p+1] \subset \dots \subset 
\mathcal{A}^{*}_{DR}(A)$, 
inducing a tower of morphisms of complexes
$$
\dots \to  
\mathcal{A}^{p,cl}(X)[-p] \to
\mathcal{A}^{p-1,cl}(X)[-p+1] \to \dots  \to \mathcal{A}^{1,cl}(X)[-1]
\to \mathcal{A}^{*}_{DR}(X),
$$
which is an incarnation of the Hodge filtration on de~Rham cohomology: the cofiber of each 
morphism $\mathcal{A}^{p,cl}(X)[-p] \longrightarrow \mathcal{A}^{p-1,cl}(X)[-p+1]$ is
the shifted complex of $(p-1)$-forms on $X$, $\mathcal{A}^{p-1}(X)[-p+1]$.

By combining \cite{feigtsyg} and \cite{good}, it can be shown 
that $\mathcal{A}^{*}_{DR}(X)$ does compute the algebraic de~Rham 
cohomology of the truncation $t_{0}(X)$. 

\begin{prop}\label{p4}
Let $X$ be a derived Artin stack locally of finite presentation over $k$. 
\begin{enumerate}
\item The natural morphism $j : t_{0}(X) \longrightarrow X$ induces an equivalence of 
complexes of $k$-modules
$$j^{*} : \mathcal{A}^{*}_{DR}(X) \simeq \mathcal{A}^{*}_{DR}(t_{0}(X)).$$

\item There exists a natural equivalence between $\mathcal{A}^{*}_{DR}(t_{0}(X))$
and the algebraic de~Rham cohomology complex of $X$ relative to $k$ in the sense
of \cite{hart} (suitably extended to Artin stacks by descent). 

\end{enumerate}
\end{prop}

The above proposition has the following important consequence. Let $\omega
 \in H^{n}(\mathcal{A}^{p,cl}(X))$ be closed $p$-form of degree  $n$ on $X$. It defines a class
in the derived de~Rham cohomology $[\omega] \in H^{n+p}_{DR}(X)$. 

\begin{cor}\label{c4}
With the above notations and under the condition that $k$ is a field, 
we have $[\omega]=0$ as soon as $n < 0$. 
\end{cor}

The above corollary follows from the proposition \ref{p4} together with 
the canonical resolution of singularities and the proper descent for algebraic de~Rham cohomology.
Indeed, proposition \ref{p4} implies that we can admit that $t_{0}(X)=X$. By the canonical resolution of 
singularities and proper descent we can check that the natural morphism 
$H^{*}_{DR}(X) \longrightarrow H^{*,naive}_{DR}(X)$ is injective, where 
$H^{*,naive}_{DR}(X)$ denotes the hyper-cohomology of $X$ with coefficients 
in the naive de~Rham complex 
$$H^{*,naive}_{DR}(X):=\mathbb{H}(X,\xymatrix{\OO_{X} \ar[r] & \Omega_{X}^{1} \ar[r]
& \dots \ar[r] & \Omega_{X}^{d} \ar[r] & \dots}).$$
But by definition the image of $[\omega]$ is clearly zero in $H^{n+p,naive}_{DR}(X)$
when $n<0$. 

A consequence of corollary \ref{c4} is that any 
closed $p$-form $\omega$ of degree $n<0$ is \emph{exact}: it lies in the image of the natural 
morphism
$dR : \mathcal{A}^{\leq p-1}_{DR}(X)[p-1] \longrightarrow \mathcal{A}^{p,cl}(X),$
where $\mathcal{A}^{\leq p-1}_{DR}(X)$ is the $(p-1)$-truncated derived de~Rham complex defined
a the cofiber of $\mathcal{A}^{p,cl}(X)[-p] \longrightarrow \mathcal{A}^{*}_{DR}(X)$. 
This has many important implications, because in general a closed $p$-form involves an infinite
number of data (because of the infinite product appearing in the definition of $\mathcal{A}^{p,cl}(X)$), 
but cocycles in $\mathcal{A}^{\leq p-1}_{DR}(X)$ only involve a finite number of data. Closed
$p$-forms of negative degrees are somehow easier to understand that their positive degree counterpart.
As an example, we quote the Darboux lemma in for shifted symplectic structures
of negative shifts (see \cite{bravjoyc,bogr}).

\subsection{Symplectic and Lagrangian structures}

We now arrive at the central notion of shifted symplectic and Lagrangian structures. We start
with the following key definition. 

\begin{df}\label{d5}
Let $X$ be a derived Artin stack locally of finite presentation over $k$ and $n\in \mathbb{Z}$.
\begin{itemize}
\item
An \emph{n-shifted symplectic structure on $X$} is the datum of a closed $2$-form 
$\omega$ of degree  $n$, such that the underlying $2$-form on $X$ is nondegenerate: 
the adjoint morphism
$$\Theta_{\omega} : \mathbb{T}_{X/k} \longrightarrow \mathbb{L}_{X/k}[n]$$
is an equivalence in $L_{\qcoh}(X)$. 
\item Let $\omega$ be an  $n$-shifted symplectic structure on $X$, and let 
$Y$ be another derived Artin stack together with a morphism $f : Y \longrightarrow X$. 
A \emph{Lagrangian structure on $f$} consists of a homotopy $h : f^{*}(\omega) \sim 0$ 
in the complex $\mathcal{A}^{p,cl}(Y)$, such that the induced morphism
$$\Theta_{\omega,h} : \mathbb{T}_{Y/k} \longrightarrow \mathbb{L}_{Y/X}[n-1]$$
is an equivalence in $L_{\qcoh}(Y)$.  
\end{itemize}
\end{df}

Some comments about the above definition. By definition $\omega$ is an element 
in $H^{n}(\mathcal{A}^{p,cl}(X))$. The morphism $\Theta_{\omega}$ is defined
by considering the image of $\omega$ in $H^{n}(\mathcal{A}^{p}(X))\simeq
H^{n}(X,\wedge^{2}_{\OO_{X}}\mathbb{L}_{X/k})=[\OO_{X},(\wedge^{2}_{\OO_{X}}\mathbb{L}_{X/k})[n]]$, 
which by duality provide a morphism 
$\Theta_{\omega} : \mathbb{T}_{X/k} \longrightarrow \mathbb{L}_{X/k}[n]$.

In the same way, the  morphism $\Theta_{h,\omega}$ is defined as follows. The
homotopy $h$ provides a homotopy in $\mathcal{A}^{p}(Y)$, which is
a homotopy to zero of the following composition in $L_{\qcoh}(Y)$
$$\xymatrix{
\mathbb{T}_{Y/k} \ar[r] & f^{*}(\mathbb{T}_{X/k}) \ar[r]^-{f^{*}(\Theta_{\omega})} & 
f^{*}(\mathbb{L}_{X/k})[n] \ar[r] & \mathbb{L}_{Y/k}[n].}$$
This homotopy to zero defines a unique morphism in $L_{\qcoh}(Y)$ from 
$\mathbb{T}_{Y/k}$ to the fiber of $f^{*}(\mathbb{L}_{X/k})[n] \longrightarrow \mathbb{L}_{Y/k}[n]$, 
which is $\mathbb{L}_{Y/X}[n-1]$. 

\begin{remark}\label{r3}
A trivial, but conceptually important remark, communicated to me by D. Calaque, 
is that the notion of a Lagrangian 
structure is a generalization of the notion of shifted symplectic structure. To see
this we let $*=Spec\, k$ be endowed with the zero $(n+1)$-shifted symplectic structure. 
Then an  $n$-shifted symplectic structure on $X$ simply is a Lagrangian structure on 
the natural morphism $X \longrightarrow Spec\, k$. 
\end{remark}

Before stating the main existence results for shifted symplectic and Lagrangian structures, 
we present some more elementary properties as well as some relations with 
standard notions of symplectic geometry such Hamiltonian action and symplectic reduction. 

\subsubsection*{Shifted symplectic structures and amplitude} Let $X$ be
a derived Artin stack locally of finite presentation over $k$. Unless
in the situation where $X$ is \'etale over $Spec\, k$, that is
$\mathbb{L}_{X/k}\simeq 0$, there can be at most one integer  $n$ such that
$X$ admits an  $n$-shifted symplectic structure. Indeed, because $X$ is locally of finite
presentation over $k$ the tangent complex $\mathbb{T}_{X/k}$ is perfect of some
bounded amplitude, and therefore cannot be equivalent in $L_{\qcoh}(X)$ to
a nontrivial shift of itself. 

From a general point of view, if $X$ is a derived scheme, or more generally 
a derived Deligne-Mumford stack, then $X$ can only admit nonpositively shifted
symplectic structure because $\mathbb{T}_{X/k}$ is in this case of nonnegative
amplitude. Dually, if $X$ is a smooth Artin stack over $k$ its tangent 
complex $\mathbb{T}_{X/k}$ has nonpositive amplitude and thus $X$ can only carry 
nonnegatively shifted symplectic structures. In particular, a smooth 
Deligne-Mumford can only admit $0$-shifted symplectic structures, and these are
nothing else than the usual symplectic structures. 

\subsubsection*{Shifted symplectic structures on $BG$ and on $\Parf$}  Let $G$ be a reductive smooth
group scheme over $k$. In our last paragraph $\S5.1$ 
we have already computed the complex of closed $2$-forms on $BG$. We deduce
from this that $2$-shifted symplectic structures on $BG$ are one to one correspondence
with nondegenerate $G$-invariant scalar products on $\frak{g}$.  

As for the derived stack of perfect complexes $\Parf$, we have seen the existence
of canonical closed $2$-form of degree $2$ whose underlying $2$-form is the trace map
$$\frac{1}{2}.Tr : \underline{End}(\mathcal{E})[1] \otimes  \underline{End}(\mathcal{E})[1]
\longrightarrow \OO_{X}[2].$$
The trace morphism is clearly nondegenerate, as this can be checked at closed points, and  for which this
reduces to the well-known fact that $(A,B) \mapsto Tr(A.B)$ is a nondegenerate
pairing on the space of (graded) matrices. 

\subsubsection*{Lagrangian intersections} Let $X$ be a derived Artin stack locally
of finite presentation over $k$ and $\omega$ an $n$-shifted symplectic structure on $X$. 
Let $Y \longrightarrow X$ and $Z \longrightarrow X$ be two morphisms of derived Artin stacks
with Lagrangian structures. Then the fibered product derived stack $Y\times_{X}Z$
carries a canonical $(n-1)$-shifted symplectic structure. The closed $2$-form 
of degree $n-1$ on $Y\times_{X}Z$ is simply obtained by pulling-back 
$\omega$ to $Y\times_{X}Z$, which by construction comes with two homotopies to zero
coming from the two Lagrangian structures. These two homotopies combine to 
a self-homotopy of $0$ in $\mathcal{A}^{2,cl}(Y\times_{X}Z)[n]$, which 
is nothing else than a well defined element in $H^{n-1}(\mathcal{A}^{2,cl}(Y\times_{X}Z))$.
This closed $2$-form of degree $(n-1)$ on $Y\times_{X}Z$ can then be checked
to be nondegenerate by a direct diagram chase using the nondegeneracy property 
of Lagrangian structures. We refer to \cite{ptvv} for more details. 

The above applies in particular when $X$ is a smooth scheme and is symplectic in the usual 
sense, and $Y$ and $Z$ are two smooth subschemes of $X$  which are Lagrangian
in $X$ in the standard sense. Then the derived scheme $Y\times_{X}Z$ carries
a canonical $(-1)$-shifted symplectic structure. This has strong consequences
for the singularities and the local structure of $Y\times_{X}Z$. For instance, 
it is shown in \cite{bravjoyc,bogr} that locally for the Zariski topology $Y\times_{X}Z$
is the derived critical locus (see below) of a function on  a smooth scheme.

\subsubsection*{Shifted cotangent stacks and derived critical loci} Let $X$ be
a derived Artin stack locally of finite presentation over $k$. For $n\in \mathbb{Z}$ we consider
the shifted tangent complex $\mathbb{T}_{X/k}[-n]$ as well the corresponding linear derived stack 
(see $\S3.3$)
$$T^{*}X[n]=\Spec\, (Sym_{\OO_{X}}(\mathbb{T}_{X/k}[-n])) = \mathbb{V}(\mathbb{T}_{X/k}[-n]) \longrightarrow 
X.$$
The derived stack $T^{*}X[n]$ is called the \emph{n-shifted cotangent stack of $X$}. In the same
way that total cotangent space of smooth schemes carries a canonical symplectic structure, the shifted
cotangent stack $T^{*}X[n]$ does carry a canonical  $n$-shifted symplectic structure (see \cite{ptvv}). Moreover, 
a closed $1$-form of degree  $n$ on $X$ defines a section $X \longrightarrow T^{*}X[n]$
which comes equiped with a natural Lagrangian structure. In particular, 
if $f \in H^{n}(X,\OO_{X})$ is a function of degree  $n$ on $X$, its differential $dR(f)$ defines
a morphism $X \longrightarrow T^{*}X[n]$ together with a Lagrangian structure. The derived critical 
locus of the function $f$ is defined to be the intersection of $dR(f)$ with the zero section 
$$\mathbb{R}Crit(f):=X\times_{0,T^{*}X[n],dR(f)}X.$$
By what we have just seen it comes equiped with a canonical $(n-1)$-shifted symplectic structure. 
We note that when $f=0$ then $\mathbb{R}Crit(f)\simeq T^{*}X[n-1]$ with its canonical 
$(n-1)$-shifted symplectic structure. In general, $\mathbb{R}Crit(f)$ is a perturbation 
of $T^{*}X[n-1]$ obtained using the function $f$. 

The case $n=0$ and $X=\Spec\, A$ a smooth affine scheme has the following explicit description.
The derived critical 
locus $\mathbb{R}Crit(f)$ can then be written as $\Spec\, B$, where 
$B$ can be explicitly represented by the following commutative dg-algebra. 
As a commutative graded algebra $B$ is $Sym_{A}(\mathbb{T}_{A/k}[1])$, where the
differential is the contraction with the $1$-form $dR(f)$. For instance
$\pi_{0}(B)$ is the Jacobian ring of $A$ with respect to $f$. 

\subsubsection*{Hamiltonian actions and symplectic reduction} Shifted symplectic 
and Lagrangian structures in degree $0$ and $1$ can also be used  to 
interpret and extend the notion of Hamiltonian action and of symplectic reduction. 
Let $X$ be a smooth scheme over $k$ (assume $k$ is a field, 
but this is not strictly necessary), equipped with a symplectic structure $\omega$.
Let $G$ be a reductive smooth group scheme over $k$ acting on $X$ and preserving the
form $\omega$. We assume that the action is Hamiltonian in the sense that there
is a moment map
$$\phi : X \longrightarrow \frak{g}^{\vee},$$
which is a $G$-equivariant morphism and which is such that $\omega$ produces 
an isomorphism of complexes of vector bundles on $X$
$$\xymatrix{
\OO_{X}\otimes_{k}\frak{g}  \ar[r]^-{a} \ar[d]_-{id} & \mathbb{T}_{X/k} \ar[r]^-{\phi} \ar[d]^-{\Theta_{\omega}} & 
\OO_{X}\otimes_{k} \frak{g}^{\vee} \ar[d]^-{id}\\
\OO_{X}\otimes_{k}\frak{g}  \ar[r]^-{\phi} & 
\mathbb{L}_{X/k} \ar[r]^-{a} & 
\OO_{X}\otimes_{k} \frak{g}^{\vee},}$$
where $a$ is the morphism induced by the infinitesimal action of $G$ on $X$.

As explained in \cite[\S2.2]{cala} (see also \cite{safr}), the moment map $\phi$ induces a morphism of quotient stacks
$$\phi : [X/G] \longrightarrow [\frak{g}^{\vee}/G]\simeq T^{*}BG[1],$$
which comes equipped with a canonical Lagrangian structure, with respect to the
standard $1$-shifted symplectic structure $\eta$ on $T^{*}BG[1]$. In the same way, 
if $\lambda \in \frak{g}^{\vee}$ with stabilizer $G_{\lambda}$, the natural 
inclusion $BG_{\lambda} \hookrightarrow [\frak{g}^{\vee}/G]$ is again 
equipped with a canonical Lagrangian structure (for instance because the 
inclusion of the orbit $G.\lambda \subset \frak{g}^{\vee}$ is a moment
map for $G.\lambda$, see \cite{cala}). As a consequence, 
the derived stack
$$[X\times_{\frak{g}^{\vee}}\{\lambda\} / G_{\lambda}] \simeq [X/G] \times_{[\frak{g}^{\vee}/G]}BG_{\lambda}$$
is a fibered product of morphisms with Lagrangian structures and thus comes equipped with 
a canonical $0$-shifted symplectic structure (see also \cite{pech} for a more direct proof
of this fact). The new feature here is that 
$\lambda$ does not need to be a regular value of the moment map $\phi$ for this to hold. When
$\lambda$ is a regular value the derived stack $[X\times_{\frak{g}^{\vee}}\{\lambda\} / G_{\lambda}]$ 
is a smooth Deligne-Mumford stack and is the standard symplectic reduction 
of $X$ by $G$. For nonregular values $\lambda$ the derived stack $[X\times_{\frak{g}^{\vee}}\{\lambda\} / 
G_{\lambda}]$ is a nonsmooth derived Artin stack. 

The above interpretation of Hamiltonian actions in terms of Lagrangian structures also has
a version in the so-called quasi-Hamiltonian setting (see \cite{qhamilton}), for which $T^{*}BG[1]$ is replaced by 
$\mathcal{L}(BG)\simeq [G/G]$. The derived stack $[G/G]$ carries a $1$-shifted symplectic
structure as soon as a nondegenerate $G$-invariant scalar product has been chosen on 
$\frak{g}$ (this follows from theorem \ref{t5} for $M=S^{1}$). The quasi-Hamiltonian structures 
are $G$-equivariant morphisms $X \longrightarrow G$ for which the induced morphism
$[X/G] \longrightarrow [G/G]$ is equipped with a natural Lagrangian structure. The quasi-Hamiltonian 
reduction can then still be interpreted as a fibered product of Lagrangian morphisms
(see \cite{cala,safr} for more details). 

\subsection{Existence results}

We arrive at the main existence result of shifted symplectic and Lagrangian structures. We will
start by stating the existence of shifted symplectic structures on derived stacks of maps from
an oriented source to a shifted symplectic target. This result can be seen as a finite dimensional
and purely algebraic version of the so-called AKZS formalism of \cite{aksz}. We will then discuss
the possible generalizations and variations in order to include boundary conditions as well as
noncommutative objects. \smallbreak

We let $X\in \dSt_{k}$ be a derived stack, possibly not represented by 
a derived Artin stack. We say that $X$ is \emph{oriented of dimension $d$}
if there exists a morphism of complexes of $k$-modules
$$or : \mathbb{H}(X,\OO_{X}) \longrightarrow k[-d],$$
which makes Poincar\'e duality to hold in the stable $\sz$-category
$L_{\perf}(X)$ of perfect complexes on $X$ (see \cite{ptvv} for details). More precisely, we ask that 
for all perfect complex $E$ on $X$ the complex $\mathbb{H}(X,E)$ is a perfect 
complex of $k$-modules. Moreover, the morphism $or$ is asked to produce a nondegenerate pairing
of degree $-d$
$$\xymatrix{
\mathbb{H}(X,E) \otimes_{k} \mathbb{H}(X,E^{\vee}) \ar[r] & \mathbb{H}(X,E\otimes_{k} E^{\vee}) \ar[r]^-{Tr}
& \mathbb{H}(X,\OO_{X}) \ar[r]^-{or} & k[-d].}$$
The equivalence induced by the above pairing $\mathbb{H}(X,E) \simeq \mathbb{H}(X,E^{\vee})^{\vee}[-d]$
is a version of Poincar\'e duality for perfect complexes on $X$. 

There are several examples of oriented objects $X$, coming from different origin. Here are the three most 
important examples.

\begin{itemize}

\item A smooth and proper scheme $X$ of relative dimension $d$ over $k$, together with a choice of a 
Calabi-Yau structure $\omega_{X/k}\simeq \OO_{X}$ is canonically an oriented 
object of dimension $d$ as above. The canonical orientation is given by Serre duality 
$H^{0}(X,\OO_{X}) \simeq H^{d}(X,\OO_{X})^{\vee}$, the image of $1 \in H^{0}(X,\OO_{X})$ 
provides an orientation $or : \mathbb{H}(X,\OO_{X}) \longrightarrow k[-d]$. 

\item Let $X$ be a smooth and proper scheme of relative dimension $d$ over $k$, 
and let $X_{DR}$ be its relative
de~Rham object over $k$ (see e.g.\ \cite{simp2}). The cohomology $\mathbb{H}(X_{DR},\OO_{X_{DR}})$
simply is alegbraic de~Rham cohomology of $X$ relative to $k$.
The trace morphism of Grothendieck provides a morphism $H^{d}(X,\Omega^{d}_{X/k})=H^{2d}_{DR}(X/k) \longrightarrow
k$, and thus an orientation $or : \mathbb{H}(X_{DR},\OO_{X_{DR}}) \longrightarrow k[-2d]$, making
$X_{DR}$ into an oriented object of dimension $2d$. The fact that this morphism $or$ provides
the required duality in  $L_{\perf}(X_{DR})$ is the usual Poincar\'e duality for 
flat bundles on $X$.

\item Let $M$ be a compact topological manifold with 
an orientation $H^{d}(M,k) \longrightarrow k$. Considered as a constant
derived stack, $M$, becomes an oriented object of dimension $d$. The required duality in 
$L_{\perf}(M)$ is here Poincar\'e duality for finite dimensional local systems on $M$.

\end{itemize}

\begin{theorem}\label{t5}
Let $X$ be a derived stack which is an oriented object of dimension $d$. Let $Y$ be 
a derived Artin stack locally of finite presentation over $k$ and endowed with an  $n$-shifted symplectic structure. 
Then, if the derived mapping stack $\Map(X,Y)$ is a derived Artin stack locally of finite presentation 
over $k$, then it carries a canonical $(n-d)$-shifted symplectic structure. 
\end{theorem}

In a nutshell the proof of the above theorem follows the following lines. There is a diagram
of derived Artin stacks
$$\xymatrix{
\Map(X,Y) & X\times \Map(X,Y) \ar[r]^-{ev} \ar[l]_-{p} & Y,}$$
where $ev$ is the evaluation morphism and $p$ the natural projection. If $\omega$ denotes 
the  $n$-shifted symplectic form on $Y$, the $(n-d)$-shifted symplectic
structure on $\Map(X,Y)$ is morally defined as 
$$\int_{or}ev^{*}(\omega),$$
where the integration is made using the orientation $or$ on $X$. The heart of the argument is to 
make sense rigorously of the formula above, for which we refer to \cite{ptvv}. 

\subsubsection*{Examples of shifted symplectic structures from the theorem \ref{t4}}
The above theorem, combined with the list of examples of oriented objects given before together
with some of the already mentioned examples of shifted symplectic structures (on $BG$, on $\Parf$ \dots)
provide an enormous number of new instances of shifted symplectic derived Artin stacks. 
The most fundamental examples are the following, for which $G$ stands for a reductive smooth group 
scheme over $k$ with a chosen $G$-invariant scalar product on $\frak{g}$.

\begin{enumerate}

\item Let $X$ be a smooth and proper scheme of relative dimension $d$ over $k$ together with 
an isomorphism $\omega_{X/k}\simeq \OO_{X}$. 
Then the derived moduli stack 
of $G$-bundles, $Bun_{G}(X):=\Map(X,BG)$ is equipped with a canonical 
$(2-d)$-shifted symplectic structure. 

\item Let $X$ be a smooth and proper scheme over $k$ of relative dimension $d$. Then the
derived moduli stack of $G$-bundles equipped with flat connexions, $\Loc_{DR}(X,G):=\Map(X_{DR},BG)$,
carries a canonical $(2-2d)$-shifted symplectic structure.

\item Let $M$ be a compact oriented topological manifold of dimension $d$. The
derived moduli stack $Bun_{G}(M):=\Map(M,BG)$, of $G$-local systems on $M$,
is equipped with a canonical $(2-d)$-shifted symplectic structure.

\end{enumerate}

When the
orientation dimension is $d=2$ the resulting shifted symplectic structures of the above three examples
are $0$-shifted. In these case, the smooth part of the moduli stack recovers some of well known 
symplectic structures on the moduli space of bundles on K3 surfaces, on the moduli space
of linear representations of the fundamental group of a compact Riemann surface \dots. 

\begin{remark}
In the examples $1-3$ above we could have replaced $BG$ by $\Parf$ with its canonical 
$2$-shifted symplectic structure. Chosing a faithful linear representation $G \hookrightarrow \mathbb{GL}_{n}$, 
produces a morphism $\rho : BG \longrightarrow \Parf$, where a vector bundle is considered as 
a perfect complex concentrated in degree $0$. The shifted symplectic structures 
on $\Map(X,BG)$ and $\Map(X,\Parf)$ are then compatible with respect to 
the morphism $\rho$, at least if the $G$-invariant scalar product on $\frak{g}$ is chosen
to be the one induced from the trace morphism on $\frak{gl}_{n}$.
\end{remark}

The theorem \ref{t5} possesses several generalizations and modifications, among which the most
important two are described below. The basic principle here is that any general form of Poincar\'e duality
should induce  a nondegenerate pairing of tangent complexes and must be interpreted
as some shifted symplectic or Lagrangian structure. In the two examples above we deal with
Poincar\'e duality with boundary and in the noncommutative setting. 

\subsubsection*{Derived mapping stack with boundary conditions} The theorem \ref{t5} can be extended
to the case where the source $X$ has a boundary as follows (see \cite{cala} for more details). 
We let $Y$ be a derived Artin stack with an $(n+1)$-shifted symplectic structure, and
$f : Z \longrightarrow Y$ a morphism of derived Artin stacks with a Lagrangian structure. On the other
hand we consider a morphism of derived stacks $j : B \longrightarrow X$ as our 
general source. We assume that $j$ is equipped with a \emph{relative orientation
of dimension $d$} which by definition consists of a morphism of complexes of $k$-modules
$$or : \mathbb{H}(X,B,\OO) \longrightarrow k[-d],$$
where $\mathbb{H}(X,B,\OO)$ is the relative cohomology of the pair $(X,B)$ defined
as the fiber of $\mathbb{H}(X,\OO_{X}) \longrightarrow \mathbb{H}(B,\OO_{B})$. The orientation
$or$ is moreover assumed to be nondegenerate in the following sense. For $E$ a perfect
complex on $X$, we denote by $\mathbb{H}(X,B,E)$ the relative cohomology of the pair
$(X,B)$ with coefficients in $E$, defined as the fiber of $\mathbb{H}(X,E) \longrightarrow \mathbb{H}(B,E)$.
The trace morphism $\mathbb{H}(X,E) \otimes \mathbb{H}(X,E^{\vee}) \longrightarrow \mathbb{H}(X,\OO_{X})$,
together with the orientation $or$
defines a canonical morphism
$$\mathbb{H}(X,E) \otimes \mathbb{H}(X,B,E^{\vee}) \longrightarrow \mathbb{H}(X,B,\OO) \longrightarrow 
k[-d-1],$$
and we ask the induced morphism
$\mathbb{H}(X,E) \longrightarrow \mathbb{H}(X,B,E^{\vee})^{\vee}[-d]$
to be an equivalence. We also ask that the induced morphism
$$\mathbb{H}(B,\OO_{B})[-1] \longrightarrow \mathbb{H}(X,B,\OO) \longrightarrow k[-d]$$
defines an orientation of dimension $(d-1)$ on $B$. 
This is the form of  relative Poincar\'e duality for the
pair $(X,B)$ with coefficients in perfect complexes. When $B=\emptyset$ we recover the
notion of orientation on $X$ already discussed for the theorem \ref{t5}. 

We denote by $\Map(j,f)$ the derived stacks of maps from
the diagram $f : Y \longrightarrow Z$ to the diagram $j : B \longrightarrow X$, 
which can also be written as a fibered product 
$$\Map(j,f) \simeq \Map(B,Y) \times_{\Map(B,Z)}\Map(X,Z).$$
The generalization of theorem \ref{t5}, under the suitable finiteness conditions
on $B$ and  $X$ is the existence of canonical $(n-d+1)$-shifted symplectic structure on 
$\Map(j,f)$ as well as a Lagrangian structure on the morphism
$$\Map(X,Y) \longrightarrow \Map(j,f).$$
The theorem \ref{t5} is recovered when $B=\emptyset$ and $Y=*$ (see remark \ref{r3}).
When $Y=*$ but $B$ is not empty, the statement is that the restriction morphism
$\Map(X,Y) \longrightarrow \Map(B,Y)$
is equipped with a Lagrangian structure (with respect to the $(n-d+1)$-shifted symplectic
structure on $\Map(B,Y)$ given by theorem \ref{t5}). Another consequence is the
existence of compositions of Lagrangian correspondences in the derived setting (see \cite{cala}). 

There are many examples of $j : B \longrightarrow X$ with relative orientations. 
First of all $X$ can be the derived stack obtained from an actual $d$-dimensional 
oriented and compact topological
manifold with boundary $B$, for which the orientation $\mathbb{H}(X,B) \longrightarrow k[-d]$
is given by the integration along the fundamental class in relative homology $[X] \in H_{d}(X,B)$. 
Another example comes from anti-canonical sections. Let $X$ be a smooth and 
projective scheme of relative dimension $d$ over $k$ and $B \hookrightarrow X$ be the derived scheme
of zeros of a section $s\in \mathbb{H}(X,\omega_{X/k}^{-1})$ of the anti-canonical sheaf.
Then the inclusion $j : B \longrightarrow X$ carries a canonical 
relative orientation of dimension $d$ obtained as follows. There is an exact triangle
of quasi-coherent complexes on $X$, $\xymatrix{\omega_{X/k} \ar[r]^-{s} & \OO_{X} \ar[r] & \OO_{B}}$,
giving rise to an exact triangle on cohomologies
$$\xymatrix{
\mathbb{H}(X,\omega_{X/k}) \ar[r] & \mathbb{H}(X,\OO_{X}) \ar[r] & \mathbb{H}(X,\OO_{B})\simeq \mathbb{H}(B,
\OO_{B}),}$$
which identifies $\mathbb{H}(X,B,\OO)$ with $\mathbb{H}(X,\omega_{X/k})$. Grothendieck trace map
furnishes a morphism $or : \mathbb{H}(X,\omega_{X/k}) \longrightarrow k[-d]$
which is a relative orientation of dimension $d$ for the morphism $B \hookrightarrow X$.
An intersecting comment here is that $B$ does not need to be smooth over $k$, and could
be a derived scheme (when $s=0$) or a nonreduced scheme. 

\textbf{Non-commutative spaces} There are noncommutative versions of the theorem \ref{t5}
concerning the existence of shifted symplectic structures on the derived stack 
$\mathcal{M}_{T}$ of objects in a given dg-category $T$ as introduced in \cite{toenvaqu}. Let 
$T$ be a smooth and proper dg-category of $k$ (see \cite{kell,toenvaqu} 
for the definition). There exists a derived stack 
$\mathcal{M}_{T}\in \dSt_{k}$ whose points over a derived $k$-algebra $A$ 
is the classifying space of perfect $T^{op}\otimes_{k}A$-dg-modules. 
The derived stack $\mathcal{M}_{T}$ is not quite a derived Artin stack, but 
is \emph{locally geometric} in the sense that it is a countable union 
of open substacks which are derived Artin stacks (the point here is that 
these open substacks are derived m-Artin stacks but the integer $m$ is not 
bounded and varies with the open substack considered). 

We assume that $T$ comes equipped with an orientation of dimension $d$, by which we mean 
a morphism
$$or : \HH(T) \longrightarrow k[-d],$$
where $\HH(T)$ is the complex of Hochschild homology of $T$ (see \cite{kell}). The morphism
is assumed to be a morphism of mixed complexes, for the natural mixed structure on $\HH(T)$
(see \cite{kell}), and the trivial mixed structure on $k[-d]$. We also assume that $or$ is nondegenerate
in the sense that for all pair of objects $(a,b)$ in $T$, the composite
$T(a,b) \otimes T(b,a) \longrightarrow \HH(T) \longrightarrow k[-d]$
induces an equivalence $T(a,b) \simeq T(b,a)^{\vee}[-d]$. It can be proved that
there exists a canonical $(2-d)$-shifted symplectic structure on $\mathcal{M}_{T}$. 
This statement is a noncoommutative analogue of theorem \ref{t5} as $\mathcal{M}_{T}$ should be
thought as the noncommutative derived mapping stack from the noncommutative space $T$ 
to $\Parf$ (according to the general philosophy of Kontsevich and al. that noncommutative spaces
are dg-categories). The proof of this noncommutative version is very close to the
proof of theorem \ref{t5} exposed in \cite{ptvv}, with suitable modification. It essentially 
consists of defining the shifted symplectic structure by means of the Chern character
of the universal object $\mathcal{E}$ on $\mathbb{M}_{T}$. The heart of the proof 
relies on the correct definition of this Chern character as a mixture of the Chern 
character in noncommutative geometry and the Chern character in derived algebraic geometry
in the style of \cite{toenvezz6}.

Finally, theorem \ref{t5} also possesses a noncommutative version with boundary conditions as follows.
We let $f : T \longrightarrow T_{0}$ be a dg-functor between nice enough dg-categories over $k$. 
We assume given a relative orientation of dimension $d$ on $f$, that is a morphism of mixed
complexes
$$or : \HH(f) \longrightarrow k[-d]$$
where $\HH(f)$ is defined as the homotopy fiber of $\HH(T) \longrightarrow \HH(T_0)$. The orientation
$or$ is also assumed to satisfy nondegeneracy conditions similar to a relative
orientation between derived stacks mentioned before. It can then be proved that the natural morphism
$\mathcal{M}_{T} \longrightarrow \mathcal{M}_{T_{0}}$
comes equipped with a natural Lagrangian structure. 

An important example is the following. We assume that $T$ is equipped with an
orientation of dimension $d$, so $\mathcal{M}_{T}$ carries a natural $(2-d)$-shifted symplectic
structure. We let $\mathcal{M}_{T}^{(1)}$ be the derived stack of morphisms in $T$, which 
comes equipped with two morphisms
$$\xymatrix{\mathcal{M}_{T} \times \mathcal{M}_{T} & \mathcal{M}_{T}^{(1)} \ar[l]_-{s,c} \ar[r]^-{t} & 
\mathcal{M}_{T},}$$
where $t$ sends a morphisms in $T$ to its cone, $s$ sends it to its source and $c$ to its target. 
The correspondence $\mathcal{M}_{T}$ can be seen to carry a canonical Lagrangian
structure with respect to the $(2-d)$-shifted symplectic structure on $\mathcal{M}_{T}^{3}$. 
This Lagrangian structure is itself induced by a natural relative orientation on the dg-functor
$$(s,t,c) : T^{(1)} \longrightarrow T\times T \times T,$$
where $T^{(1)}$ is the dg-category of morphisms in $T$. The relevance of this example
comes from the fact that the correspondence $\mathcal{M}_{T}^{(1)}$ induces
the multiplication on the so-called Hall algebra of $T$ (see \cite{kell} for a review). 
What we are claiming here is that
under the assumption that $T$ comes equipped with an orientation of dimension $d$, 
$\mathcal{M}_{T}$ becomes a monoid in the $\sz$-category of symplectic correspondences in the 
sense of \cite{cala},
which can probably also be stated by saying that $\mathcal{M}_{T}$ is 
a \emph{symplectic 2-Segal space} in the spirit of the higher Segal 
spaces of \cite{dyckkapr}. The compatibility of the
Hall algebra multiplication and the shifted symplectic structure on $\mathcal{M}_{T}$
surely is an important phenomenon and will be studied in a different work.

\subsection{Polyvectors and shifted Poisson structures}

We finish this part by mentioning few words concerning the notion of shifted Poisson structures, dual 
to that of shifted closed $2$-forms, but which is at the moment still under investigation. We present
some of the ideas reflecting the present knowledge. \smallbreak

We let $X$ be a derived Artin stack locally of finite presentation over $k$. For an integer
$n\in \mathbb{Z}$, the \emph{complex of  $n$-shifted polyvector fields on $X$} (relative to $k$)
is the graded complex (i.e. a $\mathbb{Z}$-graded object inside the category of complexes
of $k$-modules)
$$\Pol(X,n)=\displaystyle{\bigoplus_{k\in \mathbb{Z}}}\Pol(X,n)(k) := \displaystyle{\bigoplus_{k\in 
\mathbb{Z}}}\mathbb{H}(X,Sym_{\OO_{X}}^{k}(\mathbb{T}_{X/k}[-1-n])),$$
where $\mathbb{T}_{X/k}=\mathbb{L}_{X/k}^{\vee}$ is the tangent complex of $X$ relative to $k$. 
By definition a \emph{bi-vector $p$ of degree  $n$ on $X$}
is an element  $p\in H^{-n-2}(X,\Pol(X,n)(2))$, or equivalently
a morphism in $L_{\qcoh}(X)$
$$p : \OO_{X} \longrightarrow \phi_{n}^{(2)}(\mathbb{T}_{X/k})[-n],$$
where the symbol $\phi_{n}^{(2)}$ either means $\wedge_{\OO_{X}}^{2}$ if  $n$ is even
or $Sym^{2}_{\OO_{X}}$ is  $n$ is odd. 

When $X$ is a smooth scheme over $k$ and $n=0$, a bi-vector in the sense above 
simply is a section $p \in \Gamma(X,\wedge^{2}_{\OO_{X}}\mathbb{T}_{X/k})$, recovering the usual 
definition. In general, if $\omega \in H^{n}(\mathcal{A}^{2})$ is a $2$-form
of degree  $n$ on $X$, and if $\omega$ is nondegenerate, then we obtain
a bi-vector $p(\omega)$ of degree  $n$ by duality as follows. We represent
the form $\omega$ as a morphism $\mathbb{T}_{X/k} \wedge \mathbb{T}_{X/k} \longrightarrow
\OO_{X}[n]$, and we transport this morphism via the equivalence
$\Theta_{\omega} : \mathbb{T}_{X/k} \simeq \mathbb{L}_{X/k}[n]$ in order to get 
another morphism
$$(\mathbb{L}_{X/k}[n])\wedge (\mathbb{L}_{X/k}[n]) \simeq \phi_{n}^{(2)}(\mathbb{L}_{X/k})[2n] \longrightarrow \OO_{X}[n],$$
which by duality provides $p \in H^{-n}(X,\phi_{n}^{(2)}(\mathbb{T}_{X/k}))$. \smallbreak

The complexes of forms have been shown to carry an important extra structure, namely the de~Rham 
differential. In the same way, the complex $\Pol(X,n)$ does carry an extra structure dual 
to the de~Rham differential: the so-called Schouten bracket. Its definition is much harder
than the de~Rham differential, at least for derived Artin stacks which are not 
Deligne-Mumford, because polyvectors, contrary to forms (see proposition \ref{p3}), do not satisfy some form
of smooth descent (there is not even a well defined pull-back of polyvectors along a smooth morphism).
The theory of polyvector possesses a much more global nature than the theory of forms, and
at the moment there are no simple construction of the Lie bracket on $\Pol(X,n)$ 
except when $X$ is Deligne-Mumford. 

If $X=\Spec\, A$ is an affine derived scheme then $\Pol(X,n)$ has the following
explicit description. We consider $N(A)$ the normalized commutative dg-algebra
associated to $A$, and let $A'$ be a cofibrant model for $N(A)$ as a 
commutative dg-algebra over $k$. The $A'$-module of derivations from $A'$ to itself
is $\mathbb{T}_{A'/k}=\underline{Hom}(\Omega_{A'/k}^{1},A')$. This $A'$-dg-module
is endowed with the standard Lie bracket obtained by taking the (graded) 
commutator of derivations. This Lie bracket satisfies the standard Libniz rule
with respect to the $A'$-dg-module structure on $\mathbb{T}_{A'/k}$, making
$\mathbb{T}_{A'}$ into a dg-lie algebroid over $A'$
(see \cite{vezz2}). The Lie bracket 
on $\mathbb{T}_{A'}$ extends uniquely to the symmetric algebra
$$\Pol(X,n)\simeq Sym_{A'}(\mathbb{T}_{A'/k}[-1-n])$$
as a Lie bracket of cohomological degree $-1-n$ which is compatible with the multiplicative
structure. This results into a graded $p_{n+2}$-algebra structure on 
$\Pol(X,n)$, also called an $(n+2)$-algebra structure, or $(n+2)$-brace algebra structure (see e.g.\ \cite{kont3,tama}).
We will be mainly interested in a part of this structure, namely the 
structure of graded dg-lie algebra on the complex $\Pol(X,n)[n+1]$
(the \emph{graded} nature is unconventional here, as the bracket has itself
a degree $-1$: the bracket of two elements of weights $p$ and $q$ is an element
of weight $p+q-1$. In other words, the Lie operad is also considered
as an operad in graded complexes in a nontrivial manner).

This local picture can be easily globalized for the \'etale topology: when $X$ is
a derived Deligne-Mumford stack, there is a natural graded dg-lie structure
on $\Pol(X,n)[n+1]$. This is a general fact, for any nice enough derived Artin stack, due to the following 
result.

\begin{prop}\label{p5}
Let $X$ be a derived Artin stack $X$ locally of finite presentation over $k$, and $n\in \mathbb{Z}$.
We assume that $X$ is of the form $[Y/G]$ for $Y$ a quasi-projective derived scheme
and $G$ a reductive smooth group scheme acting on $Y$. Then
the
graded complex $\Pol(X,n)[n+1]$ carries a structure of a graded dg-lie algebra.
\end{prop}

At the moment the only proof of this result uses a rather involved construction 
based on $\infty$-operads (see \cite{toen5}). Also, the precise comparison between the 
graded dg-lie structure obtained in the proposition and the more explicit 
construction when $X$ is Deligne-Mumford has not been fully established yet. Finally,
it is believed that the proposition above remains correct in general, without 
the strong condition on $X$ of being a quotient of a quasi-projective derived scheme by a reductive
group. The situation with dg-lie structure on polyvector fields is thus at the moment 
not completely satisfactory. \smallbreak

The graded dg-lie structure on $\Pol(X,n)$ is a crucial piece of data for the definition of
shifted Poisson structures. 

\begin{df}\label{d6}
Let $X$ be a derived Artin stack as in proposition \ref{p5} and $n\in \mathbb{Z}$. The space of \emph{n-shifted
Poisson structure on $X$} is the simplicial set $\mathcal{P}ois(X,n)$ defined by 
$$\mathcal{P}ois(X,n):=Map_{dg\mhyphen lie^{gr}}(k[-1](2),\Pol(X,n)[n+1]),$$
where $dg\mhyphen lie^{gr}$ denotes the $\infty$-category of graded dg-lie algebras over $k$, and
$k[-1](2) \in dg\mhyphen lie^{gr}$ is the object $k$ concentrated in cohomological degree $1$, with
trivial bracket and pure of weight $2$. 
\end{df}

When $X$ is a smooth scheme over $k$, then the space $\mathcal{P}ois(X,0)$
can be seen to be  discrete and equivalent to the set of Poisson structures
on $X$ (relative to $k$) in the usual sense. Another easy case is for $X=BG$, 
for $G$ reductive, 
as $\Pol(X,n)\simeq Sym_{k}(\frak{g}[-n])^{G}$ with $\frak{g}$ being of weight $1$, 
and weight considerations show that the the graded dg-lie algebra
$\Pol(X,n)[n+1]$ must be abelian in this case. In particular, $\Pol(X,n)[n+1]$ is formal
as a graded dg-lie algebra. In particular, $BG$ admits nonzero
$n$-shifted Poisson structures only when $n=2$. When $n=2$ we have moreover
$$\pi_{0}(\mathcal{P}ois(X,2))\simeq Sym_{k}^{2}(\frak{g})^{G}.$$
 
We know little general constructions methods for $n$-shifted Poisson structures. 
It is believed that the main existence statement for $n$-shifted symplectic 
structures (see theorem \ref{t5}) has a version for $n$-shifted Poisson structures too. 
Results in that direction, but only at the formal completion of the constant map,
are given in \cite{john}. It is also believed that the dual of an $n$-shifted symplectic
structure defines a canonical  $n$-Poisson structures. Though this is 
clear at the level of forms and bi-vectors, taking into account the property of
being closed runs into several technical difficulties. Some nonfunctorial 
construction can be done locally, for instance using the Darboux theorem 
for shifted symplectic structure of \cite{bravjoyc,bogr}, but this approach 
has probably no hope to extend to more general derived Artin stacks. 
We thus leave the comparison between $n$-shifted symplectic and Poisson structures as
open questions (we strongly believe that the answers to both questions are positive).

\begin{quest}\label{conj}
Let $X$ be a derived Artin stack as in the proposition \ref{p5} so that 
the space $\mathcal{P}ois(X,n)$ is defined. Let $\mathcal{S}ymp(X,n)$
be the space of  $n$-shifted symplectic structures on $X$
(defined as a full subspace of the Dold-Kan construction
applied to $\mathcal{A}^{2,cl}(X)[n]$).

\begin{itemize}
\item
Can we define a morphism of spaces
$$\mathcal{S}ymp(X,n) \longrightarrow \mathcal{P}ois(X,n)$$
that extends the duality between  $n$-shifted nondegenerate $2$-forms and 
n-shifted bi-vectors ? 

\item Is this morphism inducing an equivalence between $\mathcal{S}ymp(X,n)$ and
the full subspace of $\mathcal{P}ois(X,n)$ consisting of nondegenerate
$n$-shifted Poisson structures ?

\end{itemize}
\end{quest} 

The general theory of  $n$-shifted Poisson structures has not been developed much 
and remains to be systematically studied. There is for instance no clear notion at 
the moment of co-isotropic structures, as well as no clear relations between 
 $n$-Poisson structures and  $n$-shifted symplectic groupoids. The theory is
thus missing some very fundamental notions, one major reason is the inherent 
complexity of the very definition of the Lie bracket on polyvector fields 
of proposition \ref{p5} making all local coordinate type argument useless. \smallbreak

To finish this section, we would like to mention the next step 
in the general theory of shifted Poisson structures. It is known
by \cite{kont2} that a smooth Poisson algebraic $k$-variety $X$  admits a canonical 
quantization by deformation, which is a formal deformation
of the category $QCoh(X)$ of quasi-coherent sheaves on $X$.
In the same way, a derived Artin stack $X$ (nice enough) endowed with 
an $n$-shifted Poisson structure should be quantified by deformations as follows
(the reader will find more details about deformation quantization in the derived setting
in \cite{toen6}).

An  $n$-Poisson structure $p$ on $X$ is by definition a morphism of 
graded dg-lie algebras, and thus of dg-lie algebras
$$p : k[-1] \longrightarrow \Pol(X,n)[n+1].$$
We put ourselves in the setting of derived deformation theory of \cite{luri} (see also 
\cite{hini}). The dg-lie algebra $k[-1]$ is the tangent Lie algebra of the formal line
$\mathbf{Spf}\, k[[t]]$, and the morphism $p$ therefore represents
an element $p \in F_{\Pol(X,n)[n+1]}(k[[t]])$, where
$F_{\frak{g}}$ denotes the formal moduli problem associated to the
dg-lie algebra $\frak{g}$. 

We invoke here the higher formality conjecture, which is a today a theorem in many 
(but not all) cases. 

\begin{conj}\label{formality}
Let $X$ be a nice enough derived Artin stack over $k$ and $n\geq 0$. Then 
the dg-lie algebra $\Pol(X,n)[n+1]$ is quasi-isomorphic
to the dg-lie algebra $\HH^{E_{n+1}}(X)[n+1]$, where 
$\HH^{E_{n+1}}$ stands for the iterated Hochschild 
cohomology of $X$.
\end{conj}

Somme comments about conjecture \ref{formality}.

\begin{itemize}

\item The higher Hochschild cohomology $\HH^{E_{n+1}}(X)$ is
a global counter-part of the higher Hochschild cohomology of \cite{pira}. 
As a complex it is defined to be 
$$\HH^{E_{n+1}}(X):=End_{L_{\qcoh}(\mathcal{L}^{(n)}X)}(\OO_{X}),$$
where $\mathcal{L}^{(n)}(X):=\Map(S^{n},X)$ is the higher 
dimension free loop space of $X$ (suitably completed
when  $n$ is small). The sheaf $\OO_{X}$ is considered
on $\mathcal{L}^{(n)}(X)$ via the natural morphism $X \longrightarrow \mathcal{L}^{(n)}(X)$
corresponding to constant maps. It is proven in \cite{toen5} that 
$\HH^{E_{n+1}}(X)$ endows a natural structure of a $E_{n+2}$-algebra, and thus
that $\HH^{E_{n+1}}(X)[n+1]$ has a natural dg-lie algebra structure, 
at last when $X$ is a quotient stack of a derived quasi-projective scheme 
by a linear group (see \cite{fran} for the special case of affine derived schemes).

\item The conjecture \ref{formality} follows from the main result of \cite{toen5} when $n>0$
(and $X$ satisfies enough finiteness conditions). When $n=0$ and $X$ is a smooth
Deligne-Mumford stack the conjecture is a consequence of Kontsevich's formality theorem (see \cite{kont2}). 
When $X$ is a derived Deligne-Mumford stack Tamarkin's proof of Kontsevich's formality 
seems to extend to also provide a positive answer to the conjecture (this
is already observed implicitly in \cite{kont3}, but has also been explained
to me by Calaque). 
Finally, for $n=0$ and $X$ is a derived Artin stack which is not Deligne-Mumford,
the conjecture is wide open. 

\item The case $n=1$ and $X$ a smooth scheme of the conjecture 
appears in \cite{kapu}. The case where $X$ is a smooth affine variety appears implicitly in 
\cite{kont3} and has been known from experts at least for the case of polynomial algebras. We also 
refer to \cite{calawill} for related results in the context of commutative dg-algebras.   

\end{itemize}

The consequence of the conjecture \ref{formality} is the existence of 
deformations quantization of shifted Poisson structures. Indeed, 
let $n\geq 0$ and $p$ an  $n$-shifted Poisson structure on $X$. 
By conjecture \ref{formality} we get out of $p$ a morphism of dg-lie algebras
$p : k[-1] \longrightarrow \HH^{E_{n+1}}(X)[n+1],$
and thus an element 
$$p \in F_{\HH^{E_{n+1}}(X)[n+1]}(k[[t]]).$$
When $n<0$ the same argument provides an element 
$$p \in  F_{\HH^{E_{-n+1}}(X)[-n+1]}(k[[t_{2n}]]),$$
where $t_{2n}$ is now a formal variable of cohomological degree $2n$. 

The element $p$ as above defines quantizations by deformations thanks to the
following theorem whose proof will appear elsewhere 
(we refer to \cite{fran} for an incarnation of this result in the topological context).

\begin{theorem}\label{t6}
For $n\geq 0$, the formal moduli problem $F_{\HH^{E_{n+1}}(X)[n+1]}$ 
associated dg-lie algebra $\HH^{E_{n+1}}(X)[n+1]$, controls formal
deformations of the $\sz$-category $L_{\qcoh}(X)$ considered as an
$n$-fold monoidal stable $k$-linear $\sz$-category.
\end{theorem}

The element $p$ defined above, and the theorem \ref{t5}, defines
a formal deformation $L_{\qcoh}(X,p)$ of $L_{\qcoh}(X)$ as an $|n|$-fold monoidal
$\sz$-category, which by definition is the deformation quantization
of the  $n$-shifted Poisson structure $p$. 

\begin{remark}\label{rfinal}
\begin{itemize}

\item Theorem \ref{t6} refers to a rather evolved notion of deformation of $n$-fold monoidal
linear $\sz$-categories, based on a higher notion of Morita equivalences. In the affine case, 
this is incarnated by the fact that $E_{n+1}$-algebras must be considered as
$(n+1)$-categories with a unique object (see e.g.\ \cite{fran}). The precise notions and
definitions behind theorem \ref{t6} are out of the scope of this survey.

\item When $X=BG$, for $G$ reductive, and the $2$-Poisson structure
on $X$ is given by the choice of an element $p\in Sym^{2}(\frak{g})^{G}$, 
the deformation quantization $L_{\qcoh}(X,p)$ is a formal 
deformation of the derived $\sz$-category of representations of $G$
as a $2$-fold monoidal $\sz$-category. This deformation is the quantum
group associated to $G$ and the choice of $p \in Sym^{2}(\frak{g})^{G}$. 

\item The derived mapping stacks $X=\Map(M,BG)$ are often 
endowed with  $n$-shifted symplectic structures (see theorem \ref{t5}). By 
\ref{conj} these are expected to correspond to  $n$-shifted Poisson structures on $X$, 
which can be quantified by deformations as explained above. 
These quantization are very closely related to 
quantum invariants of $M$ when $M$ is of dimension $3$ (e.g.
Casson, or Donaldons-Thomas invariants). In higher dimension
the quantization remains more mysterious and will be studied in 
forthcoming works.

\end{itemize}
\end{remark}

%\section{References}

\iffalse

\cite{And}, \cite{FrQu}, \cite{Jon}.

\enddocument